\documentclass[graybox]{svmult}

\usepackage{mathptmx}       
\usepackage{helvet}         
\usepackage{courier}        
\usepackage{type1cm}        
\usepackage{makeidx}         
\usepackage{graphicx}        

\usepackage{multicol}        
\usepackage[bottom]{footmisc}

\usepackage{color}
\usepackage[cmyk,usenames,dvipsnames]{xcolor}
\usepackage{latexsym }
\usepackage{amsmath, amssymb,amsfonts,amscd}
\usepackage{graphicx,hyperref,enumerate,verbatim}
\usepackage{epsfig}
\usepackage{pdfsync}
\usepackage[export]{adjustbox}
\usepackage{cite}
\usepackage{listings}
\usepackage{float}

\definecolor{purple}{rgb}{0.65, 0, 1}
\definecolor{orange}{rgb}{1,.5,0}
\definecolor{green}{rgb}{0, 1, 0}
\definecolor{orange}{rgb}{1,.5,0}
\definecolor{gray}{rgb}{.6,.6,.6}

\newcommand{\R}{{\mathord{\mathbb R}}}
\newcommand{\NN}{{\mathord{\mathbb N}}}
\newcommand{\Rd}{{\mathord{\mathbb R}^d}}
\newcommand{\grad}{\nabla}

\newcommand{\W}{\mathcal{W}}
\newcommand{\argmin}{\operatornamewithlimits{argmin}}
\def\e{\varepsilon}
\def\epsilon{\varepsilon}
\def\P{\mathcal{P}}
\def\S{\mathcal{S}}
\def\E{E}

\makeindex             

\begin{document}

\title*{Aggregation-diffusion equations: dynamics, asymptotics, and singular limits}
\author{Jos\'e A. Carrillo, Katy Craig, and Yao Yao}
\institute{Jos\'e A. Carrillo \at Imperial College London\\
 \email{carrillo@imperial.ac.uk}
\and 
Katy Craig \at University of California, Santa Barbara\\
\email{kcraig@math.ucsb.edu}
\and
Yao Yao \at Georgia Institute of Technology\\
 \email{yaoyao@math.gatech.edu}}

\maketitle

\abstract{Given a large ensemble of interacting particles, driven by nonlocal interactions and localized repulsion,   the mean-field limit leads to a class of nonlocal, nonlinear partial differential equations known as  \emph{aggregation-diffusion equations}. Over the past fifteen years, aggregation-diffusion equations have become widespread in biological applications and have also attracted significant mathematical interest, due to their competing forces at different length scales. These competing forces  lead to rich dynamics, including symmetrization, stabilization, and metastability, as well as sharp dichotomies separating well-posedness from finite time blowup.  In the present work, we review known analytical results for aggregation-diffusion equations and consider singular limits of these equations, including the slow diffusion limit, which leads to the constrained aggregation equation, and localized aggregation and vanishing diffusion limits, which lead to metastability behavior. We also review the range of numerical methods available for simulating solutions, with special attention devoted to recent advances in deterministic particle methods. We close by applying such a method -- the \emph{blob method for diffusion} -- to showcase key properties of the dynamics of aggregation-diffusion equations and related singular limits. }



\section{Introduction}

Many phenomena in the life sciences, ranging from the microscopic to  macroscopic level, exhibit surprisingly similar structures. Behavior at the microscopic level, including ion channel transport, chemotaxis, and angiogenesis,  and behavior at the macroscopic level, including   herding of animal populations, motion of human crowds,  and bacteria orientation, are both largely driven by long-range attractive forces, due to electrical, chemical or social interactions, and short-range repulsion, due to dissipation or finite size effects.

Various modeling approaches at the agent-based level, from  cellular automata to Brownian particles, have been used to describe these phenomena. To pass from microscopic models \cite{DCBC,CS07,CDFLS} to continuum descriptions requires analysis of the mean-field limit,  as the number of agents becomes large \cite{BCC11,J,CCHmean,CCHS,JWz}. This  approach leads to a continuum kinematic equation for the evolution of the density of individuals $\rho(x,t)$ known as the \emph{aggregation equation},
\begin{align} \label{aggeqn0}
\frac{\partial\rho}{\partial t}+\nabla \cdot \left(\rho u\right)=0 ,
\quad   u=-\nabla W*\rho\,,
\end{align}
where $W:\R^d \longrightarrow \R$ is a symmetric  
 interaction potential. Typical examples of interaction potentials $W$ are given by
repulsive-attractive power laws
\begin{align}\label{repattpower} W(x) = \frac{|x|^A}{A} -  \frac{|x|^B}{B} , \quad 2\geq A>B>-d  ,
\end{align}
or Morse potentials \cite{DCBC},
$$
W(x) = -C_A e^{-|x| /\ell_A} + C_R e^{-|x| /\ell_R} , \quad C=C_R/C_A >1, \ 
\ell=\ell_R/\ell_A<1 , \ C\ell^d<1.
$$

The continuum equation (\ref{aggeqn0}) may also be modified to include  linear or nonlinear diffusion terms, which arise from two possible modeling assumptions: noise at the level of interacting particles (linear
diffusion) or scaling hypotheses on a repulsion potential that  depend on the inter-particle distance (nonlinear diffusion)
\cite{BodnarVelazquez,O90,TopazBertozzi2}. The latter case can be seen via the following formal argument: take $W_\nu=W+2\nu \delta_0$ with $W$ an attractive potential. In this case, $W_\nu$ induces strongly localized repulsion of strength $\nu$, along with a nonlocal attractive force via the potential $W$. Formally, the corresponding PDE is given by
\begin{equation}\label{aux}
\rho_t  = \nu\Delta \rho^2 + \nabla\cdot(\rho\nabla(\rho*W)) \,.
\end{equation} 
This equation can be obtained from particle approximations \cite{O90,BodnarVelazquez,CCS} in a limit in which the potential has a repulsive part that becomes concentrated at the origin as the number of particles increases, with an overall mean-field scaling for the forces. See \cite{carrillo2017blob} and the references therein for the rigorous relation between \eqref{aggeqn0} and \eqref{aux}.

Incorporating either linear or nonlinear diffusion into the aggregation equation (\ref{aggeqn0}), we arrive at a partial differential equation of the form
\begin{equation*}
\frac{\partial \rho}{\partial t} = -\nabla \cdot (\rho u) =
\nabla \cdot \left [ \rho \nabla \left ( U' \left ( \rho \right )
+ W\ast \rho \right) \right ],
\end{equation*}
where the function $U(\rho)$ determines the type of diffusion. In particular, $U(s)=s\log s$ corresponds to linear
diffusion and $U(s)={s^m}/{(m-1)}$, $m>0$ corresponds to nonlinear diffusion \cite{Va}.  Nonlocal partial differential equations of this form are becoming widespread in both  modeling and theory, and a key question from the biological perspective is how to identify the proper mechanisms for collective motion among the many potential equations describing the behavior \cite{KCBFL,Carillo09_Kinetic_Attraction-Repulsion,MT}.

In this chapter, we consider a class of partial differential equations with nonlocal interactions and nonlinear diffusion, known as \emph{aggregation-diffusion equations}
\begin{equation}\label{aggregation}
\rho_t  = \Delta \rho^m + \nabla\cdot(\rho\nabla(\rho*W)) , m \geq 1.
\end{equation}
Formally, equations of this form   have a gradient flow structure with respect to the Wasserstein metric $d_2$ on the space of probability measures with finite second moment $\P_2(\Rd)$  \cite{JKO,O2,Villani,AGS,cmcv-03}. In particular, defining the free energy
\begin{equation}\label{energy}
\E[\rho] = \underbrace{ \frac{1}{m-1}\int_{\mathbb{R}^d}\rho^m dx}_{ \S_m[\rho]} + \underbrace{\frac{1}{2}\int_{\mathbb{R}^d} \rho(\rho*W)dx}_{\W[\rho]} ,
\end{equation} 
one may rewrite the aggregation diffusion equation \eqref{aggregation}  as
\begin{align} \label{formal W2 grad flow}
 \frac{d \rho}{d t}  =  - \grad_{d_2} E(\rho) ,\qquad \grad_{d_2} E(\rho)  = - \grad \cdot \left(\rho \grad \frac{\delta E}{\delta \rho} \right) \ ,
 \end{align}
 where  $\grad_{d_2}$ denotes the gradient with respect to the Wasserstein metric and $\grad \frac{\delta E}{\delta \rho}$ denotes the first variation of the free energy with respect to $\rho$. The first term in the free energy (\ref{energy}) is the entropy $\S_m[\rho]$, which gives rise to the diffusion term, where we use the convention $\frac{1}{m-1} \rho^{m-1} = \log \rho$ for $m=1$. The second term in the energy (\ref{energy}) is the interaction energy $\W[\rho]$, which gives rise to the aggregation term.  

An archetypical example of aggregation-diffusion equation (\ref{aggregation}) is the Keller-Segel model of chemotaxis in mathematical biology, which describes the collective motion of cells (usually bacteria or slime mold) that are attracted by a self-emitted chemical substance  \cite{ks}. Let $\rho(x,t)$ denote the population density of a cell colony on a two dimensional surface subject to Brownian motion, and assume these cells have an additional drift velocity due to their tendency to move towards higher concentrations of a chemical attractant $c(x,t)$. The cells themselves are continually emitting this chemical attractant, and it diffuses across the surface and decays with rate $\alpha$. This leads to the following system of parabolic equations:
\begin{equation}
\label{ks_system}
\begin{cases}\rho_t = \Delta \rho - \nabla\cdot (\rho \nabla c) ,\\
\epsilon c_t = \Delta c - \alpha c + \rho.
\end{cases}
\end{equation}
In fact, the same equation was also proposed by Patlak, as a mathematical model for random walk with persistence and external bias \cite{patlak}.

Since the chemical attractant reaches its equilibrium much faster than the bacteria density, it is common to take  $\epsilon\to 0$, simplifying the second equation so that \eqref{ks_system} becomes a parabolic-elliptic system. This simplification allows us to write $c = -W * \rho$, where $W$ is the Newtonian potential (if $\alpha=0$) or the Bessel potential (if $\alpha>0$). Finally, substituting $c$ into the first equation gives us a single  equation for $\rho(x,t)$ in the form of \eqref{aggregation} with $m=1$.  Variations of the Keller-Segel model, such as the consideration of volume effects, lead to a range of aggregation-diffusion equations of the form \eqref{aggregation}, see \cite{CC06}. 
 
 We begin in section  \ref{wellposedsec} by reviewing analytical results for aggregation-diffusion equations, focusing on conditions that ensure solutions are globally well-posed or blow-up in finite time. The majority of these analytical results consider interaction kernels of one of the following forms:
 \begin{itemize}
\item Power-law kernel: 
\begin{equation}\label{power_kernel}
W_k(x) = \begin{cases} \frac{|x|^k}{k} & \text{ if }k\neq 0\\ \ln|x| & \text{ if } k=0.\end{cases}
\end{equation}
\item Integrable kernel: $W \in L^1(\mathbb{R}^d)$.
\end{itemize}
In both cases, the associated free energy functional plays an important role in the study of well-posedness of solutions and properties of steady states. Formally taking time derivative along a solution, we deduce
\[
\frac{d}{dt} E[\rho] = -\int_{\mathbb{R}^d} \rho \left|\nabla (\frac{m}{m-1}\rho^{m-1} + \rho*W)\right|^2 dx \leq 0.
\]
This energy dissipation inequality (in the integral form) can be made rigorous for weak solutions and can also be seen as a consequence of the equation's underlying Wasserstein gradient flow structure. We begin, in subsection \ref{subsection:KS}, by considering the classical Keller-Segel equation in two dimensions, providing heuristic arguments illustrating the dichotomy between well-posedness and finite time blowup, as well as summarizing rigorous results. We then discuss the case for more general aggregation diffusion equations in subsection \ref{subsection:wellposed} and consider long time behavior of solutions in subsection \ref{subsection:diffusiondominated}.

In section \ref{singularlimits}, we consider singular limits of aggregation diffusion equations in two limiting regimes. We begin in subsection \ref{slowdiff} by considering the \emph{slow diffusion limit}, as the diffusion exponent $m \to +\infty$, which leads to the constrained aggregation equation. Then, in subsection \ref{singularmeta}, we discuss singular limits that affect the balance between aggregation and diffusion, causing solutions to exhibit metastable behavior. In subsection \ref{metaatt}, we consider the localized attraction limit, and in subsection \ref{metadiffvanish}, we consider the vanishing diffusion limit.

In section \ref{numericsmainsection}, we review recent work on numerical methods for aggregation diffusion equations, including the recent blob method for diffusion, developed by Carrillo, Craig, and Patacchini \cite{carrillo2017blob}. We conclude  in section \ref{numericssection}  by applying this numerical method to illustrate properties of the singular limits discussed in section \ref{singularlimits}. 

 
\section{Well-posedness, steady states, and dynamics}  \label{wellposedsec}

We now give a brief overview of analytical results for aggregation diffusion equations \eqref{aggregation}, describing conditions that ensure well-posedness or finite-time blow-up of solutions , existence or non-existence of steady states, and long time  behavior of solutions. We begin  in section~\ref{subsection:KS} by considering the classical two-dimensional \emph{Keller-Segel equation}, where the interaction kernel   in equation \eqref{aggregation} is given by the Newtonian potential, $W(x) = \frac{1}{2\pi} \ln |x|$. The analysis in this particular case will serve as a guideline to understand the equation's behavior for more general interaction kernels. In section~\ref{subsection:wellposed}, we review results classifying well-posedness and finite blow-up  for general interaction kernels. Finally, in section~\ref{subsection:diffusiondominated}, we discuss the steady states and dynamics of solutions in the diffusion-dominated regime.


\subsection{Keller-Segel equation with linear diffusion in $\mathbb{R}^2$} \label{subsection:KS}
We begin with the classical Keller-Segel equation, 
 \begin{equation}\label{eq:KS}
\rho_t = \Delta \rho + \nabla\cdot\Big(\rho\nabla\Big(\frac{1}{2\pi} \ln|x| * \rho\Big)\Big)\quad \text{ in } \mathbb{R}^2\times [0,T).
\end{equation}
As the primary interest of the present work is aggregation diffusion equations of the form   \eqref{aggregation},  we will not discuss the vast literature on the parabolic-parabolic Keller-Segel equation; see instead \cite{Horstmann,BBTW}.

A fundamental property of the Keller-Segel equation   is that solutions are subject to a remarkable dichotomy depending on their mass $M = \int \rho(x)dx$. In particular, solutions exist globally in time if $M>8\pi$, whereas they blow-up in finite time when $M<8\pi$.  The fact that the critical value of the mass is $8 \pi$ can be seen from the following formal scaling argument, which we will generalize to a range of interaction kernels $W$ in the following section~\ref{subsection:wellposed}.

As described in the introduction, the Keller-Segel equation \eqref{eq:KS} is formally the Wasserstein gradient flow of the energy 
\begin{align} \label{KSenergy}
E[\rho] = \S_m[\rho] + \W[\rho] = \int_{\mathbb{R}^2} \rho \log \rho dx + \frac{1}{4\pi }\int_{\mathbb{R}^2} \rho(\ln|x|*\rho)dx.
\end{align}
In particular, solutions of the Keller-Segel equation are characterized by the fact that they move in the direction of steepest descent of the energy (\ref{KSenergy}).
Given $\rho\in L^1\cap L^\infty(\mathbb{R}^d)$,  consider the following dilation of $\rho$, which concentrates the density while preserving its mass:
\[
\rho_\lambda(x) := \lambda^2 \rho(\lambda x), \quad \lambda \gg 1 .
\] 
Then, $E[\rho_\lambda]$ and $E[\rho]$ are related in the following way:
\begin{equation*}
\begin{split}
E[\rho_\lambda] &= \S_m[\rho_\lambda] + \W[\rho_\lambda] = (\S_m[\rho] + 2 M \log\lambda) + \left(\W[\rho] - \frac{ M^2\log\lambda}{4\pi}\right)\\
&= E[\rho] +  M \log \lambda \left(2-\frac{M}{4\pi}\right).
\end{split}
\end{equation*}
Consequently, 
\begin{align*}
\lim_{\lambda\to\infty} E[\rho_\lambda] = \begin{cases} +\infty &\text{ when } M< 8 \pi , \\  - \infty &\text{ when } M> 8 \pi . \end{cases}
\end{align*}
Thus, when $M<8\pi$, it is not energy favorable for the solution to concentrate to a $\delta$-function, so one  expects global well-posdness. On the other hand, when $M>8\pi$, it is possible for the solution to blow-up in finite time. 

Motivated by this formal argument, we now summarize the rigorous results.

 \vspace{-5mm} \paragraph{\textbf{Subcritical mass, $M<8\pi$}}  \vspace{-3mm}
Global well-posedness for subcritical mass  was obtained by Blanchet, Dolbeault, and Perthame \cite{DP, bdp}, improving an an earlier result of J\"ager and Luckhaus \cite{JL}, which showed global existence for smaller mass. The key idea behind this approach is to use the logarithmic Hardy-Littlewood-Sobolev inequality \cite{CL}
\begin{equation}\label{logHLS}
\int_{\mathbb{R}^2} \rho \log \rho dx + \frac{2}{M} \int_{\mathbb{R}^2} \rho (\rho * \ln |x|) dx \geq -C(M).
\end{equation}
Combining \eqref{logHLS} with the fact that  $M < 8 \pi$ and the energy decreases along solutions, $E[\rho(t)] \leq E[\rho_0]$, one can obtain an a priori upper bound of the entropy $\int \rho \log \rho dx$. This uniform in time equi-integrability allows one to use  iterative arguments to upgrade the bound to a uniform $L^\infty$ bound,  ensuring global well-posedness.

In addition to global well-posedness,  long-time behavior for solutions with subcritical mass has also been well studied. Blanchet, Dolbeault, and Perthame \cite{bdp} showed that the solution will indeed follow the scaling of the heat equation as $t\to\infty$, converging to a self-similar profile in the rescaled variables. Furthermore, the  rate of convergence is exponential \cite{BDEF, CD,FM}.

 \vspace{-5mm} \paragraph{\textbf{Supercritical mass, $M>8\pi$}}  \vspace{-3mm}
 In the case of supercritical mass,  any initial data $\rho_0 \in L^1_+((1+|x|^2)dx)$ will experience finite-time blow-up \cite{DP, bdp}. This was proved by Blanchet, Dolbeault, and Perthame  by tracking the evolution of the second moment $M_2[\rho(t)] := \int \rho(x,t) |x|^2 dx$. During the existence of the solution, one has
\[
\frac{d}{dt} M_2[\rho(t)] = 4M \left(1-\frac{M}{8\pi}\right),
\]
so that the second moment becomes negative in finite time, indicating that the solution must break down by this time. Note that the above virial argument does not give the nature of the blow-up. Subsequent studies showed that solutions must concentrate at least $8\pi$ of their mass into a single point at the blow-up time \cite{SS, BM}, as well as considering finer properties of the blow-up profile \cite{HV, senba,RS}.

 \vspace{-5mm} \paragraph{\textbf{Critical mass, $M=8\pi$}}  \vspace{-3mm}
At the critical mass, the free energy functional $E[\rho]$ is invariant  under dilations, and there exists a one-parameter family of steady states $\bar \rho_\lambda(x) := \frac{M}{\pi} \frac{\lambda}{(\lambda + |x|^2)^2}$, all with infinite second moment, which are the unique global minimizers of the free energy, up to a translation. While solutions of the Keller-Segel equation exist globally in time \cite{BCM},  solutions with finite initial second moment exhibit an infinite time aggregation that is stable under perturbations \cite{GhMas}. On the other hand, if the initial data has infinite second moment and is sufficiently close to a stationary solution $\bar\rho_\lambda$, then it converges to $\bar\rho_\lambda$, with  quantitative rate \cite{BCC, CF}.


\subsection{Well-posedness v.s. blow-up for general interaction kernels} \label{subsection:wellposed}
We now describe how the same dichotomy that separates global well-posedness from finite time blowup for solutions of the classical Keller-Segel equation can be adapted to general aggregation diffusion equations \eqref{aggregation}, in arbitrary dimensions.  We begin by considering the case of attractive power-law kernels $W_k$ of the form \eqref{power_kernel}. Since blow-up  occurs locally near a specific point, the regimes of global well-posedness vs. finite time blowup are primarily determined by the singularity of interaction kernel near the origin. We conclude by discussing results for more general interaction kernels.

As in the Keller-Segel equation case, let us again consider the behavior of the free energy functional $E[\rho]$ under mass-preserving dilation $\rho_\lambda(x) = \lambda^d \rho(\lambda x)$. As before, both the entropy $\S_m[\rho]$  and the interaction energy $\W_k[\rho]$ possess simple homogeneity properties under this transformation:
\begin{equation}\label{scale_entropy}
\S_m[\rho_\lambda] = \begin{cases}
\lambda^{(m-1)d} \S_m[\rho] &\text{if } m> 1,\\
\S_m[\rho] + dM\log\lambda &\text{if } m = 1,
\end{cases}
\end{equation}
and
\begin{equation}\label{scale_interaction}
\W_k[\rho_\lambda] = \begin{cases}
\lambda^{-k} \W_k[\rho] &\text{if } k\neq 0,\\
\W_k[\rho] - M^2 \log\lambda &\text{if } k=0.
\end{cases}
\end{equation}

Comparing the scaling properties of $\S_m$ and $\W_k$ for different values of diffusion exponent $m$ and interaction range $k$ motivates the consideration of three regimes, separated by the critical diffusion exponent $m_c = 1-k/d$. We now summarize what is known in each regime, regarding well-posedness vs. blowup.

 \vspace{-5mm} \paragraph{\textbf{Diffusion dominated regime, $m>m_c$}}  \vspace{-3mm}
In the diffusion dominated regime, Calvez and Carrillo \cite{CC06} and Sugiyama \cite{sugiyama} considered the Newtonian potential in arbitrary dimensions $d \geq 3$, proving global well-posedness and a uniform-in-time $L^\infty$ bound  for  any initial data $\rho_0 \in L^1\cap L^\infty(\mathbb{R}^d)$. Bedrossian, Rodr\'iguez, and Bertozzi generalized this result to  general power-law kernels that are no more singular than Newtonian, $2-d\leq k \leq 0$ \cite{BRB, BR}. For  kernels that are more singular than Newtonian, $-d<k<2-d$, Zhang recently proved that solutions  remain globally bounded if either $k>1-d$  or $m<2$ \cite{zhang}. 

 \vspace{-5mm} \paragraph{\textbf{Aggregation dominated regime, $1 \leq m<m_c$}}  \vspace{-3mm}
In the aggregation-dominated regime, most results consider the case when the interaction potential is Newtonian, $k = 2-d$. Sugiyama proved that  for arbitrarily small mass $M>0$, there exist solutions with mass $M$ that blow up in finite time \cite{sugiyama} .  On the other hand, Bedrossian showed that, if the initial data is sufficiently flat (even if the mass is large),  solutions exist globally, and dissipate with porous medium equation scaling as $t\to\infty$ \cite{bedrossian2}.

Subsequently, Bian and Liu showed that solutions with small $L^p(\Rd)$ norm, $p~=~\frac{d(2-m)}{2}$, exist globally in time, where $p$ is chosen so that diffusion and aggregation are balanced under the scaling that keeps the $L^p$ norm invariant  \cite{BL}. The special case $m=\frac{2d}{2+d}$ was  studied further by Chen, Liu, and Wang, who showed that solutions with small $L^p$ norm exist globally, whereas solutions with large $L^p$ norm must blow-up in finite time  \cite{CLW}. Chen and Wang then generalized this result to all   $\frac{2d}{2+d}<m<2-\frac2d$ \cite{CW2} .  Finally, for general interaction potentials that are no more singular than Newtonian, $2-d\leq k \leq 0$, Bedrossian, Rodr\'iguez, and Bertozzi showed that there exist solutions with arbitrarily small mass that blow up in finite time \cite{BRB}. 

 \vspace{-5mm} \paragraph{\textbf{Fair competition regime, $m=m_c$}}  \vspace{-3mm}
In the fair competition regime, aggregation diffusion equations exhibit a dichotomy in terms of the mass of solutions that is similar to the classical Keller-Segel equation in $\mathbb{R}^2$. In this case, the Hardy-Littlewood-Sobolev inequality with optimal constant $C(k,d)$,
\[
\iint_{\mathbb{R}^d\times \mathbb{R}^d} \rho(x)\rho(y) W_k(x-y) dxdy \leq C(k, d)M^{2-m_c} \|\rho\|_{m_c}^{m_c}
\]
 plays the same role as the logarithmic Hardy-Littlewood-Sobolev inequality \eqref{logHLS} for the Keller-Segel equation. In particular,
combining this inequality with the fact that the free energy functional decreases along solutions gives a critical value of the mass $M_c = M_c(k,d)$, for which one has an a priori bound on $\|\rho\|_{m_c}$ if $M<M_c$, which then leads to uniform in time $L^\infty$ bound. The results on well-posedness v.s. blow-up can be summarized as follows:

\begin{itemize}
\item[-] \textbf{Subcritical Mass, $M < M_c$.} Blanchet, Carrillo, and Lauren{\c{c}}ot proved that, for Newtonian interaction $k = 2-d$, solutions remain bounded globally in time \cite{BCL},  and there exists self-similar solutions decaying in time with the porous medium equation.  Bedrossian, Rodriguez, and Bertozzi generalized this result to power-law kernels no more singular than Newtonian, $2-d\leq k\leq 0$. Furthermore, Bedrossian proved that,  
 under certain conditions of the initial data, solutions decay to the self-similar spreading solutions of the porous medium equation with an explicit convergence rate \cite{bedrossian2}.

\item[-]  \textbf{Supercritical Mass, $M >M_c$.} For kernels that are no more singular than Newtonian, there exist solutions with mass $M$ that blow-up in finite time \cite{BCL, BRB}. Bedrossian and Kim showed that, in the Newtonian case $k=2-d$,   all radial solutions blow up in finite time \cite{BK}. However, it remains  unknown whether all non-radial solutions must also blow up.

\item[-]  \textbf{Critical Mass $M = M_c$} In the Newtonian case, Blanchet, Carrillo and Lauren{\c{c}}ot showed that solutions are globally well-posed, and their $L^\infty$ norm is globally bounded in time \cite{BCL}. Furthermore, there is a family of compactly supported stationary solutions that are dilations of each other and all of  which are global minimizers of the free energy. Calvez, Carrillo, and Hoffmann extended the latter result to general power-law interaction potentials $-d<k<0$ \cite{CalCarrHo, CCH2}.  Finally, in the Newtonian case, Yao showed that every radial solution with compact support will converge to some stationary solution in this family \cite{yao}, though the asymptotic behavior of non-radial solutions remains unclear.

\end{itemize}

 \vspace{-5mm} \paragraph{\textbf{General interaction potentials}}  \vspace{-3mm}
 
We close this section by discussing what is known about well-posedness of aggregation diffusion equations for general interaction potentials $W$, which are not necessarily of power-law form. If the interaction potential is purely attractive, Bedrossian, Rodr\'iguez, and Bertozzi fully develop the well-posedness theory for kernels with singularity up to and including Newtonian \cite{BRB, BR}. On the other hand, if the interaction potential is a repulsive-attractive power-law potential (\ref{repattpower}) for $A \geq 2-d$, Carrillo and Wang prove a global-in-time $L^\infty$ bound for all $m\geq1$, under the assumption that solutions exist locally in time  \cite{CW}.

General well-posedness results were recently obtained by Craig and Topaloglu, without any assumptions on the specific structure of the attractive and repulsive parts of the interaction potential \cite{CraigTopalogluHeight}. In particular, under weak regularity assumptions on the interaction potential, which roughly correspond to being no more singular than the Newtonian potential, they prove that solutions of aggregation diffusion equations exist and provide sufficient conditions for global in time uniqueness. In particular, this extends the existence theory for power-law interaction kernels $W_k(x) = |x|^k/k$ to include unbounded initial data, as long as it belongs to domain of the energy, $\rho_0 \in D(E)$, as well as obtaining existence for  repulsive-attractive power-law potentials (\ref{repattpower}), $2 \geq A > B \geq 2-d$, complementing previous work by Carrillo and Wang.


\subsection{Steady states and dynamics in the diffusion-dominated regime}\label{subsection:diffusiondominated}

In this subsection, we restrict our attention to the diffusion-dominated regime, where  global well-posedness is  known. The following questions naturally arise:  Are there stationary solutions for all masses? And if so, are they unique up to translations? Do they determine the long time asymptotics of solutions? In this section, we review the existing work on  these questions. We begin, in section \ref{scalingintegrable}, by showing how the same dilation scaling argument that separated global well-posedness from finite time blowup in the diffusion vs. aggregation dominated regimes leads to a second dichotomy in terms of existence of global minimizers. For kernels that grow at least as fast as a power-law $W_k$ at infinity, global minimizers exist for   $m$ in the diffusion dominated regime. On the other hand, for \emph{globally integrable} interaction kernels, the value of $m$ does not depend on the specific structure of the kernel at hand, and the critical diffusion exponent $m_c =2$ separates existence from nonexistence of global minimizers. In section \ref{existenceminimizerssummary} we summarize the precise statements of these results, and in section \ref{steadystatesdynamics}, we discuss in which regimes it is known that solutions of aggregation diffusion equations  asymptotically converge to these steady states.

\subsubsection{A formal scaling argument: power-law kernels vs integrable kernels} \label{scalingintegrable}
Due to the  gradient flow structure of aggregation diffusion equations with respect to the  energy $E[\rho(t)]$, a natural approach for studying long time behavior of solutions is to begin by finding global minimizers of the energy. As explained in section \ref{subsection:wellposed}, considering the energy's scaling under dilations $\rho_\lambda(x) = \lambda^d \rho(\lambda x)$ reveals that in the  diffusion-dominated regime, it is not energy favorable for the density to become concentrated. Consequently, if the energy lacks a global minimizer for a given mass, it must be due to lack of compactness --- that is, there is a minimizing sequence $\{\rho_n\}_{n\in\NN}$ which is not tight, up to translation. 

As in the case for well-posedness vs. blowup, considering the scaling of the energy under dilations reveals a dichotomy separating integrable kernels and kernels that grow like a power-law at infinity. In the pure power-law case $W_k$, $-d < k \leq 0$, equations \eqref{scale_entropy} and \eqref{scale_interaction} for the scaling of the energy show that, when $m>m_c:=1-k/d$,  $E[\rho_\lambda]$ is increasing as $\lambda \to 0^+$. In other words, it is not energy favorable for the mass to spread to infinity. Thus, for power-law kernels $W_k$ in the diffusion dominated regime, we expect a global minimizer for any given mass. Furthermore, if the kernel grows faster than a power law $W_k$ at infinity,  then the long range attractive interactions are even stronger, and a similar scaling argument immediately yields that, for any $m\geq 1$, we again expect a global minimizer for any mass. 

On the other hand, if $W$ is a globally integrable kernel, the interaction energy $\W[\rho]$ scales differently. In particular, we obtain
\[
 \W[\rho_\lambda] = \frac{\lambda^{d}}{2} \int_{\mathbb{R}^d} \rho(\rho * \tilde W_\lambda) dx, \quad \tilde W_\lambda(x) = \lambda^{-d} W(\lambda^{-1}x) . 
\]
If $W$ is integrable, then as $\lambda\to 0$, $\tilde W_\lambda$ approaches a Dirac mass in distribution, and for bounded and continuous densities $\rho$,  $\lim_{\lambda\to 0^+} \lambda^{-d} \W[\rho_\lambda] =  \frac{1}{2} \left( \int W dx \right) \int \rho^2 dx$. As a result, for $0<\lambda\ll 1$, we obtain
\[
\W[\rho_\lambda] =  \frac{\lambda^{d}}{2} \left( \int_{\mathbb{R}^d} W  dx \right) \int_{\mathbb{R}^d} \rho^2 dx + o(\lambda^{d}).\]
Comparing this with the scaling for $\S_m[\rho_\lambda]$ from equation \eqref{scale_entropy}, we formally obtain that $m=2$ is the critical power separating the energies for which it is favorable ($m<2$) vs. unfavorable ($m>2$) for the mass to spread to infinity. Thus, we expect this critical exponent to also determine non-existence vs. existence of global minimizers.

At the critical power $m=2$, the balance between diffusion and aggregation is more delicate.  In this case, 
\[
E[\rho_\lambda] = \lambda^d \int_{\mathbb{R}^d} \rho^2 dx \Big(1 + \frac{1}{2}\int_{\mathbb{R}^d} W  dx\Big) + o(\lambda^d).
\]
Therefore, if $\int_{\mathbb{R}^d} W  dx < -2$ , we formally expect existence of global minimizers, whereas  $\int_{\mathbb{R}^d} W  dx > -2$ leads to non-existence.


\subsubsection{Existence / non-existence of global minimizers} \label{existenceminimizerssummary}
We now turn to the precise statements of what is known concerning existence vs. non-existence of global minimizers. Note that once one obtains  existence of a global minimizer, applying Riesz's rearrangement inequality directly yields that such global minimizer must be radially decreasing for all purely attractive kernels. 

If the interaction potential is a power-law $W_k$, $-d< k\leq 0$,  and $m$ is in the diffusion-dominated regime, $m>m_c := 1-k/d$, then for any   mass $M$,   there exists a global minimizer for the energy $E[\rho]$ in the class 
\[\mathcal{Y} := \{\rho\in L_+^1(\mathbb{R}^d) \cap L^m(\mathbb{R}^d), \|\rho\|_1 = M, \int_{\mathbb{R}^d} x \rho(x) = 0\}.
\] 
Lions was the first to show this, in the case the interaction is potential is Newtonian and $d \geq 3$, via a concentration compactness argument \cite{lions, lions2}. Subsequently, Bedrossian generalized the result to all purely attractive potentials no more singular than Newtonian potential \cite{bedrossian1}. Carrillo, Castorina, and Volzone extended the result to the Newtonian potential in two dimensions  \cite{CCV}. For kernels more singular than Newtonian, recent work by Carrillo, Hoffmann, Mainini, and Volzone  showed existence of global minimizers for power law kernels $W_k$ with $-d<k<2-d$. Furthermore, the same work proved that, for all power law kernels $-d<k \leq 0$,  global minimizers are compactly supported \cite{CHMV}.

If $W$ is a integrable, bounded, and purely attractive kernel, Bedrossian showed that, for all $m >2$ and any mass $M$,  there exists a global minimizer of \eqref{energy}   \cite{bedrossian1}. At the critical power $m=2$,  existence vs. non-existence of global minimizers   depends on the value of $\int W dx$. Bedrossian and Burger, Di Francesco, and Franek showed that, for any mass,     a global minimizer of \eqref{energy} exists if and only if $\int W dx < -2$ \cite{bedrossian1, BDF}. Subsequently, Kaib showed that for $1<m<2$, there is a compactly supported global minimizer if $-\int W dx$ is sufficiently large \cite{kaib}. More recently, Carrillo, Delgadino, and Patacchini showed that, for $m=1$, there is no steady state (thus no global minimizer) for any bounded interaction kernel \cite{CDP}.
 
We conclude by observing that, in all the above results where $W$ is no more singular than the Newtonian potential at the origin, whenever a global minimizer for a given mass is known to exist, standard arguments ensure that it must be a steady state in the weak sense. For power-law kernels more singular than Newtonian, in the diffusion-dominated regime, Carrillo, Hoffmann, Mainini, and Volzone showed that global minimizers are steady states  for all $1-d \leq k < 2-d$; for more singular kernels, $-d < k < 1-d$, they also show that, under additional constraints on $m$, global minimizers have sufficient regularity to be steady states \cite{CHMV}.

\subsubsection{Steady states and dynamics} \label{steadystatesdynamics}
Once existence of a global minimizer is known, it is natural to ask whether it is a global attractor of the evolution equation. Note that a necessary condition for this to be true is that steady states are unique, up to a translation.
Compared to the existence theory, much less is known about uniqueness of steady states. Using continuous Steiner symmetrization techniques, Carrillo, Hittmeir, Volzone, and Yao showed that, for purely attractive kernels, all stationary solutions in $L_+^1 \cap L^\infty(\Rd)$  are radially decreasing  \cite{CHVY}. Therefore, if the kernel is attractive, it suffices to study the uniqueness question among the radial class of functions. 

In the particular case of the Newtonian interaction potential, uniqueness of steady states among radial functions was first shown by Lieb and Yau \cite{LY}, using the specific structure of the Newtonian kernel to obtain an explicit ordinary differential equation for the mass distribution function of solutions. For general power-law kernels in the diffusion-dominated regime, uniqueness of steady states is known when $d=1$ \cite{CHMV} but remains unclear for higher dimensions, aside from the Newtonian case. For integrable, purely attractive kernels, Kaib and Burger, Di Francesco, and Franek proved uniqueness of steady states   at the critical power $m=2$ when $\int W dx < -2$  using the Krein-Rutman theorem, under the additional assumptions that $W$ is smooth and $W'(r)$  vanishes  only at 0 \cite{BDF, kaib}. Aside from these special cases, the uniqueness of steady states remains open. 

We close by reviewing what is known about the dynamics of solutions in the diffusion-dominated regime---in particular, when solutions converge to the previously described steady states and global minimizers. Formally, if uniqueness of steady states, up to translation, is known for a given mass, one would expect that any solution would converge to a translation of it. On the other hand, if there are no global minimizers, one would expect the solution to spread to infinity. However, in practice, long-time asymptotics are only known in special cases. Clarifying the precise conditions on the interaction kernel that lead to convergence towards equilibrium or spreading to infinity   is a challenging open problem related to sharp conditions on the confinement of the mass.

 For the Newtonian potential in two dimensions in the diffusion-dominated regime $m>1$, Carrillo, Hittmeir, Volzone, and Yao showed that every solution with finite initial second moment must converge to a (unique) stationary solution, which has the same mass and center of mass as the initial data \cite{CHVY}. The proof relies on a global-in-time bound of the second moment, leveraging the structure of the two dimensional Newtonian potential. When $d\geq 3$, although the uniqueness of steady states is also known \cite{CHVY}, without any uniform-in-time mass confinement bounds, it is unclear how to show that the mass cannot escape to infinity as $t\to\infty$. However, in the case of radial solutions, Kim and Yao succeeded in showing that that all radial solutions converge to the global minimizer exponentially fast, by building sub/super solutions for the mass concentration functions  \cite{KY} .

For an integrable, purely attractive kernel, at the critical power $m=2$, Di Francesco and Jaafra obtained the following results on the asymptotic behavior in dimension one  \cite{DJ}. If $\int_{\mathbb{R}} W dx < -2$, then the (unique) steady state is locally asymptotically stable in 2-Wasserstein distance, and if $\int_{\mathbb{R}} W dx > -2$ (where non-existence of global minimizer is known), all solutions with finite energy must decay to zero locally in $L^2$ and almost everywhere in $\mathbb{R}$ as $t\to\infty$. When $m\neq 2$, also in dimension one, numerical results by Burger, Fetecau, and Huang  suggest that solutions are attracted to the global minimizer when $m>2$ (with some coarsening and metastability, see next section), and when $1<m<2$, the steady states (if any) have a limited basin of attraction   \cite{BFH}.


\section{Singular limits} \label{singularlimits}

We now turn from the well-posedness theory and long-time behavior discussed in the previous sections to consider dynamics of solutions of aggregation diffusion equations in two limiting regimes. In section \ref{slowdiff}, we consider the behavior of solutions in the \emph{slow diffusion limit}, sending the diffusion exponent $m \to +\infty$. In this case, the effect of diffusion transforms into a hard height constraint on the density, leading to the \emph{constrained aggregation equation}. In section \ref{singularmeta}, we consider singular limits that affect the balance between aggregation and diffusion, leading to metastable behavior, as solutions spend long time scales converging towards a local equilibrium, before quickly transitioning to converge to a global equilibrium.


\subsection{Constrained aggregation and the slow diffusion limit} \label{slowdiff}
We begin by considering the behavior of solutions of aggregation diffusion equations as the diffusion exponent $m \to +\infty$.
Formally, this leads to an aggregation equation with a height constraint $\rho \leq 1$ on the density,
\begin{equation} \label{constrainedagg}
\begin{cases}
\rho_t  =  \nabla\cdot(\rho\nabla(\rho*W)) & \text{ if }\rho<1,\\
``\rho \leq 1 \text{ always.''}
\end{cases}
\end{equation}

This relation between the above constrained aggregation equation and aggregation diffusion equations  can be formally understood by noting that the diffusion exponent $m\geq 1$ controls the strength of diffusion at different heights of the density. Since  
\begin{align*}
\lim_{m\to\infty}\rho^m = \begin{cases} 0 &\text{ for } \rho < 1, \\ +\infty &\text{ for } \rho >1 , \end{cases}
\end{align*} 
in the $m\to\infty$ limit, the diffusion term $\Delta \rho^m$ vanishes in the set $\{\rho<1\}$, whereas it becomes strong enough in the set $\{\rho>1\}$  to prevent the density from growing above  height 1. In what follows, we will clarify the  sense in which we enforce the height constraint. (See Figure \ref{varyingdiffexp} for a numerical illustration of how degenerate diffusion with $m$ large approximates a height constraint.)

For various partial differential equations with a degenerate diffusion term $\Delta \rho^m$, it is well known that the asymptotic limit $m\to\infty$ imposes a height constraint $\rho\leq 1$, and the evolution of the free boundary $\partial\{\rho=1\}$ is often determined by a  Hele--Shaw type free boundary problem. In particular, for the porous medium equation $\rho_t = \Delta \rho^m$, this is known as the Mesa Problem, and it has been proved that the solution pair $(\rho_m, p_m)$ (where $p_m  := \frac{m}{m-1}\rho^{m-1}$ is the pressure) converges as $m\to\infty$ to a solution of the Hele--Shaw problem   \cite{CF87, GQ1, GQ2}. More recently, the analogous limit has been considered for the porous medium equation with growth $\rho_t = \Delta \rho^m + \rho \Phi(p_m)$, where the limiting equation  describes the growth of a tumor with a restriction on the maximal cell density \cite{PQV, PQTV, MPQ, KP}.

Another viewpoint for understanding the $m\to\infty$ limit in the case of aggregation diffusion equations is  through the gradient flow formulation. Recall that the aggregation diffusion equation \eqref{aggregation} is formally the Wasserstein gradient flow of the free energy functional 
\begin{equation*}
E_m[\rho] = \frac{1}{m-1}\int_{\mathbb{R}^d}\rho^m dx + \frac{1}{2}\int_{\mathbb{R}^d} \rho(\rho*W)dx =: \S_m[\rho] + \W[\rho].
\end{equation*} 
For any fixed $\rho$, it is easy to check that $\lim_{m\to\infty} E_m[\rho] = E_\infty[\rho]$, which is given by the constrained interaction energy
\begin{equation*}
E_\infty[\rho] = 
\begin{cases}
\frac{1}{2}\int_{\mathbb{R}^d} \rho(\rho*W)dx & \text{ if } \|\rho\|_{L^\infty} \leq 1\\
+\infty & \text{ otherwise}.
\end{cases}
\end{equation*}
Thus we formally expect the slow diffusion limit of aggregation diffusion equations should correspond to the Wasserstein gradient flow of $E_\infty[\rho]$, which is indeed the constrained aggregation equation (\ref{constrainedagg}). In fact, this is the precise sense in which we impose the height constraint ``$\rho \leq 1$ always.''

In the case when  the nonlocal interaction energy $\frac{1}{2}\int_{\mathbb{R}^d} \rho(\rho*W)dx$ in $E_\infty$ is replaced by a local potential energy $\int \rho(x) V(x) dx$, with a $\lambda$-convex potential $V(x)$, the gradient flow for $E_\infty$ was introduced by Maury, Roudneff-Chupin, Santambrogio, and Venel as a model for pedestrian crowd motion \cite{MRS, MRSV}. They showed that this gradient flow satisfies the transport equation $\rho_t + \nabla\cdot(\rho\vec{v})=0$  in the weak sense, where the velocity field $\vec{v}(\cdot,t)$ is given by the $L^2$ projection of $\nabla V$ onto the set of admissible velocities that do not increase the density in the saturated zone $\{\rho(\cdot,t)=1\}$. Building upon this work, Alexander, Kim, and Yao showed that the gradient flow of $E_\infty$ (with local potential $V$) can be approximated by the gradient flows  of $E_m$, i.e. nonlinear Fokker-Planck equations, as $m\to\infty$, proving convergence in the Wasserstein distance with an explicit  rate depending on $m$ \cite{AKY} . 

Building upon this work in the local case, Craig, Kim, and Yao studied the gradient flow of the constrained interaction energy $E_\infty$, when $W$ is an attractive Newtonian interaction potential \cite{CKY}. Craig showed that the energy $E_\infty$ is \emph{$\omega$-convex}, ensuring that its gradient flow is well-posed and can be quantitatively approximated by the discrete time Jordan-Kinderlehrer-Otto scheme \cite{CraigNonconvex,JKO}. Craig, Kim, and Yao then  showed that a patch solution of the constrained aggregation equation (where the initial data is given by a characteristic function) remains a patch for all time, and the evolution of the patch boundary satisfies a Hele-Shaw type free boundary problem. In addition, in two dimensions, these patch solutions converge to a characteristic function of a disk in the long-time limit. A key element in the proof was the approximation of the constrained aggregation equation by a sequence of nonlinear Fokker-Planck equations in the slow diffusion limit, where the local potential $V$ depended on the choice of initial data. However, for nonconvex interaction potentials, including the Newtonian interaction potential $W$ under consideration, convergence of aggregation diffusion equations to the constrained aggregation equation as $m \to +\infty$ remained open.

More recently, Craig and Topaloglu \cite{CraigTopalogluHeight} proved that minimizers and gradient flows of $E_m$ do converge to minimizers and gradient flows of $E_\infty$ in the $m\to\infty$ limit for a general class of interaction potentials $W$, including both attractive and repulsive-attractive power-law potentials with singularity up to and including the Newtonian potential. This work was largely inspired by the following related question in geometric shape optimization: can we characterize \emph{set-valued} or \emph{patch }minimizers $E_\infty$, i.e. minimizers which are of the form $\rho = 1_\Omega$ for $\Omega \subseteq \Rd$ of volume $M$? When the interaction potential is a repulsive-attractive power law (\ref{repattpower}), competition between the attraction parameter $A$ and the repulsion parameter $B$ determines existence, nonexistence, and qualitative properties of minimizers, providing a counterpoint to the well-studied nonlocal isoperimetric problem \cite{ChMuTo2017}. Due to the fact that solutions of aggregation diffusion equations approximate the constrained aggregation equation in the slow diffusion limit, numerical simulations of aggregation diffusion equations for large $m$ can be used to explore qualitative properties of set-values minimizers of $E_\infty$. (See Figures \ref{evolutiontoequilibriumx4mx2} and \ref{differentrepulsioncriticalmass} for numerical simulations illustrating critical mass behavior.)

Several open questions remain for the height constrained aggregation equation, stemming from the fact that posing the equation as a Wasserstein gradient flow   offers only a very weak notion of solution. In the particular case of the attractive Newtonian interaction potential with patch initial data, Craig, Kim, and Yao's result provides a true characterization in terms of a partial differential equation \cite{CKY}. It remains open whether this result can extended to more general initial data, building on recent work in the local case \cite{MPQ, KP,kim2017singular}. An obstacle is the very low regularity of solutions and the variety of possible merging phenomena that may occur as the solution interacts with the height constraint. (See Figure \ref{heightconstrNewt} for a numerical illustration of possible merging behavior.) Likewise, it is unknown how to extend this  result to  more general interaction potentials.

\subsection{Singular limits and metastability} \label{singularmeta}

As a counterpoint to the previous section, which considered singular limits under which the effect of diffusion transforms into a height constraint, in the section, we consider singular limits affecting the balance between aggregation and diffusion, leading to various types of metastability behavior. 
In the absence of diffusion, much is known about the long-time behavior of solutions of the aggregation equation,
\begin{align} \label{aggeqn}
\rho_t = \nabla \cdot (\rho\nabla (\rho*W)) ,  
\end{align}
for both purely attractive interaction potentials $W$ and interaction potentials that are repulsive at short length scales and attractive at longer length scales. Solutions approach compactly supported stationary states, and the regularity of these equilibria is determined by the the repulsive strength of the interaction potential at the origin \cite{BCLR2,CCP,CDM,CFP}. Furthermore, solutions may experience either finite or infinite time blow-up, concentrating at Dirac masses or on other lower dimensional sets \cite{BCP,BertozziCarrilloLaurent,BertozziLaurentRosado,BertozziGarnettLaurent,BertozziLaurentLeger,5person}.

A natural question is how robust this behavior is under different perturbations, and in particular, how asymptotics are affected by the addition of linear or nonlinear diffusion. When are the above described stationary states preserved? And if steady states exist, how stable are they?

The appearance of metastable behavior in nonlinear aggregation-diffusion equations is apparent from numerical simulations when the interactions have nearly finite radius of perception, i.e. when the interaction kernel decays quickly at $+\infty$ \cite{EK,BFH,CCH,CCWW}.   Moreover, for compactly supported attractive or repulsive-attractive interaction potentials with nonlinear diffusions, clustered solutions can form and the support of steady states can contain several disconnected components.  \cite{MT,ABPP}. 
Can the local effect of the nonlinear diffusion as repulsion lead to clusterization for fully attractive potentials? If the fully attractive potential satisfies $W'(r)>0$ everywhere, then we know that all steady states must be radially decreasing \cite{CHVY}, so this effect cannot happen for large times. So, can this clusterization happen as a metastable behavior in the dynamics? 

While these questions remain mostly open, we now describe two different directions of recent progress. First,  we will consider  metastability that appears for porous medium diffusion ($m=2$) when the attraction potential is localized via dilation. Next, we will discuss metastability  of aggregation diffusion equations for linear diffusion ($m=1$) as the amount of diffusion vanishes, analogous to the well-known vanishing viscosity limit in classical fluids equations.
\subsubsection{Metastability in the localized attraction limit} \label{metaatt} 
We begin by considering the behavior of aggregation diffusion equations when the interaction potential localizes as $\delta \to 0$,
\begin{align} \label{agglocatt}
\rho_t = \nu \Delta \rho^m + \nabla \cdot (\rho\nabla(\rho*W_\delta)),   \ m>1, \  W_\delta(x) = \delta^{-d} W(\delta^{-1}x), \ 0<\delta\ll 1 .
\end{align}
For simplicity, we consider the case when $W$ is a bounded, strictly attractive interaction kernel. The attraction of $W_\delta$ localizes as $\delta \to 0$, since the dilation  decreases attractive forces in the long range and increase them in a $\delta$-neighborhood.

For this class of equations, Carrilllo, Hittmeir, Volzone, and Yao showed that, if a stationary state exists, it must be radially decreasing \cite{CHVY}. However,  numerical results illustrate that the dynamical solution develops into non-radially decreasing clusters which are approximate steady states, where the length scale of each cluster goes to zero as $\delta \to 0$. Eventually, these clusters coalesce into a single cluster, which ultimately converges to a radially decreasing steady state;  see Figures \ref{varyinginteractionrange}  and \ref{varyinginteractionrangelongtime}.  Similar metastability phenomena have been observed for fixed  $\delta >0$ in several related numerical studies \cite{CCH,BFH,CCWW}. The structure, duration, and existence of metastable states is strongly dependent on the smoothness of the interaction kernel $W$ and the size of the diffusion coefficient $\nu$. In particular, Figure \ref{varyingdiffcoeff} illustrates how increasing the diffusion exponent can prevent the formation of steady states, and Figures \ref{metalogfig} and \ref{metanewtfig} illustrate that a singularity in the interaction kernel $W$ can also inhibit their formation.

Little is known about how the time scales of this metastability behavior relate to the localization parameter $\delta$, the diffusion exponent $m>1$, the diffusion coefficient $\nu$, and the choice of interaction potential $W$.  However, in the special case of aggregation-diffusion equations on a bounded interval, with no-flux boundary conditions, one can show that, for initial data given by a constant function, the length scale of the bumps in the resulting metastable state can be computed by analyzing the instability of the constant solution.

\subsubsection{Metastability and the vanishing diffusion limit} \label{metadiffvanish}

A natural context in which to examine robustness of solutions of the aggregation equation (\ref{aggeqn}) under perturbations is to consider behavior of  aggregation diffusion equations, with a small amount of linear diffusion,
\begin{align} \label{aggvanishdiff}
\rho_t = \nu \Delta \rho + \nabla \cdot (\rho\nabla(\rho*W)) \quad \text{ in } \mathbb{R}^d.
\end{align}
On one hand, for  $\lambda$-convex interaction potentials $W$,  classical $\Gamma$-convergence techniques yield that solutions of (\ref{aggvanishdiff}) converge to solutions of the inviscid aggregation equation (\ref{aggeqn})  as $\nu \to 0$, on any bounded time interval $[0,T]$  \cite{AGS, Serfaty}. Furthermore, Cozzi, Gie, and Kelliher extended this result to the case when $W$ is the attractive Newtonian potential, on any time interval $[0,T]$, as long as the corresponding the inviscid solution remains well-defined \cite{CGK}.

On the other hand, the analogous convergence result fails in the long time limit for all smooth, bounded attractive kernel $W$. For such kernels, solutions of the inviscid aggregation equation (\ref{aggeqn}) always concentrate to one or more Dirac masses as $ t \to +\infty$: the solution converges to a single Dirac mass if $\rho_0$ is initially within the perception range of $W$, otherwise it may cluster into multiple delta functions \cite{motsch2014heterophilious}. However, Carrillo, Delgadino, and Patacchini recently showed that, for all such interaction kernels, equation (\ref{aggvanishdiff}) has no steady states for any $\nu >0$ \cite{CDP}; see Figure \ref{lineardiffusionmetastablefig}. Consequently, there is no hope of proving convergence to any nonzero equilibria as $t\to\infty$ for any $\nu>0$. (Note that this result on nonexistence of steady states strongly uses the boundedness of $W$ and non-compactness of $\mathbb{R}^d$).

Therefore, when the interaction potential $W$ is bounded, equation (\ref{aggvanishdiff}) illustrates a typical example of non-commutative  double limits: $\nu\to 0$ and $t\to \infty$. For small $\nu>0$, we expect solutions to remain close to the equilibrium of the inviscid aggregation equation (\ref{aggeqn}) for a large time interval $[0,T]$, but  to eventually dissipate away. Consequently, numerical simulations of (\ref{aggvanishdiff}) for small $\nu>0$ may be quite misleading: it may appear that solutions are converging towards a stationary state before  eventually becoming evident that mass is spreading at an extremely slow rate. It is still unknown how this time scale depends on $\nu>0$ and how solutions dynamically transition from near the inviscid steady states to then glue together and spread to $+\infty$  \cite{CCS}. 

 Related metastability behavior has also been observed for interaction kernels that are unbounded away from the origin. For example, Evers and Kolkolnikov discovered that, in the case of the repulsive-attractive interaction potential $W(x) = \frac{x^4}4-\frac{x^2}2$ in $d=1$, the effect of the diffusive regularization in equation (\ref{aggvanishdiff}) is to produce a extremely slow equilibration of the weights between the two stable points of the dynamics.
 To see this, recall that for the inviscid aggregation equation (\ref{aggeqn}), any convex combination of Dirac masses at distance one apart is a stationary solution (since $W'(x) = x^3 - x$ has a zero at $x=1$), and  the mass need not be evenly distributed. However, for all $\nu >0$, solutions of the aggregation diffusion equation (\ref{aggvanishdiff}) seek to equilibrate mass. Numerical simulations indicate that there is a unique equilibrium for all $\nu >0$, which converges as $\nu \to 0$ to two Dirac masses with equal weight; see Figure \ref{metax4mx2}. Still, proving uniqueness of equilibria for the $\nu >0$ case and understanding the timescales on which these effects take place remains entirely open.


\section{Numerical methods} \label{numericsmainsection}

As described in the previous sections, two key themes in the analysis of aggregation diffusion equations are the gradient flow structure and the competition between local and nonlocal effects. The gradient flow structure provides Lyapunov functionals, stability estimates, and a natural framework for considering singular limits, as in the case of the vanishing diffusion limit or the slow diffusion to height constraint limit. The delicate balance between aggregation and diffusion equations   leads to a range of asymptotic behavior, from global existence to finite time blowup, and rich structure of equilibrium configurations.

The same key themes arise in the development of numerical methods for aggregation diffusion and height constrained aggregation equations. Due to the competition between aggregation and diffusion or aggregation and a height constraint, numerical methods must be able to cope with solutions with very low regularity and   to perform well enough globally in time to accurately approximate equilibrium profiles. With regard to the Wasserstein gradient flow structure, numerical methods must at a minimum respect the space to which solutions belong ---  preserving positivity and conserving mass --- 
as well as ideally reproducing the energy decreasing property and the rate of energy decrease.

One of the first numerical methods developed for aggregation diffusion equations directly leveraged the Wasserstein gradient flow structure, and in particular, a time discretization of the gradient flow problem due to Jordan, Kinderlehrer, and Otto \cite{JKO}. In analogy with the implicit Euler method for approximating gradient flows in Euclidean space, a Wasserstein gradient flow (\ref{formal W2 grad flow}) with initial data $\rho_0 \in \P_2(\Rd)$ may be approximated by solving the sequence of infinite dimensional minimization problems, 
\begin{align} \label{JKOdef} \rho^{n} =\argmin_{\nu \in \P_2(\Rd) } \left\{\frac{1}{2(t/N)} d_2^2(\nu,\rho^{n-1}) + E(\nu) \right\} , \quad  n= 1, \dots, N .
\end{align}
Under weak regularity and coercivity assumptions on $E$, we have
\[ \lim_{N \to +\infty} d_2(\rho^N, \rho(t)) =0 , \]
 where $\rho(t)$ is the solution of the Wasserstein gradient flow at time $t$ \cite{AGS}. Furthermore, provided that the energy satisfies generalized convexity requirements, such as $\lambda$-convexity or $\omega$-convexity, one may obtain quantitative estimates on the rate of convergence \cite{AGS, CraigNonconvex,CC,CraigJKO}.

 A first benefit of using the JKO scheme numerically is that, by definition of the sequence, the energy is decreasing in $\rho^n$, which provides automatic stability estimates along the sequence. Furthermore, since the minimization problem is constrained over $\P_2(\Rd)$, this approach conserves mass and preserves positivity of solutions. More generally, the JKO scheme can be easily extended to any initial data with mass $\int \rho_0 = M  >0$ by replacing $\P_2(\Rd)$ with the set of finite measures with mass $M$ and finite second moment. Finally, since this approach to numerical methods most closely respects the equation's underlying gradient flow structure, such methods are a natural choice for simulating singular limits of aggregation diffusion equations, including height constrained aggregation equations \cite{BCL16,carlier2017convergence}.

As the JKO scheme (\ref{JKOdef}) is already discrete in time, to apply it numerically, one simply needs to determine a method for discretizing the infinite dimensional minimization problem (\ref{JKOdef}) in space. The key difficulty is computing the Wasserstein distance term accurately and efficiently. One of the first techniques took advantage of the Monge formulation of the Wasserstein term, discretizing the space of Lagrangian maps via finite difference \cite{ESG2005, CM} or finite element approximations \cite{CRW}, or directly using a discretization of the Monge-Amp\'ere operator \cite{BCMO}. A more recent approach is to leverage the Benamou-Brenier dynamical reformulation of the Wasserstein distance, which can then either be discretized via Benamou and Brenier's original ALG-2 method \cite{BB00, BCL16} or via modern operator splitting techniques \cite{PPO14,CCWW}. Aside from these two main approaches, several other methods for discretizing the Wasserstein term have also been applied to simulate aggregation diffusion equations \cite{BlanchetCalvezCarrillo,CG, GT,WW}, including, most recently, using the Kantorovich formulation of the Wasserstein distance, with entropic regularization to improve computational efficiency \cite{carlier2017convergence}.
 
A second common approach to developing numerical methods for aggregation diffusion equations is to discretize the equations directly via classical Eulerian methods, based on similar schemes for related hyperbolic, kinetic, and degenerate parabolic equations. Leveraging the analogy with degenerate diffusion equations, Carrillo, Chertock, and Huang developed a  finite volume method for aggregation diffusion equations  \cite{CCH, Filbet}. Building on this, recent work by Bailo, Carrillo, and Hu proposes  and implicit in time fully discrete entropy decreasing finite volume schemes for general aggregation-diffusion equations  \cite{BCH}. Alternatively, using the analogy with hyperbolic conservation laws, Sun, Carrillo, and Shu proposed a finite element method, using a dicontinuous Galerkin approach with gives higher order convergence to smooth solutions \cite{ZCS}. For the particular case of attractive Newtonian aggregation, Liu, Wang, and Zhou    introduced a change of variables by which the equation shares structural similarities with Fokker-Planck equations and used this to develop an implicit in time finite difference scheme \cite{LiuWangZhou}. Each of the previous methods preserves positivity of solutions, and numerical experiments indicate they also succeed in dissipating energy, providing  stability for numerical solutions and allowing the methods to accurately capture asymptotic behavior. 

A third approach to numerical methods for aggregation diffusion equations, and the approach most closely related to the  microscopic  interacting particle system underlying the partial differential equation, is to consider particle approximations for aggregation diffusion equations. Much of the existing work in this direction has leveraged the structural similarity between aggregation diffusion equations and equations from classical fluids, particularly the Navier-Stokes equation in vorticity form \cite{BertozziLaurentLeger, CB}. In the absence of diffusion, the simplest type of particle method for the aggregation equation 
\begin{align} \label{aggeqn2}
\rho_t  =  \nabla\cdot(\rho\nabla(\rho*W))\quad \text{ in } \mathbb{R}^d,
\end{align}
proceeds by discretizing the initial datum $\rho_0$ as a finite sum of $N$ Dirac masses,
\begin{align} \label{partinitial1}
 \rho_0 \approx \rho_0^N = \sum_{i = 1}^N \delta_{x_i} m_i, \qquad x_i \in \Rd, \quad m_i \geq 0 , 
 \end{align}
where $\delta_{x_i}$ is a Dirac mass centered at $x_i \in \Rd$, and then evolving the locations of the Dirac masses according to the velocity field from equation (\ref{aggeqn2}),
 \begin{align} \label{partsol1} 
 	\rho^N(t) =  \sum_{i = 1}^N \delta_{x_i(t)} m_i, \quad \dot{x}_i = - \sum_{j \neq i } \grad W(x_i - x_j) m_j .
\end{align}
 
For $W \in C^1(\Rd)$,  the particle solution $\rho^N(t)$ is then a weak solution of the aggregation equation (\ref{aggeqn}) with initial data $\rho_0^N$, and in particular, it is a Wasserstein gradient flow of the interaction energy $\W[\rho]$, so that the energy $\W[\rho]$ automatically decreases along $\rho_N(t)$. Furthermore, results on the mean-field limit for the aggregation equation (\ref{aggeqn}) ensure that, for a range of $\lambda$-convex or power law interaction kernels $W$, as the approximation of the initial data improves, $\rho_0^N \xrightarrow{d_2} \rho_0$, the particle solution converges to the exact solution of the aggregation equation $\rho^N(t) \xrightarrow{d_2} \rho(t)$ on bounded time intervals \cite{CarrilloChoiHauray,5person,Jabin}. The particle method (\ref{partsol1})  provides a semi-discrete numerical method, and in order to obtain a fully discrete scheme, one may use a variety of fast solvers to integrate the system of ODEs.

In order to develop methods with higher-order accuracy and capture competing effects in repulsive-attractive systems, recent work has considered enhancements of  standard particle methods inspired by techniques from classical fluid dynamics, including \emph{vortex blob methods} and \emph{linearly transformed particle methods}  \cite{CB,CCCC, GHL}. Bertozzi and the second author's blob method for the aggregation equation obtained a higher order accurate method for singular interaction potentials $W$ by convolving $W$ with a mollifier $\varphi_\e(x) = \varphi(x/\epsilon)/\epsilon^d$, $\epsilon >0$. In terms of the Wasserstein gradient flow perspective this translates into regularizing the interaction energy $(1/2) \int \rho (W*\rho) \,d x$ as $(1/2) \int \rho (W*\varphi_\e*\rho)\,dx $.

To extend particle methods to aggregation \emph{diffusion} equations, one has to confront a fundamental obstacle: unlike in the pure aggregation case, solutions of aggregation diffusion equations with particle initial data (\ref{partinitial1}) do not remain particles. One way to address this issue is to simulate the aggregation and diffusion terms separately, via a splitting scheme, using a classical finite volume method for the diffusion term \cite{YaoBertozzi}. Another natural approach in the case of linear diffusion ($m=1$) is to consider a stochastic particle method, in which Brownian motion is added to the differential equation for the motion of the particles (\ref{partsol1}) \cite{Jabin,JW,Liu1,Liu2}. The main practical disadvantages of such stochastic methods is that the simulations must be averaged over a large number of runs to compensate for the randomness of the approximation and such methods have not been extended to the case of degenerate diffusion $m>1$. 

More recently, two \emph{deterministic} particle methods for aggregation diffusion equations have been introduced, inspired by classical particle-in-cell methods in fluid, kinetic, and plasma physics equations \cite{Russo2,CR,DegondMustieles,Russo,LacombeMasGallic,LionsMasGallic,MGallic, Chertock}.
 The first, due to Carrillo, Huang, Patacchini, Sternberg, and Wolansky, approximates one dimensional aggregation diffusion equations by discretizing the energy using non-overlapping balls centered at the particles \cite{CHPW,CPSW}. The second, due to Carrillo, Craig, and Patacchini, extends naturally to all dimensions $d \geq 1$ by considering a  regularization of the energy similar to Craig and Bertozzi's blob method for the aggregation equation, leading to a \emph{blob method for diffusion} \cite{carrillo2017blob}. 
 
 We now turn to a precise description of this blob method for diffusion, which we apply in the next section to generate our numerical examples. Given a Gaussian mollifier $\varphi_\epsilon(x)$, the regularized energy is defined by
 \begin{align*} 
 \E_\e(\rho) =  
	\displaystyle \int \frac{(\varphi_\e*\rho)^{m-1}}{m-1}\,d \rho + \frac12 \int (W*\rho) \rho .
\end{align*}
Unlike gradient flows of the original aggregation diffusion energy $\E(\rho)$ defined in equation (\ref{energy}),   gradient flows of $\E_\e(\rho)$ for $\e >0$ with particle initial data  remain particles for all time, and the evolution of the particles is determined by the following system of ordinary differential equations
 \begin{align} \label{partsol2} 
	\dot{x}_i &= - \sum_{j =1}^N \grad W(x_i - x_j) m_j  + \sum_{j=1}^N  \grad \varphi_\e(x_i- x_j)  m_j \left( \sum_{k=1}^N \varphi_\e(x_j-x_k)m_k \right)^{m-2} \\
	&\qquad  + \left(\sum_{j=1}^N \varphi_\e(x_i-x_j) m_j \right)^{m-2} \left(\sum_{j=1}^N \grad \varphi_\e(x_i-x_j)m_j \right) . \label{partsol3}
\end{align}
As $\e \to 0$, Carrillo, Craig, and Patacchini show that, for any lower-semicontinuous interaction potential $W$,  the regularized energies $\E_\e$ $\Gamma$-converge to the unregularized energy $\E$ for all $m \geq 1$. Furthermore, for $W$ $\lambda$-convex and $m \geq 2$, they show that gradient flows of the regularized energies $\E_\e$ are well-posed. Finally, provided that sufficient a priori estimates hold along the flow, gradient flows of the regularized energies converge to the solution of the aggregation diffusion equation with initial data $\rho_0$ as $\e \to 0$ and $d_2(\rho_0^N, \rho_0) \to 0$. More recently, Craig and Topaloglu have demonstrates  numerically that  this method also provides a robust approach for simulating aggregation diffusion equations, by sending the diffusion exponent $m \to +\infty$ as the regularization and the discretization of the initial data are refined \cite{CraigTopalogluHeight}.

\section{Simulations via the blob method for diffusion} \label{numericssection}
In this section, we apply the blob method for diffusion to illustrate several properties of the singular limits discussed in section \ref{singularlimits}. We begin in section \ref{numericalimplement} by describing  the details of our numerical implementation, which include various refinements over previous works, such as regridding to reduce the number of particles required for convergence \cite{carrillo2017blob, CraigTopalogluHeight}. In section \ref{heightnumerics}, we provide several numerical examples of the slow diffusion limit and properties of the constrained aggregation equation, particularly critical mass behavior relating to open problems in geometric shape optimization \cite{BuChTo2016,CraigTopalogluHeight,FrankLieb}. In section \ref{metanumerics}, we give numerical examples illustrating the relationship between singular limits and metastability behavior, both as aggregation becomes localized and as diffusion vanishes.

\subsection{Numerical implementation} \label{numericalimplement}
We now describe our numerical implementation of the blob method for diffusion, which we perform in Python, using the Numpy, SciPy, and Matplotlib libraries \cite{NumPy,SciPy,matplotlib}. As explained in the previous section, there are two key steps in  the blob method for diffusion: first, one must approximate the initial data $\rho_0$ by a sum of Dirac masses (\ref{partinitial1}); then, one must evolve the locations of those Dirac masses by solving a system of ordinary differential equations (\ref{partsol2}-\ref{partsol3}).

We begin by describing  the approximation of the initial data (\ref{partinitial1}). In the following, we typically consider examples in which the solution is entirely supported on the spatial domain $[-1,1]$. Unless otherwise specified, we approximate the initial data on the computational domain $[-1.1,1.1]$ by partitioning the domain into $N$ intervals  and placing a Dirac mass at the center of each interval, weighted according to the integral of the initial data over the interval. We let $h$ denote the width of these initial intervals.

In order to solve the system of ordinary differential equations (\ref{partsol2}-\ref{partsol3}), we take the mollifier $\varphi_\epsilon$ to be a Gaussian 
\begin{align*}
 \varphi_\e(x) = \frac{1}{(4\pi \epsilon^2)^{d/2}} e^{-|x|^2/4 \epsilon^2} .
 \end{align*}
We then solve the system of ordinary differential equations using the SciPy \lstinline{solve_ivp} implementation of the backward differentiation formula (BDF) method. Similarly to analogous work on blob methods in the fluids case \cite{AndersonGreengard}, we observe that the numerical error due to the choice of time discretization is of lower order than the error due to the regularization and spatial discretization. When the interaction kernel $W$ has a Newtonian or stronger singularity (e.g. $W(x) = |x|$ or $W(x) = \log |x|$ in one dimension), we also mollify the interaction potential by convolution with the mollifier, as Bertozzi and the second author demonstrated this provides higher order rates of convergence \cite{CB}.

In order to pass from the particle approximation
\[ \rho_N(t) =\sum_{i=1}^N \delta_{x_i(t)}m_i \]
 to a density that can be compared with exact solutions and visualized, we convolve our particle approximation with the mollifier $\varphi_\epsilon$, leading to the following approximate density, that is defined on all of Euclidean space:
 \begin{align} \label{approxden} \rho_{\e,N}(x) = \sum_{i =1}^N \varphi_\e(x-x_i(t)) m_i . 
 \end{align}
 
As in previous work on the Euler equations \cite{BealeMajda1, BealeMajda2} and the aggregation equation \cite{CB}, we observe the fastest rate of convergence if the regularization parameter $\epsilon$ scales according to the initial grid spacing $h$ according to
 \begin{align} \label{epsscaling} h^{1- p} \leq \epsilon   , \text{ for } 0 < p \ll 1. \end{align}
In the fluids community, it is well known that this relationship \eqref{epsscaling} between the interparticle distance and the regularization needs to hold globally in time for the numerical simulation to agree well with exact solutions. In previous work on aggregation and aggregation diffusion equations \cite{CB, carrillo2017blob, CraigTopalogluHeight}, this was accomplished by taking the spatial discretization $h$ very small. In the present work, we accomplish this by using the formula for the approximate density (\ref{approxden}) to re-initialize our particle approximation (\ref{partinitial1}) whenever the maximum interparticle distance exceeds $1.5h$. Numerical examples illustrate that such remeshing is much more computationally efficient than taking $h$ very small.

In the following numerical examples, we  choose our initial data to be given by characteristic functions or Barenblatt profiles, which we define as follows
 \begin{align*}
 1_\Omega(x) &= \begin{cases} 1 &\text{ if } x \in \Omega \\ 0 &\text{ otherwise.} \end{cases} \\ 
 	\rho_{\alpha}(x,\tau) &=  
 \tau^{-d\beta}\left(\kappa -  \frac{\beta}{2} \left( \frac{\alpha-1}{\alpha} \right) \tau^{-2\beta} |x|^2 \right)_+^{(\alpha-1)^{-1}} ,  \ \beta = \frac{1}{2+d(\alpha-1)} ,   
 \end{align*}
 with $\alpha >1$, $\tau >0$, and $\kappa = \kappa(\alpha,d) >0$ chosen so that $\int \rho_{\alpha} \,dx= 1$.

\subsection{Numerical examples: height constrained aggregation} \label{heightnumerics}

We begin with several numerical examples illustrating the slow diffusion limit and properties of the height constrained aggregation equation. See section \ref{slowdiff} for a discussion of this singular limit and the constrained aggregation equation.

\begin{figure}[H]
\vspace{-.25cm}
\begin{center}
 {\bf Equilibria for Varying Diffusion Exponent, $W(x) = |x|^4/4- |x|$}  \\ 
\vspace{.1cm}
 {\bf  mass = 0.6} \\
\includegraphics[height=4cm,trim={.6cm .7cm .6cm .7cm},clip]{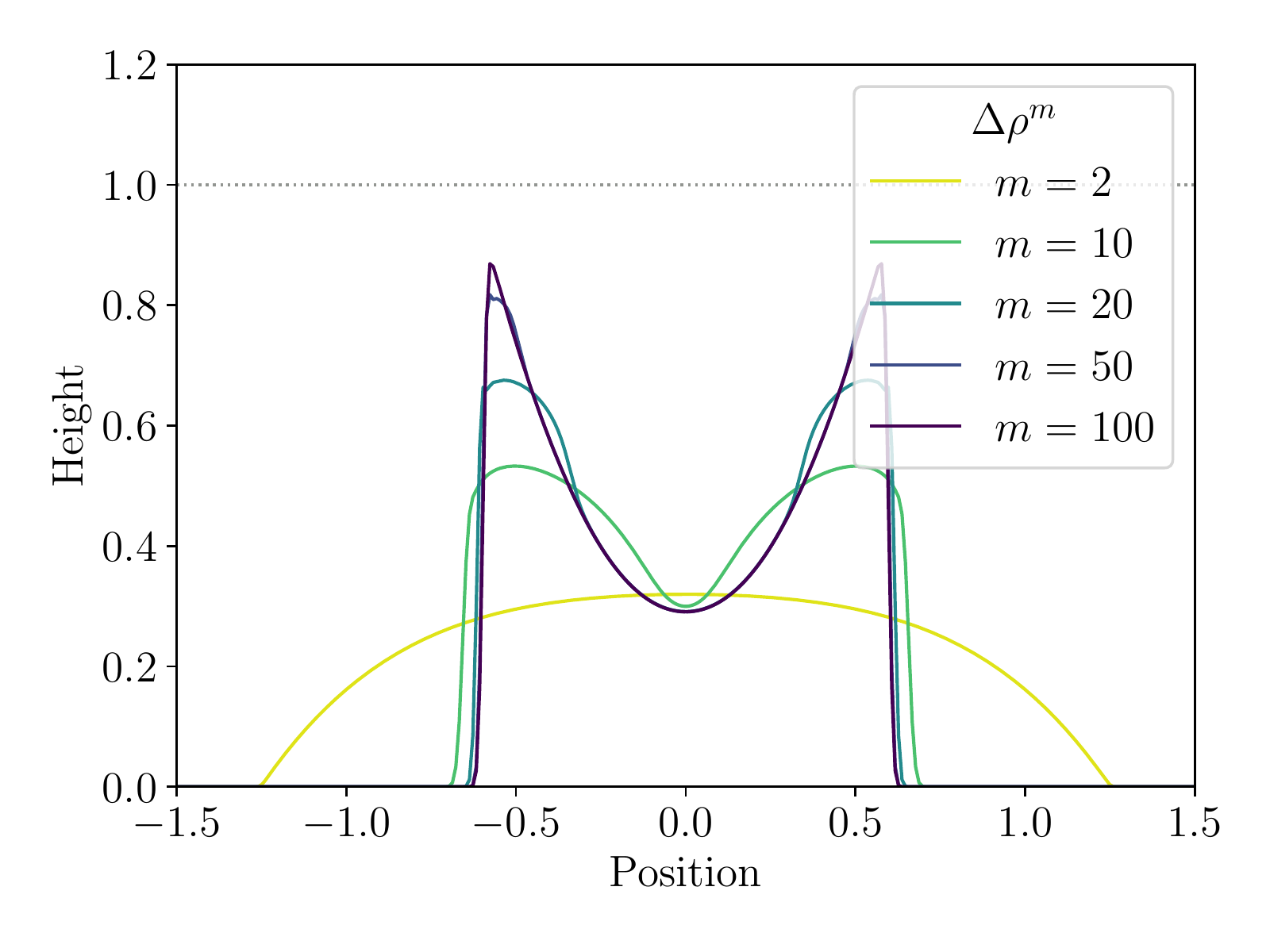}

 {\bf  mass = 1.0} \\
\includegraphics[height=4cm,trim={.6cm .7cm .6cm .7cm},clip]{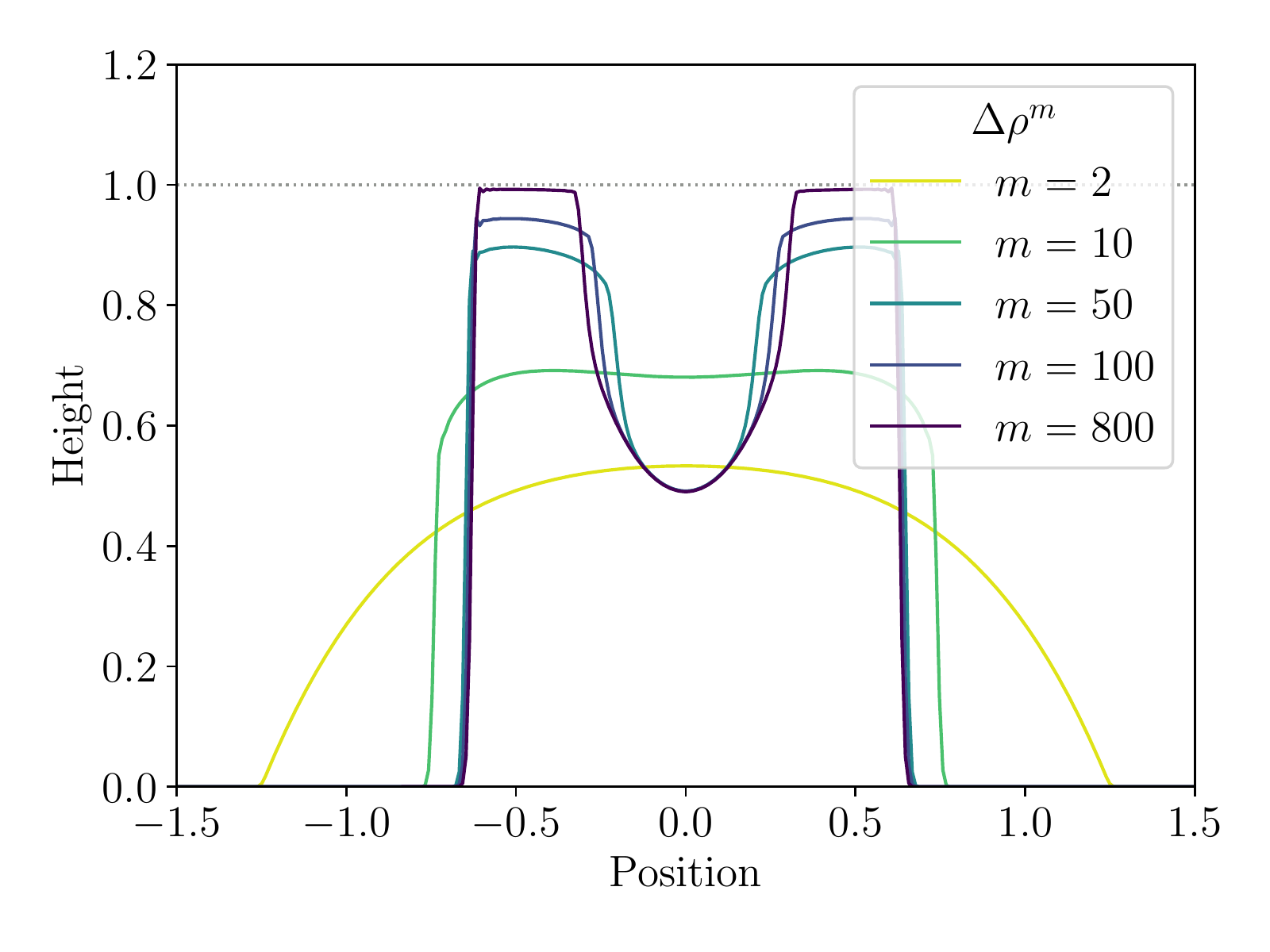}

 {\bf  mass = 1.4} \\
\includegraphics[height=4cm,trim={.6cm .7cm .6cm .7cm},clip]{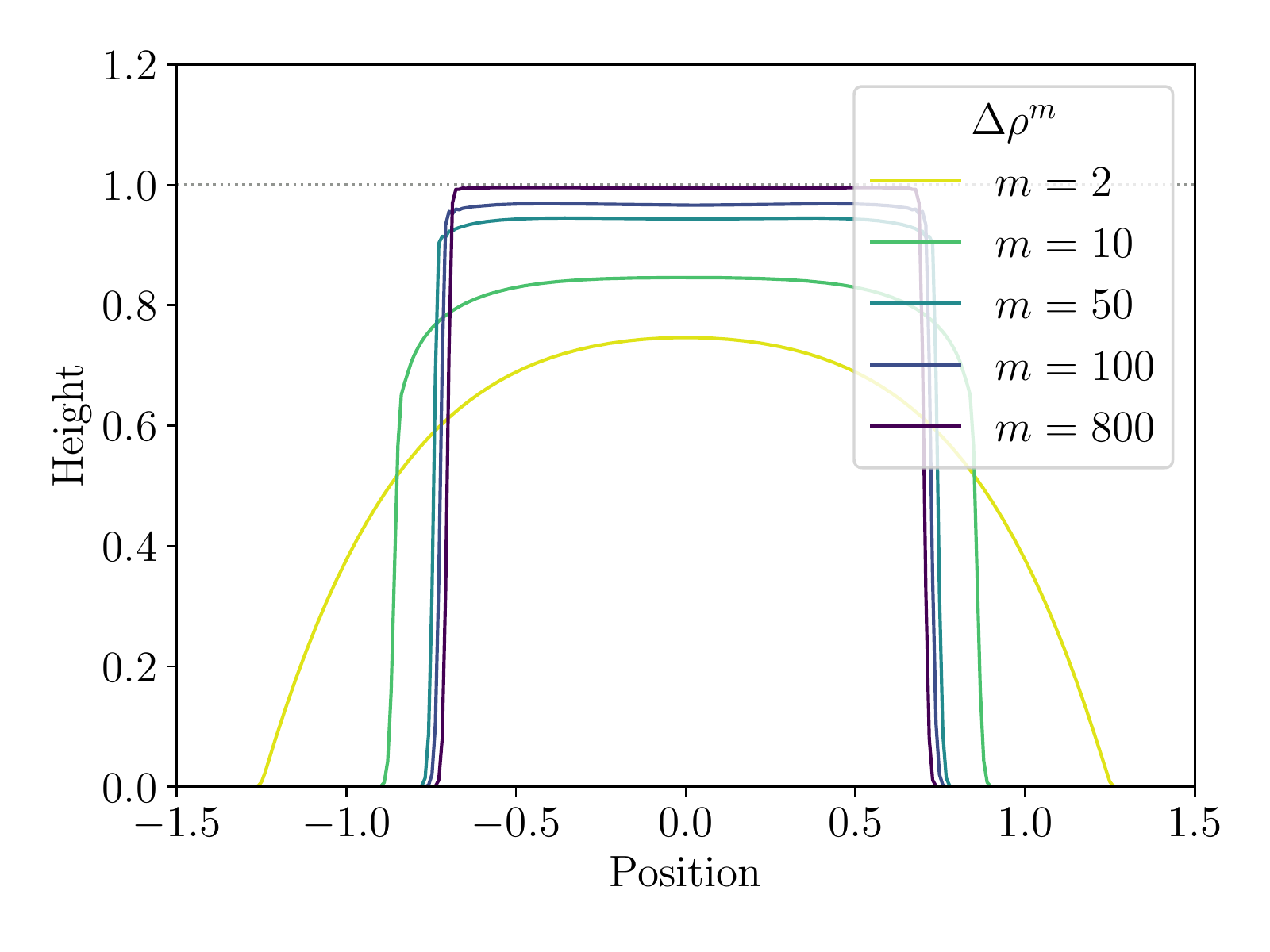}

\vspace{-.25cm}
\caption{We simulate asymptotic behavior of aggregation diffusion equations with repulsive attractive interaction kernel and for varying diffusion exponents. In all three examples, diffusion exponent $m=800$ is sufficient to impose the height constraint. The mass  of the solution determines whether the height constraint is active.}
\label{varyingdiffexp}
\end{center}
\vspace{-.75cm}
\end{figure}

In Figure \ref{varyingdiffexp}, we consider the repulsive attractive interaction kernel $W(x) = |x|^4/4 - |x|$ with degenerate diffusion, illustrating how the asymptotic behavior of the solution depends on the diffusion exponent $m$ for three different choices of mass of the initial data.  For all three choices of the mass, we observe that, once the diffusion exponent reaches $m = 800$, the degenerate diffusion acts effectively as a height constraint. In each case, the mass of the initial data determines whether the height constraint is active. Motivated by this result, in all future simulations, we approximate a height constraint via degenerate diffusion with diffusion exponent $m=800$.    We compute the equilibrium profile in Figure \ref{varyingdiffexp} by solving the aggregation diffusion equation with initial data given by a multiple of a Barenblatt profile $M\rho_2(x,0.15)$,  where  $M$ determines the mass. We compute the evolution up to time $T = 6.0$, with maximum time step size $k = 10^{-3}$. We discretize the domain with $N=1000$ particles for $m=2$ and $N = 500$ particles for $m >2$, and we regularize the diffusion and singular interaction terms with $\epsilon = 0.9$.

\begin{figure}[H]
\vspace{-.25cm}
\begin{center}
 {\bf Constrained Aggregation: Evolution to Equilibrium, $W(x) = |x|^4/4- |x|$}  \\ 
\vspace{.1cm}
 {\bf  mass = 0.6  } 
 
\includegraphics[height=4cm,trim={.6cm 0cm .6cm .8cm},clip]{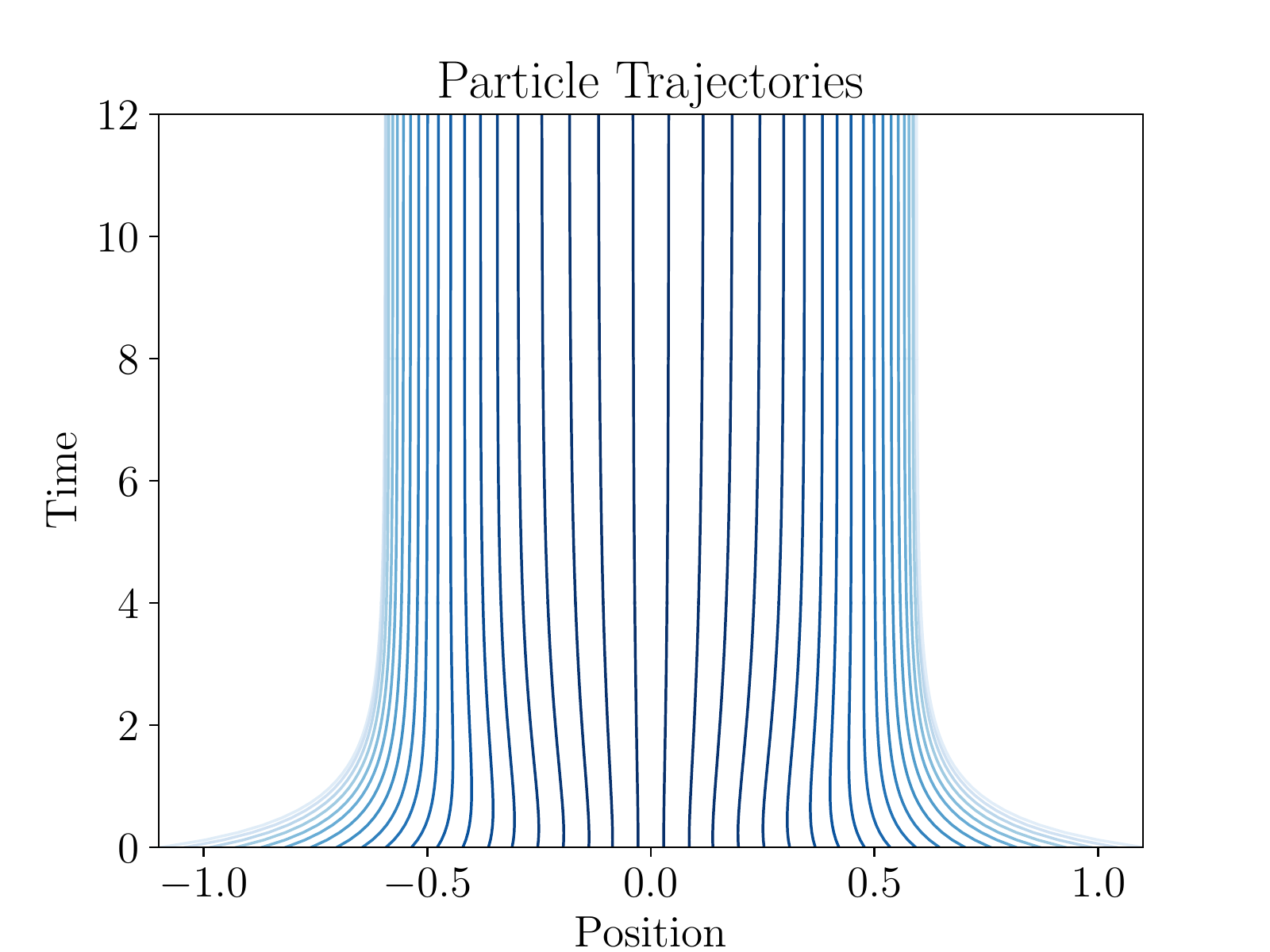} \includegraphics[height=4cm,trim={.6cm .7cm .6cm .6cm},clip]{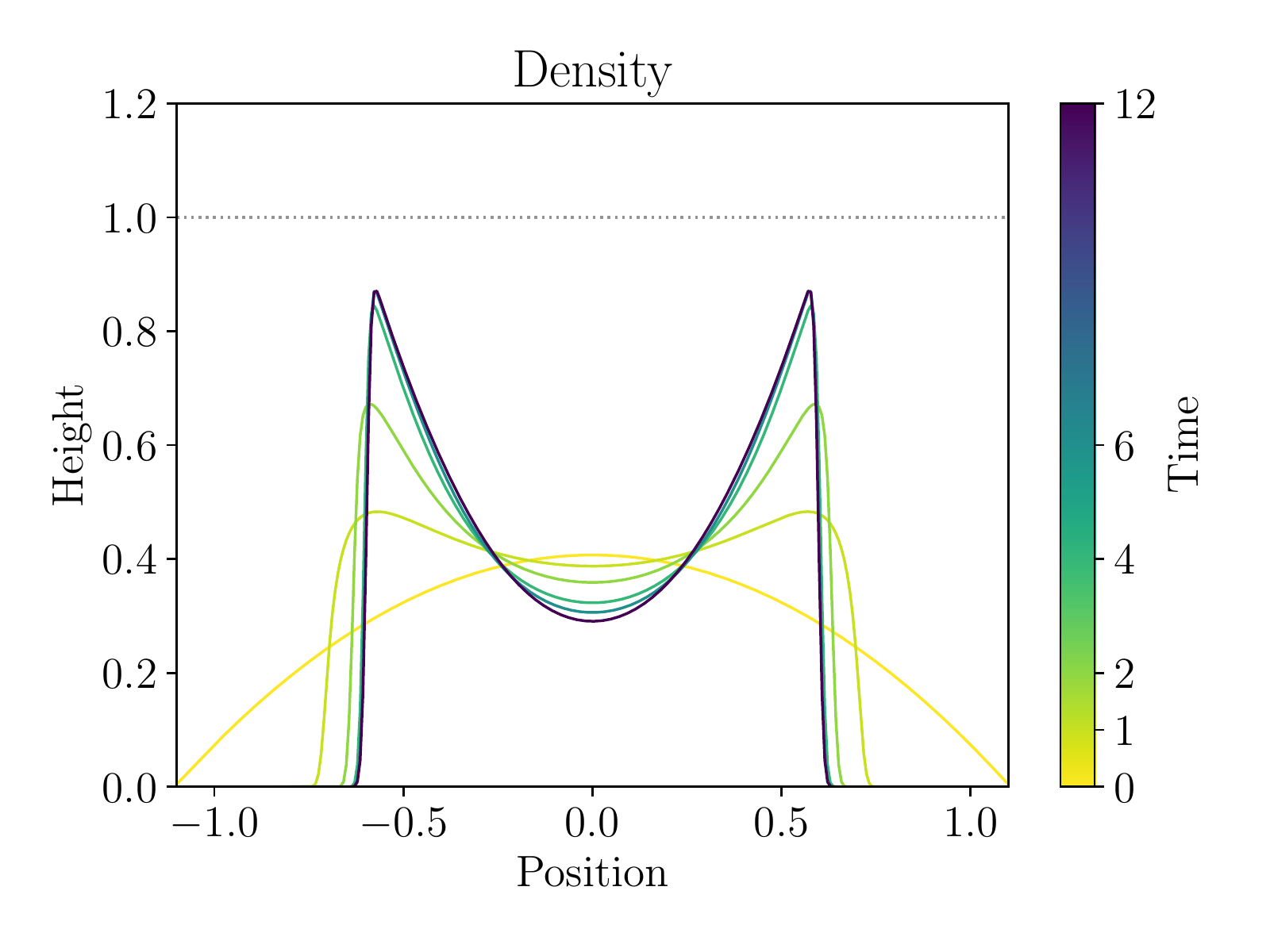}  \\

{\bf mass = 1.0}

\includegraphics[height=4cm,trim={.6cm 0cm .6cm .8cm},clip]{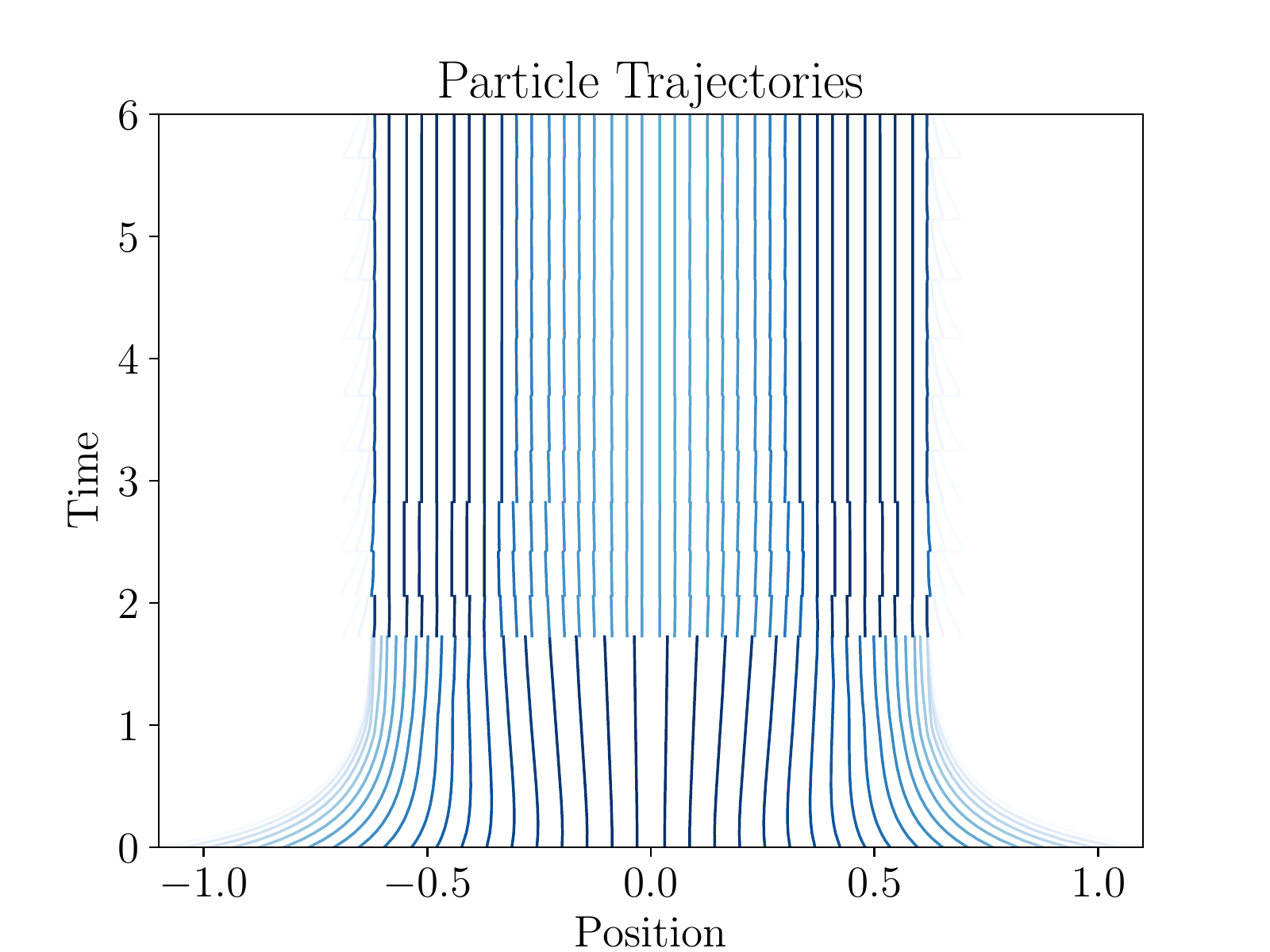} \includegraphics[height=4cm,trim={.6cm .7cm .6cm .6cm},clip]{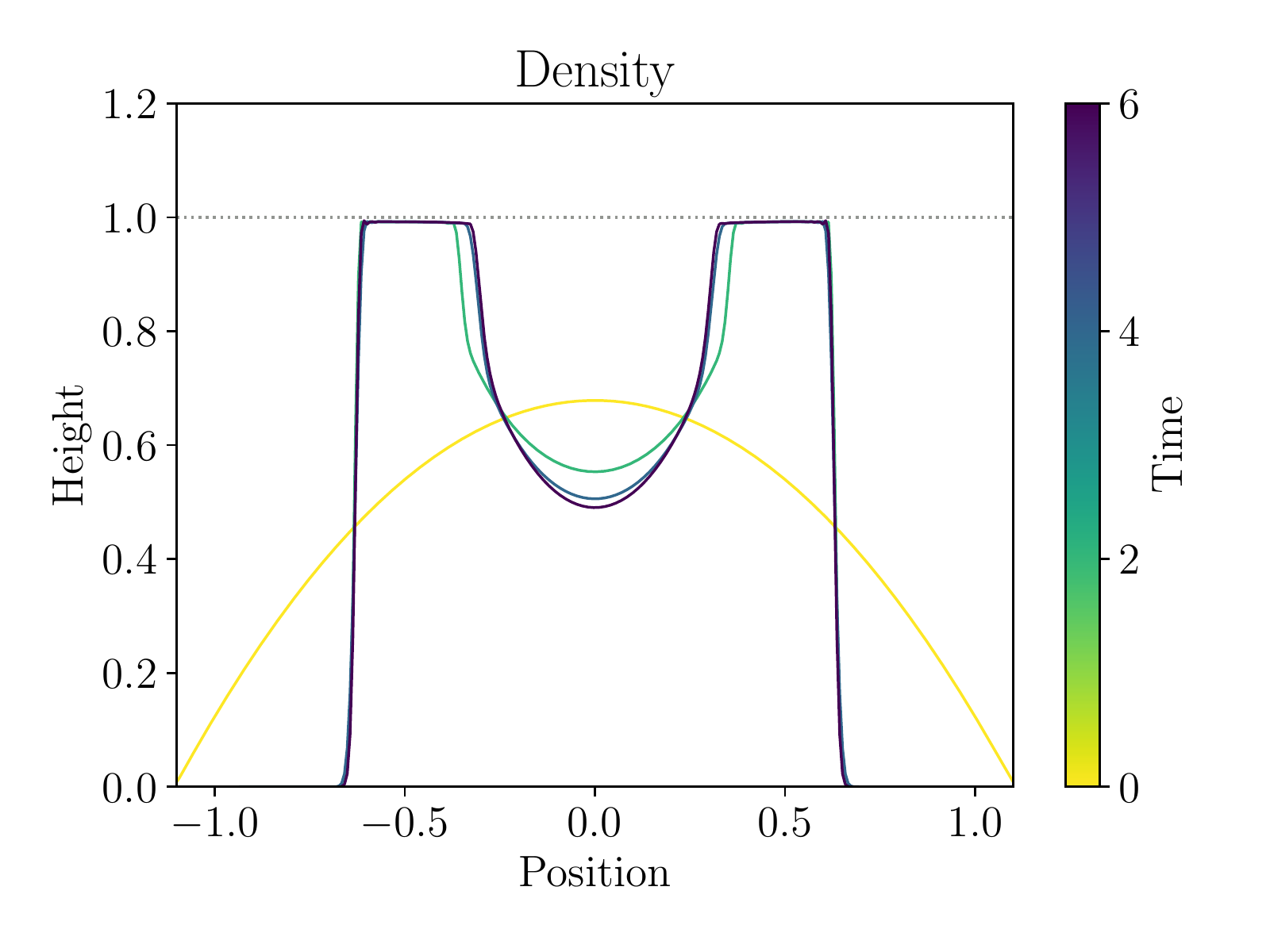}

{\bf mass = 1.4}

\includegraphics[height=4cm,trim={.4cm 0cm .6cm .8cm},clip]{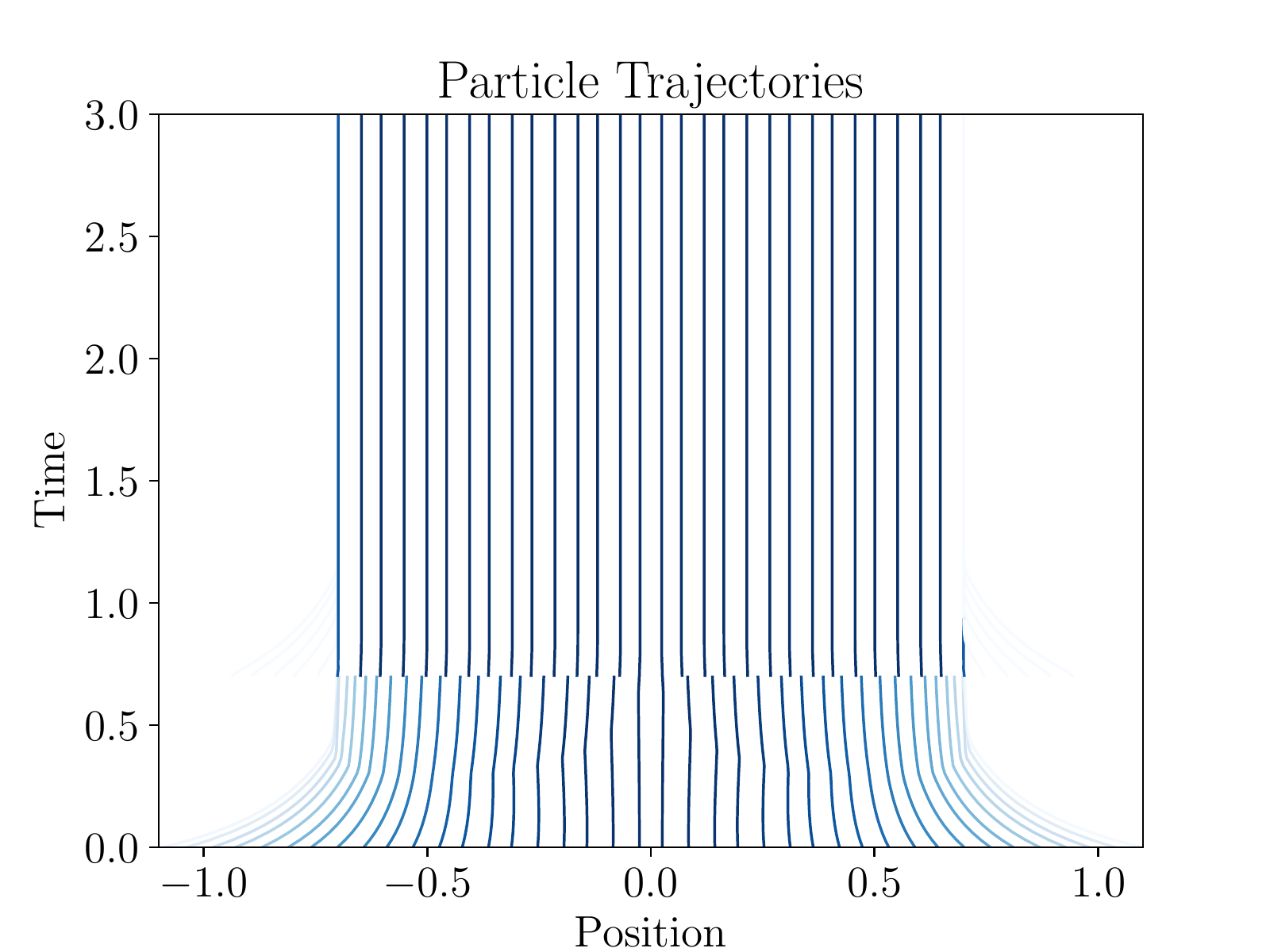} \includegraphics[height=4cm,trim={.6cm .7cm .6cm .6cm},clip]{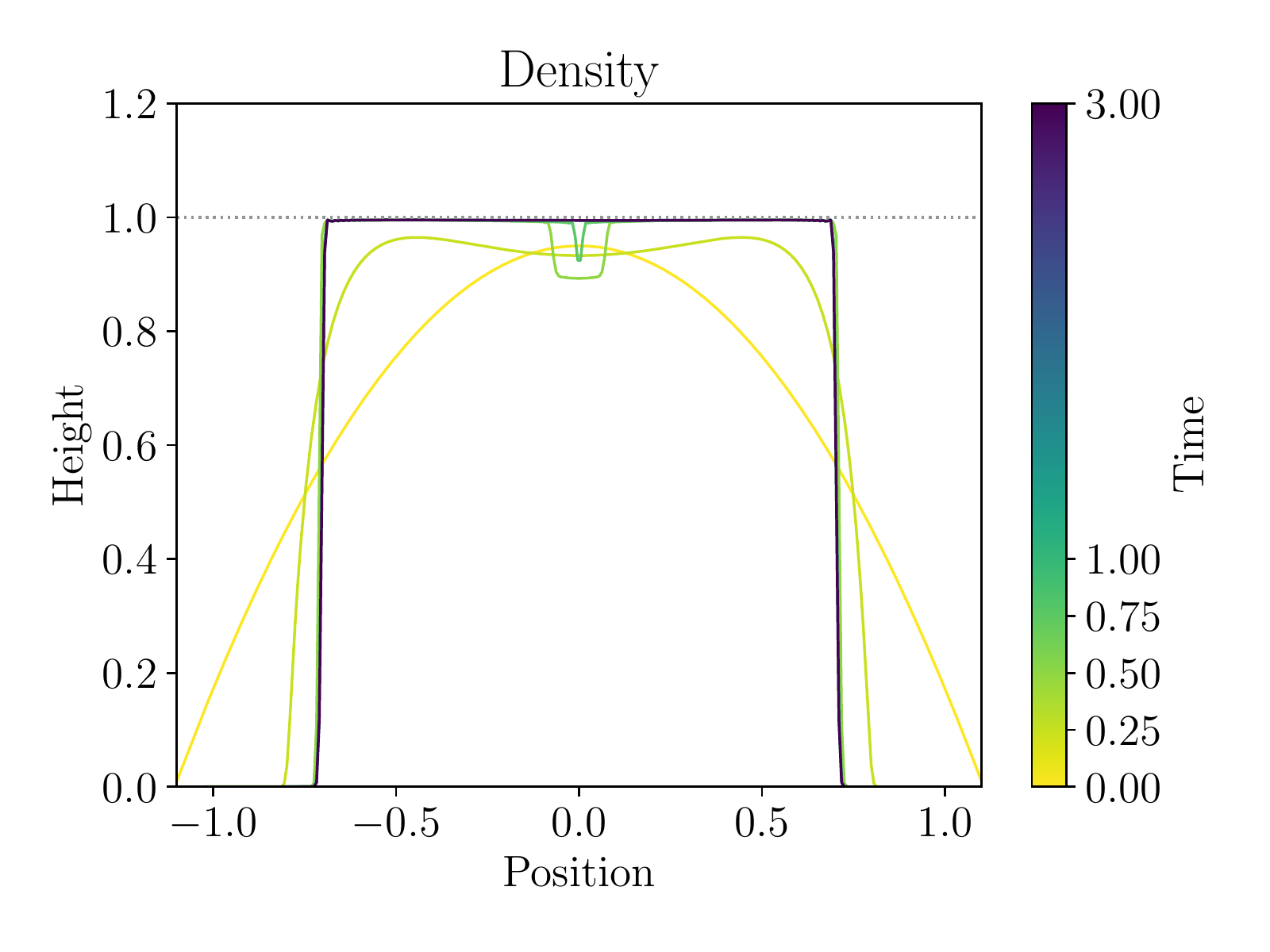}

\caption{We simulate evolution to equilibrium of height constrained aggregation equations with a repulsive attractive interaction kernel. The small mass equilibrium is in liquid phase, the large mass equilibrium is in solid phase (a characteristic function on a set), and in between the equilibrium is in an intermediate phase.}
\label{evolutiontoequilibriumx4mx2}
\end{center}
\vspace{-.75cm}
\end{figure}

\begin{figure}[h]
\begin{center}
 {\bf Constrained Aggregation: Equilibria for Varying Masses, $W(x) = |x|^4/4 - |x|^p/p$}  \\ 
\vspace{.1cm}

 {\bf p=0} \\ \hspace{-.55cm}
 \includegraphics[height=4cm,trim={.6cm .4cm .6cm .6cm},clip, valign=t]{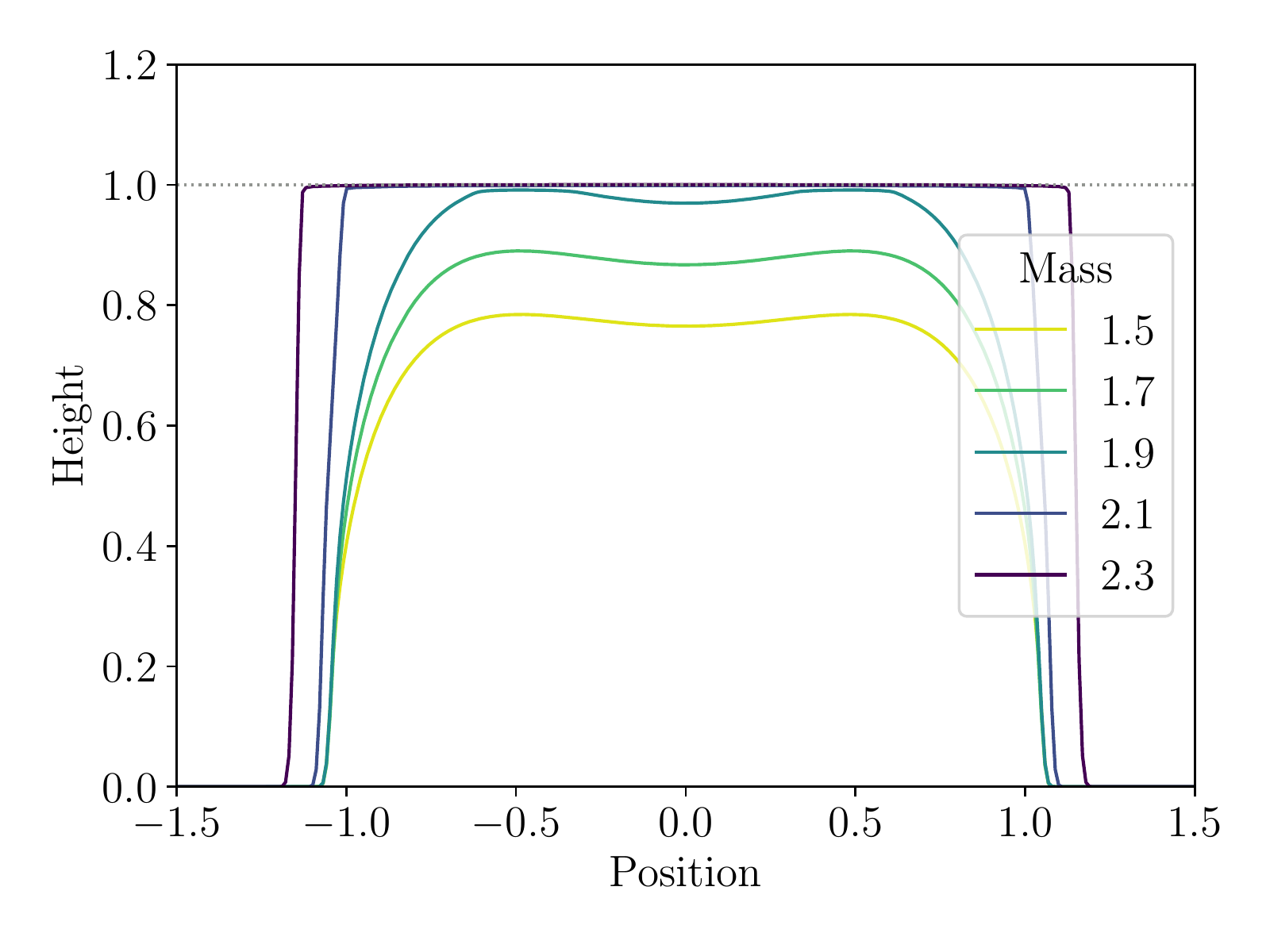} 
 
 {\bf p = 1} \\ \hspace{-.55cm}
\includegraphics[height=4cm,trim={.6cm .4cm .6cm .6cm},clip, valign=t]{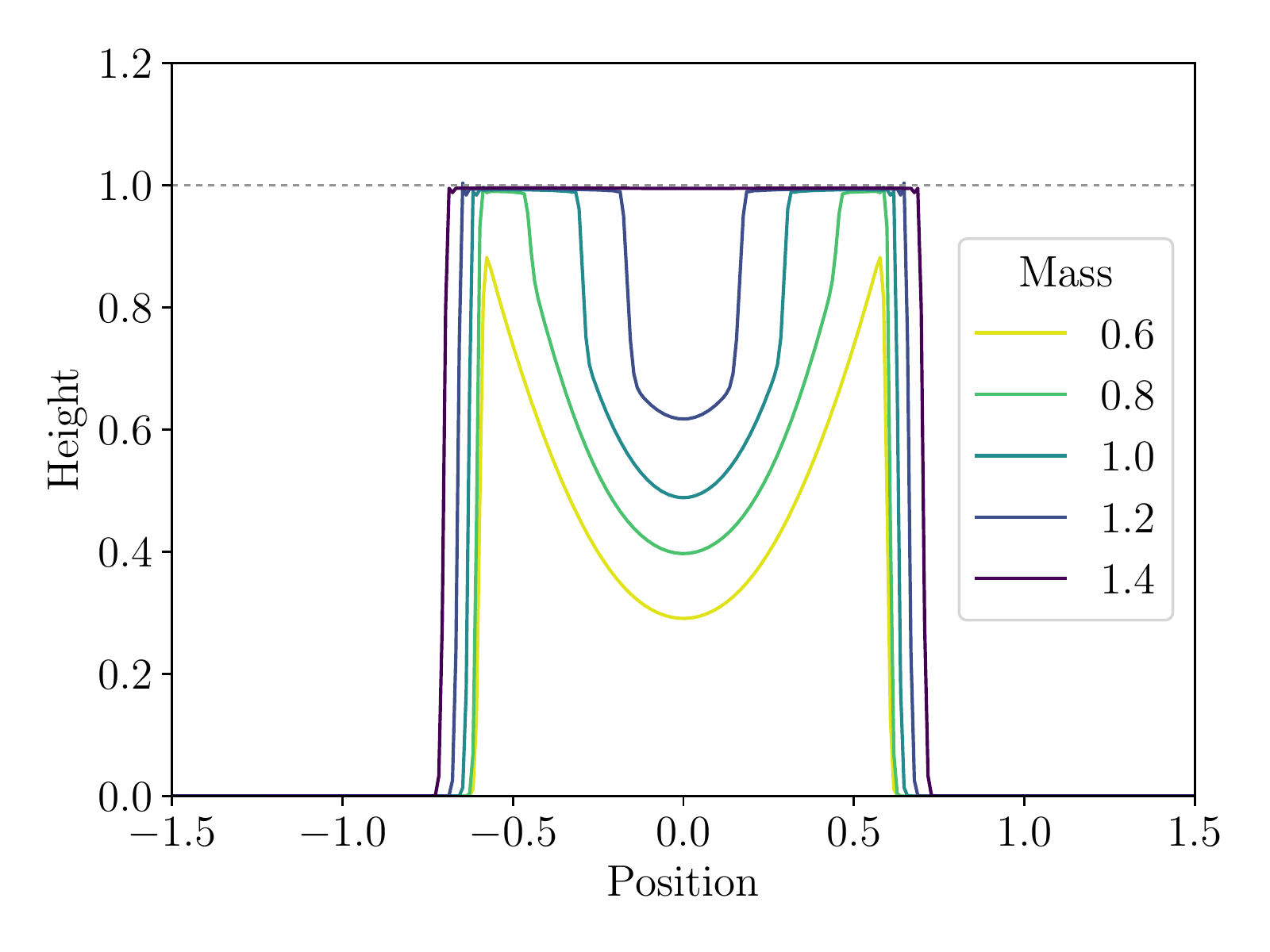} 

{\bf p=2} \\  \hspace{-.55cm}
\includegraphics[height=4cm,trim={.6cm .4cm .6cm .6cm},clip, valign=t]{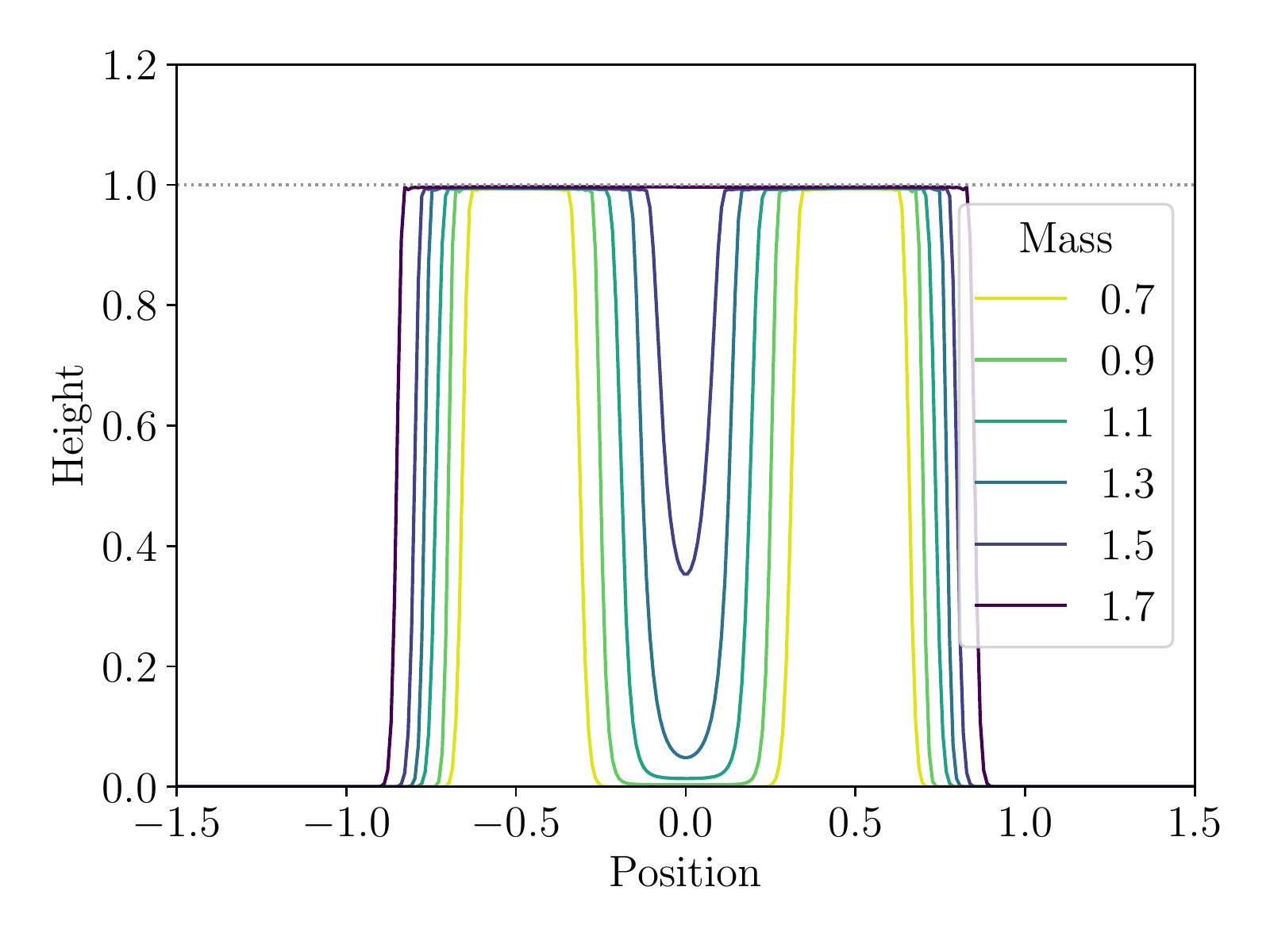}

\vspace{-.25cm}
\caption{We illustrate how the mass of the initial data affects the equilibrium configuration. For $p=0$ and $p=1$, as the mass increases, the equilibrium transitions from a liquid phase to an intermediate phase to a solid phase. However, for $p=2$, the equilibrium begins in solid phase, transitions to an intermediate phase, and then returns to solid phase. This provides numerical evidence that Lieb and Frank's result on the existence of liquid phase equilibria for singular interaction potentials is sharp, in the sense that it requires this strength of short range repulsion in order to hold \cite{FrankLieb}.}
\label{differentrepulsioncriticalmass}
\end{center}
\vspace{-.75cm}
\end{figure}

In Figure \ref{evolutiontoequilibriumx4mx2}, we illustrate the long time behavior of solutions to the height constrained aggregation equation for the same repulsive attractive interaction kernel and masses considered in  Figure \ref{varyingdiffexp}. To simulate the height constraint, we take $m=100$ for mass 0.6 and $m =800$ for masses 1.0 and 1.4.
In the left column, we plot the trajectories of particles computed via the blob method. In the right column, we plot the reconstructed density. For small mass 0.6, the equilibrium  is in a liquid phase  ($|\{ \rho_\infty = 1\}|=0$); for large mass 1.4, the equilibrium is in a solid phase ($|\{ \rho_\infty = 1 \}| = \int \rho_\infty$); and an intermediate phase exists for mass 1.0 $(0<|\{ \rho_\infty = 1\}|<\int \rho_\infty$). This agrees with previous numerical simulations of height constrained equilibria \cite{CraigTopaloglu}, as well as analytical results on the existence of set valued minimizers \cite{FrankLieb,BuChTo2016}. See  Figure \ref{differentrepulsioncriticalmass} for further analysis of  how equilibria of the height constrained problem depend on the mass, and  in particular how the strength of the short range repulsion affects this relationship.

\begin{figure}[h]
\vspace{-.25cm}
\begin{center}
 {\bf Constrained Aggregation: Evolution to Equilibrium, $W(x) = |x|$}  \\ 
\vspace{.1cm}
 
\adjustbox{valign=t}{\includegraphics[height=4cm,trim={.6cm .7cm 3.5cm .6cm},clip,valign = t]{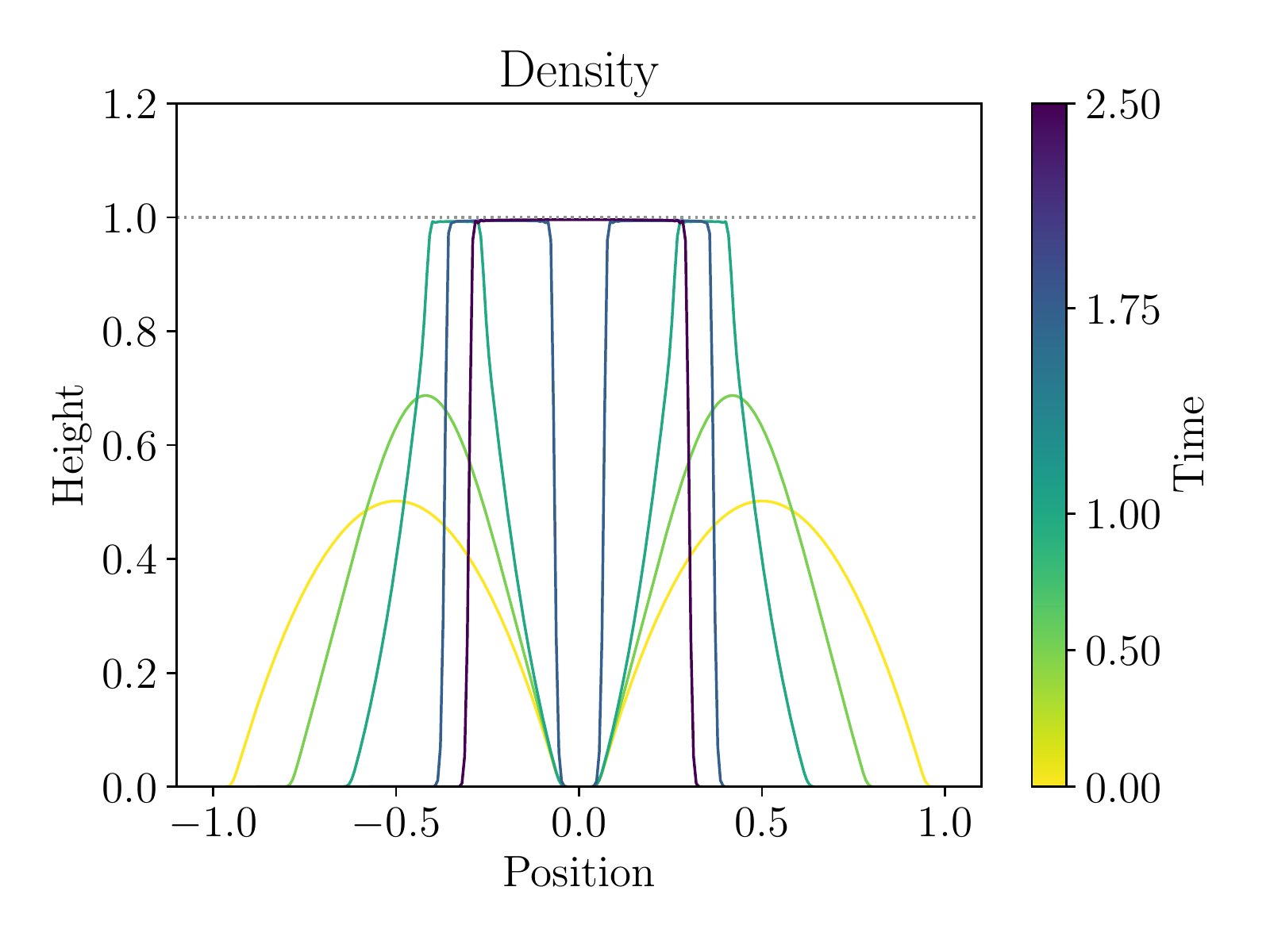}  \includegraphics[height=3.79cm,trim={2.2cm 1.3cm .6cm .58cm},clip,valign = t]{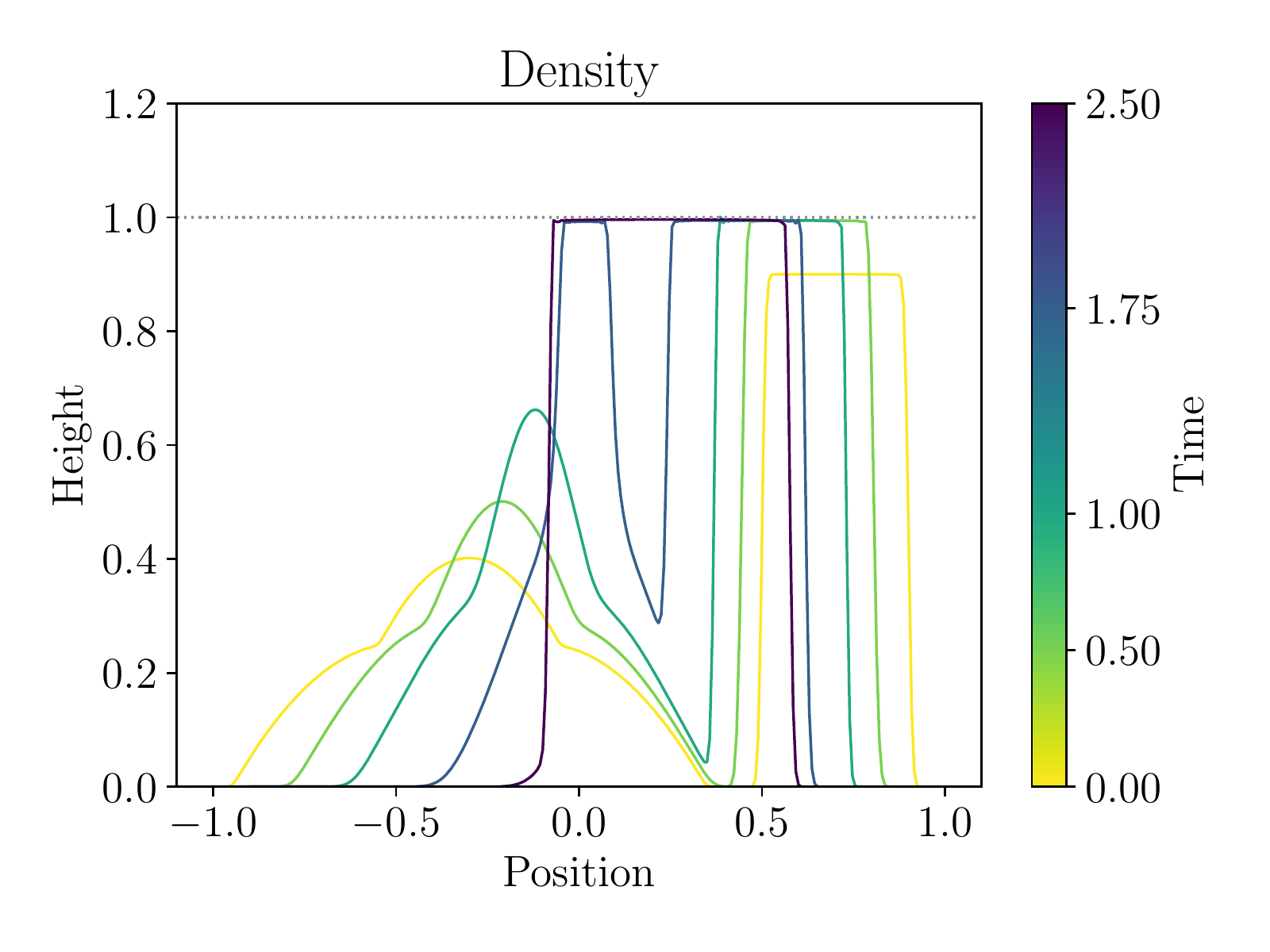}} \\
\vspace{-.25cm}

\caption{Evolution to equilibrium of height constrained aggregation with Newtonian interaction for two choices of initial data, illustrating different possible merging behavior as density reaches height constraint.}
\label{heightconstrNewt}
\end{center}
\vspace{-.75cm}
\end{figure}

In Figure \ref{heightconstrNewt}, we simulate the evolution to equilibrium of the height constrained aggregation equation with attractive Newtonian interaction. We observe different possible merging behaviors as the density reaches the height constraint. 
On the left, we consider initial data given by the sum of two Barenblatt profiles $(0.3)\rho_2(x,0.01)$ centered at $\pm 0.5$. In this case, the bumps first hit the height constraint and then merge. On the right, we consider initial data given by the sum of two Barenblatts $(0.15)\rho_2(x,0.01)$ centered at $x=0.1$ and $x=0.5$ and a characteristic function $(0.9)1_{[-0.2,0.2]}(x)$ centered at $x=0.7$. In this case, merging begins before the height constraint becomes active. We compute the evolution with a maximum time step size $k = 10^{-4}$, discretizing the domain with $N = 400$ particles and regularizing the diffusion and singular interaction terms with $\epsilon = 0.99$.

In Figure \ref{differentrepulsioncriticalmass}, we illustrate equilibria of the the height constrained aggregation equation for repulsive attractive interaction kernel $W(x) = |x|^4/4 - |x|^p/p$, showing how these equilibria depend on the repulsion exponent $p$ and the mass of the solution. For strong short range repulsion, as in the cases $p=0$ and $p=1$, we observe the following behavior as the mass increases:  the equilibrium is initially in liquid phase ($|\{ \rho_\infty = 1\}|=0$), transitions to an intermediate phase $(0<|\{ \rho_\infty = 1\}|<\int \rho_\infty$), and ultimately reaches a solid phase  ($|\{ \rho_\infty = 1 \}| = \int \rho_\infty$), where it is a characteristic function on an interval. However, for weaker short range repulsion, as in the case $p=2$, we observe that the equilibrium is in the solid phase for all sufficiently small masses, then transitions to an intermediate phase, and ultimately returns to a solid phase. This provides numerical evidence that Lieb and Frank's result on the existence of liquid phase equilibria for interaction kernels that are at least as singular as the Newtonian potential is sharp, in that it requires this strength of short range repulsion in order to hold \cite{FrankLieb}.

To compute the equilibrium profiles illustrated in Figure \ref{differentrepulsioncriticalmass}, we compute solutions of aggregation diffusion equations for $m = 800$ with patch initial data $\left(\frac{M}{2.1} \right)1_{[-1.05,1.05]}$, where the scalar $M$ determines the mass of the initial data. The maximum time for which we compute the evolution depends on the value of $p$ and the mass of the initial data, as described in Table \ref{maxtimetable}. For $p =0$, we take spatial discretization and regularization $N = 300$ and $\epsilon = 0.99$. For $p=1$, we take $N = 500$, $\epsilon = 0.9$. For $p=2$, we take $N = 600$, $\epsilon = 0.85$.

\begin{table}[H] 
\caption{Values of maximum time $T$ and minimum time step $k$ for each simulation in Figure \ref{differentrepulsioncriticalmass}.}
\begin{center}
\textbf{p = 0}  \hspace{1.3cm} \textbf{p = 1} \hspace{1.3cm} \textbf{p = 2} \\
\begin{tabular}{|c|c|c|}  
\hline
mass & T & k \\ \hline
1.5 & 1 & $10^{-3}$ \\
1.7 & 1 & $10^{-3}$  \\
1.9 & 1 & $10^{-3}$  \\
2.1 & 0.25 & $10^{-4}$ \\
2.3 &  0.25 & $10^{-4}$ \\  &\ & \\ \hline
\end{tabular}
\begin{tabular}{|c|c|c|}
\hline
mass & T & k \\ \hline
0.6 & 12 & $10^{-3}$ \\
0.8 & 6 & $10^{-3.5}$  \\
1.0 & 6 & $10^{-3.5}$  \\
1.2 & 3 & $10^{-3.5}$ \\
1.4  &  3 & $10^{-3.5}$ \\  & \ &  \\ \hline
\end{tabular}
\begin{tabular}{|c|c|c|}
\hline
mass & T & k \\ \hline
0.7 & 28 & $10^{-4}$ \\
0.9 & 20 & $10^{-4}$ \\
1.1 & 12 & $10^{-4}$ \\
1.3 & 12 & $10^{-4}$ \\
1.5 & 2.5 & $10^{-5}$ \\
1.7 & 0.75 & $10^{-5}$ \\ \hline
\end{tabular}
\end{center}
\label{maxtimetable}
\end{table}%

\subsection{Numerical examples: metastability} \label{metanumerics}
\vspace{0cm}
We now turn to numerical examples illustrating metastability behavior in singular limits, considering both localized aggregation and vanishing diffusion. See section \ref{singularmeta} for a discussion of the interplay between  singular limits and metastability.

\begin{figure}[H]
\vspace{0cm}
\begin{center}
 {\bf Initial Dynamics for Varying Interaction Range, $W(x) = e^{-x^2/4 \delta^2}/\sqrt{4 \pi \delta^2}$}  \\ 
\vspace{.1cm}

 {\bf $\boldsymbol\delta$ = 0.02} \\
  \includegraphics[height=4cm,trim={.6cm .7cm .6cm .6cm},clip, valign=t]{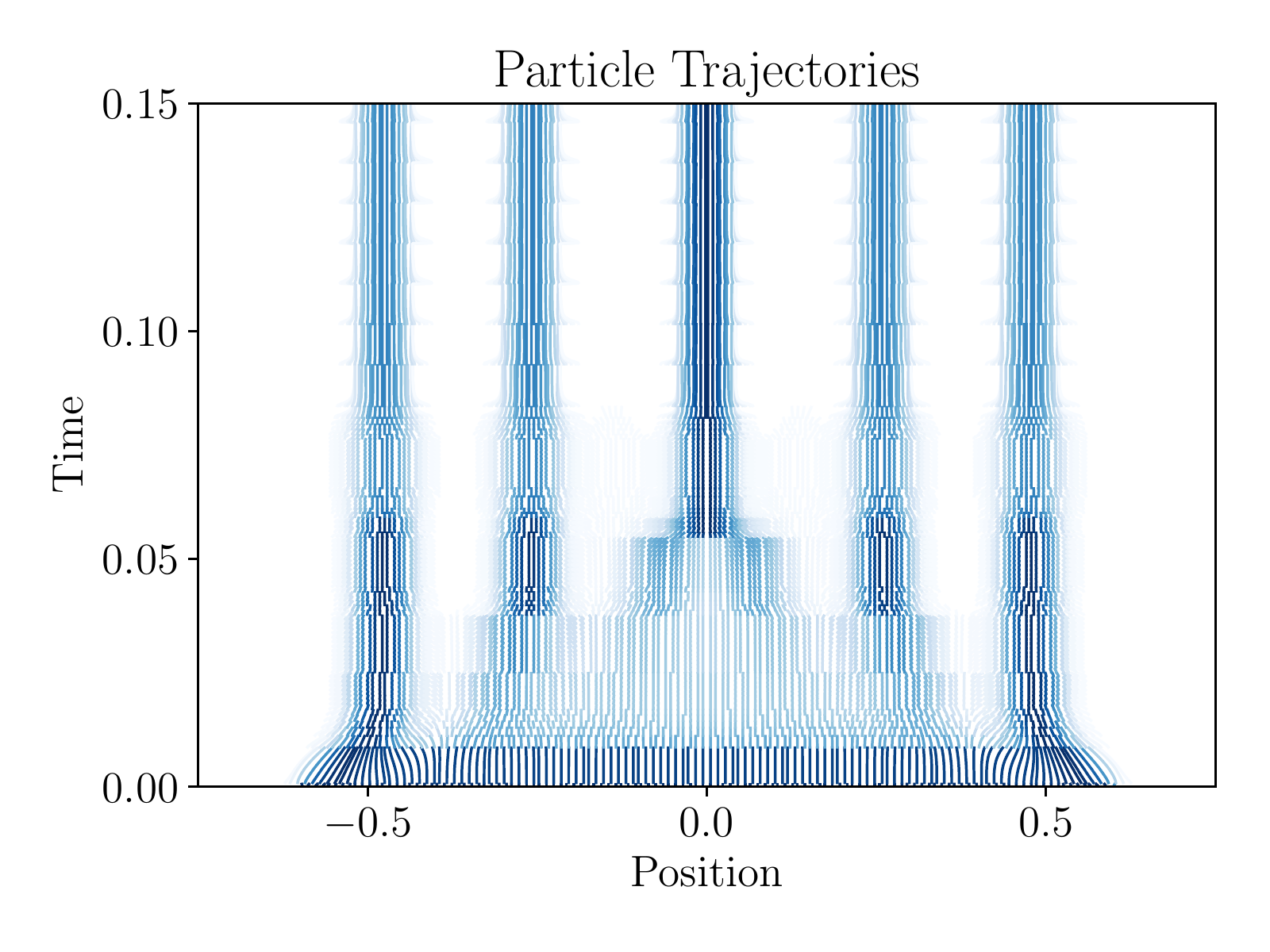}  
 \includegraphics[height=4cm,trim={.6cm .7cm .6cm .6cm},clip, valign=t]{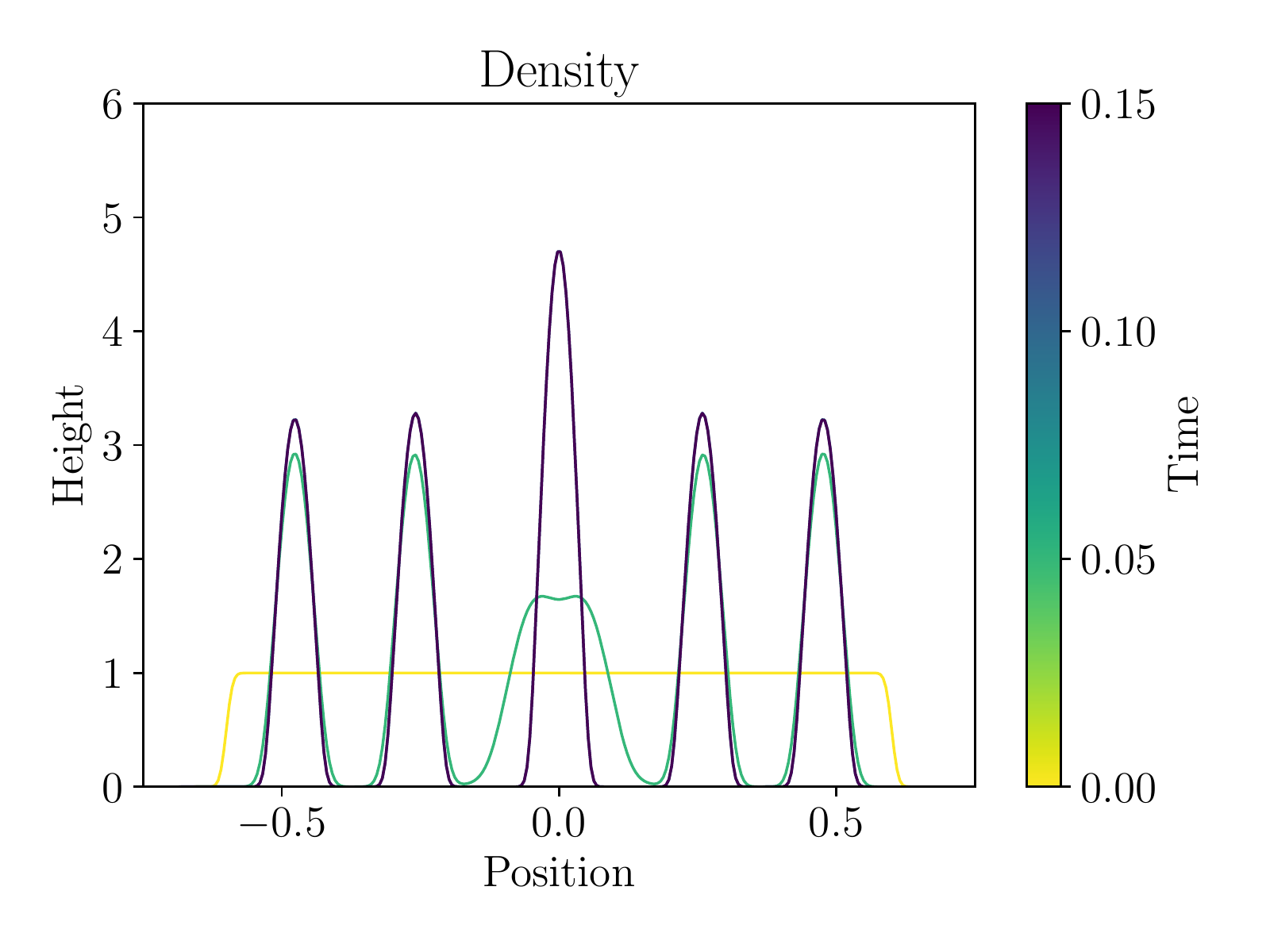} 

\vspace{-.25cm}
 {\bf $\boldsymbol\delta$ = 0.03} \\
 \includegraphics[height=4cm,trim={.6cm .7cm .6cm .6cm},clip, valign=t]{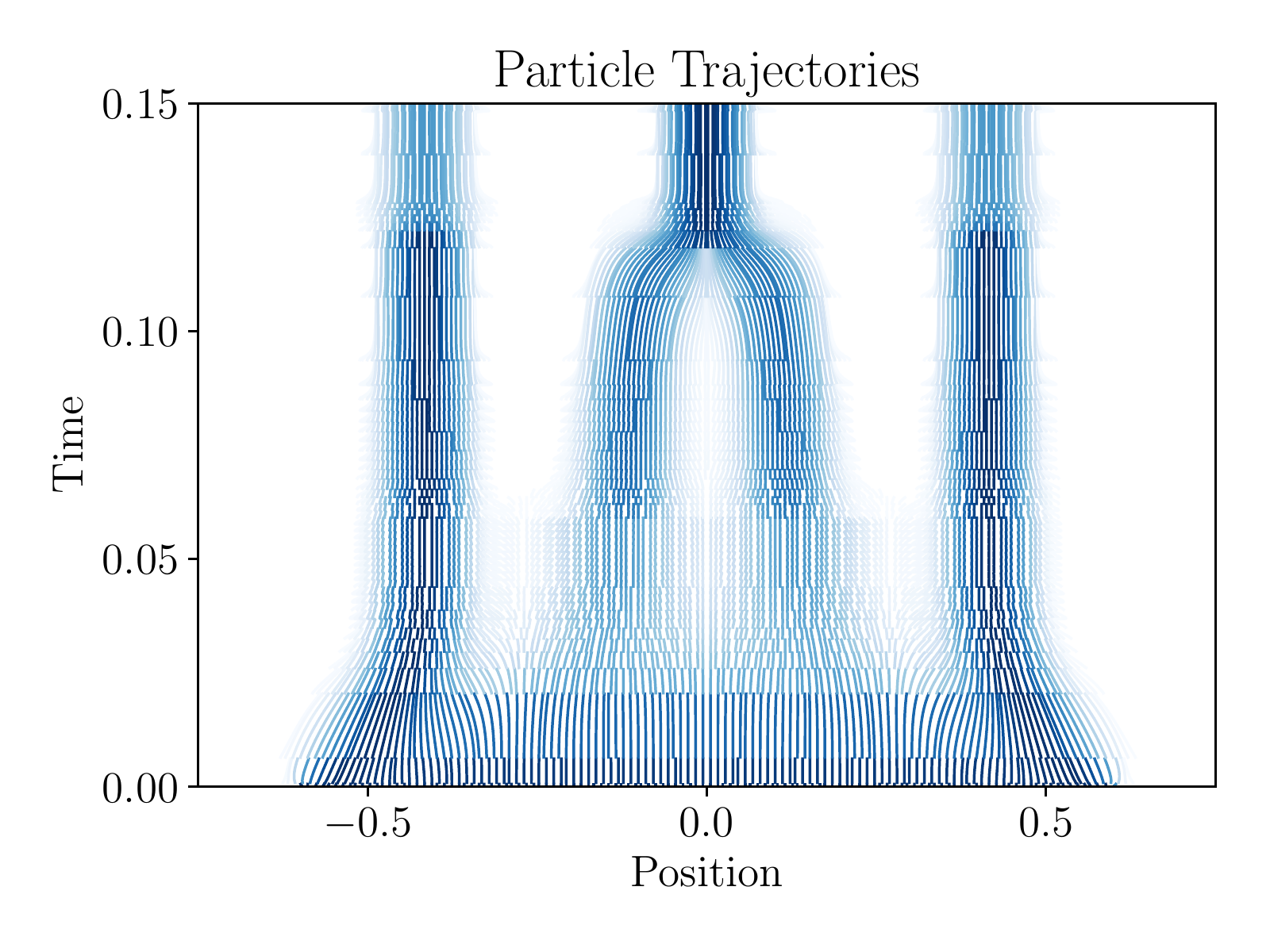}  \includegraphics[height=4cm,trim={.6cm .7cm .6cm .6cm},clip, valign=t]{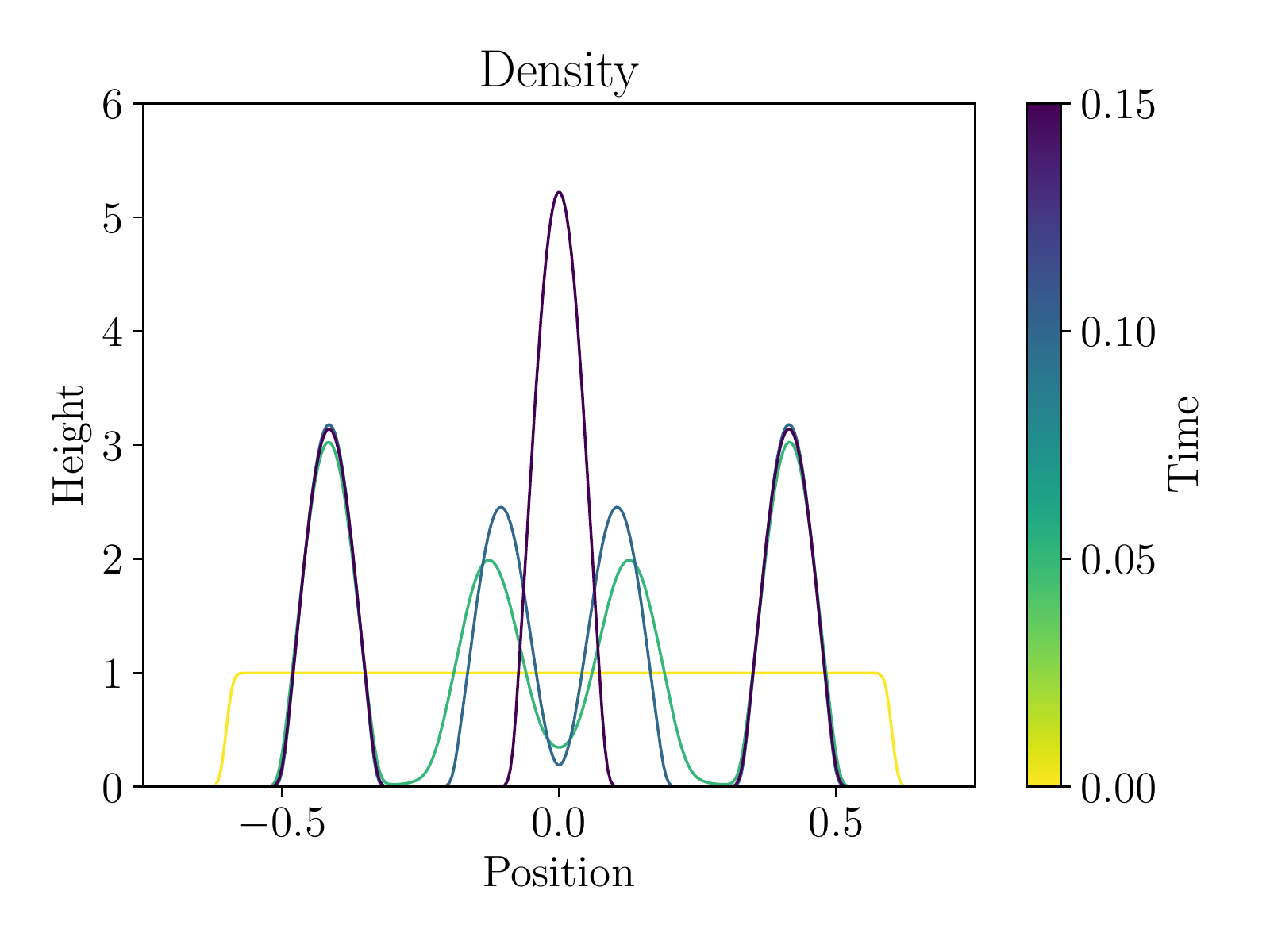} 

\vspace{-.25cm}
 {\bf $\boldsymbol\delta$ = 0.04} \\
 \includegraphics[height=4cm,trim={.6cm .7cm .6cm .6cm},clip, valign=t]{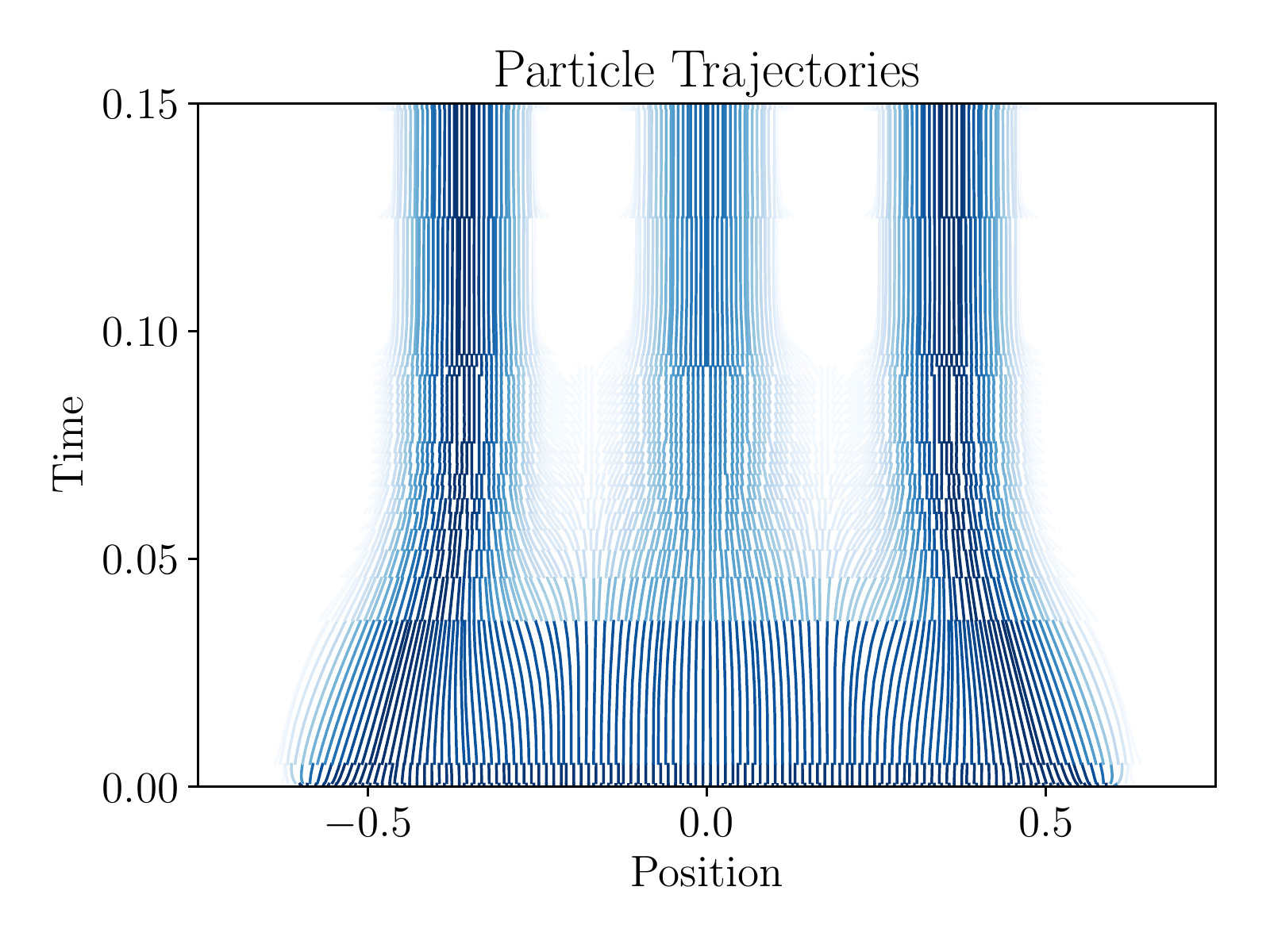}  \includegraphics[height=4cm,trim={.6cm .7cm .6cm .6cm},clip, valign=t]{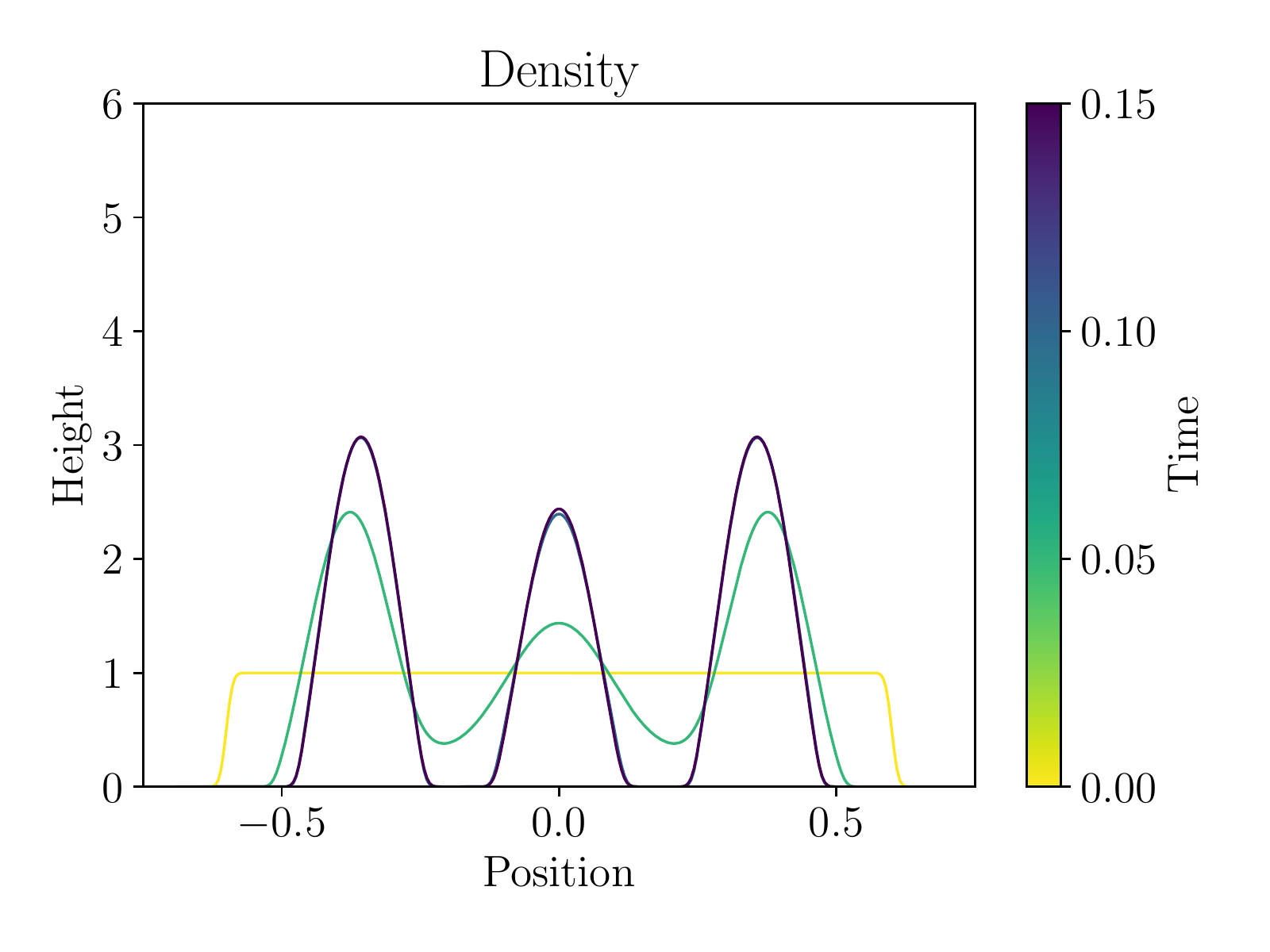} 

\vspace{-.25cm}
 {\bf $\boldsymbol\delta$ = 0.05} \\
  \includegraphics[height=4cm,trim={.6cm .7cm .6cm .6cm},clip, valign=t]{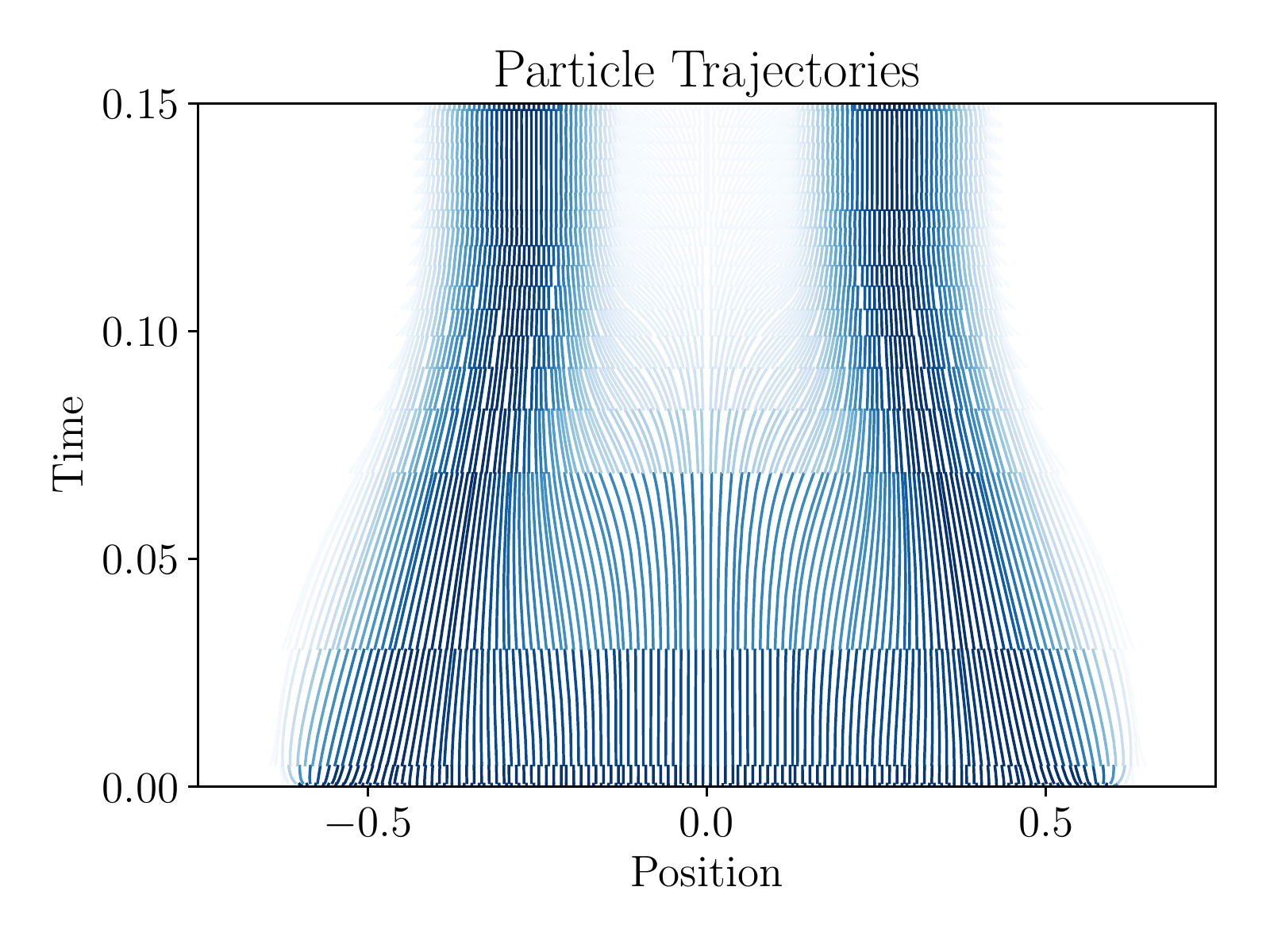}  
 \includegraphics[height=4cm,trim={.6cm .7cm .6cm .6cm},clip, valign=t]{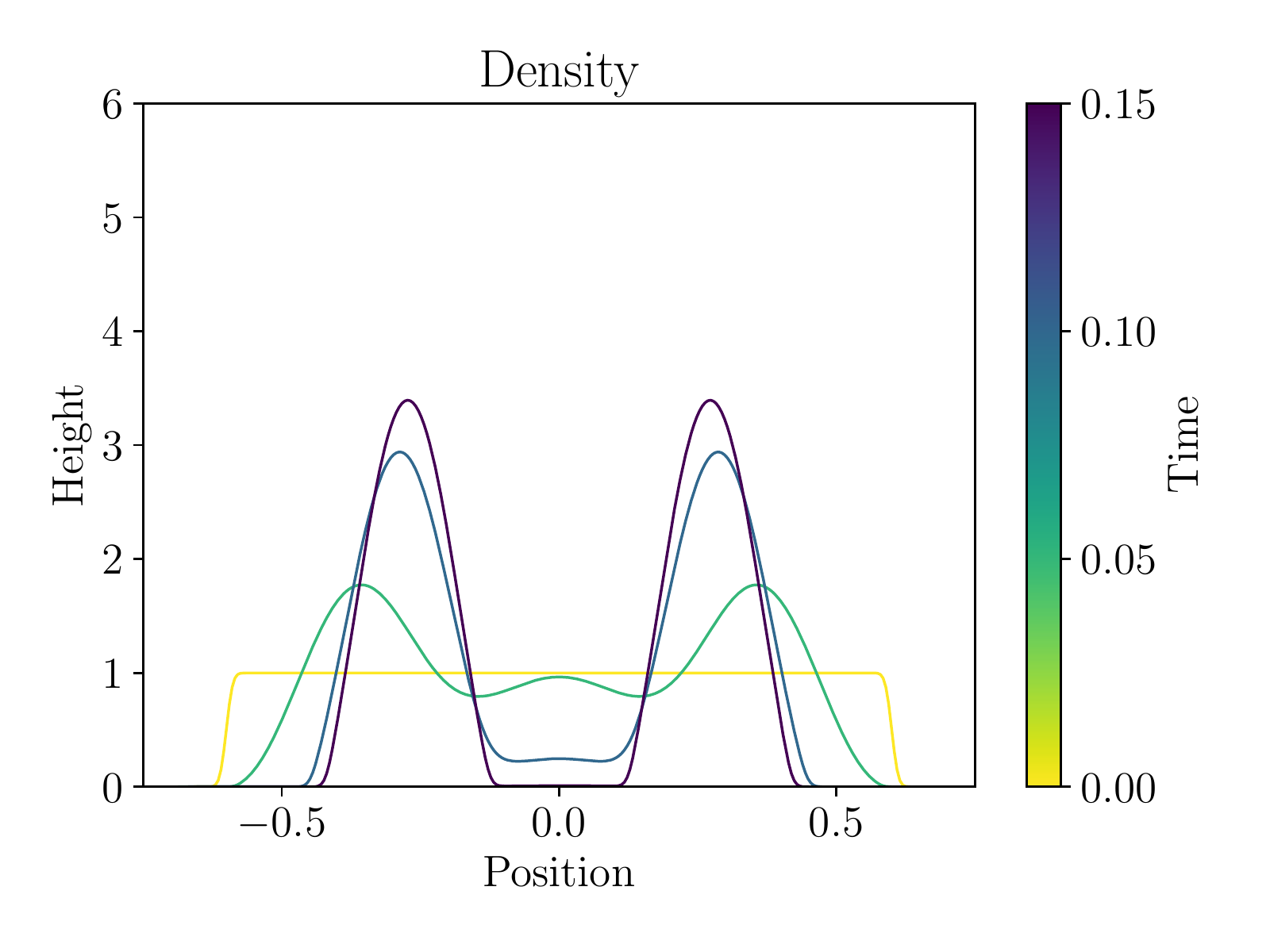} 

\vspace{-.25cm}
\caption{For aggregation diffusion equations with  porous medium diffusion ($m=2$), the localization of the interaction kernel $\delta$ affects the number of clusters in each metastable state.}
\label{varyinginteractionrange}
\end{center}
\vspace{-2cm}
\end{figure}

In Figure \ref{varyinginteractionrange}, we consider the initial dynamics of solutions for a localized,  attractive Gaussian interaction potential. The diffusion is of porous medium type $m=2$, weighted with diffusion coefficient $\nu = 0.25$, as in equation (\ref{agglocatt}). In each case, we take the initial data to be $1_{[-.6,.6]}$, with $N = 600$ particles, regularization $\epsilon = 0.9$, and maximal time step $k = 10^{-4}$. For $\delta = 0.02$, we briefly observe the formation of six bumps, before two quickly merge and form a metastable five bump configuration. Next, for $\delta = 0.03$, we observe the formation of a metastable state with four bumps, which at time $0.1$ quickly transitions into a second metastable state with three bumps. For $\delta = 0.04$, a three  bump metastable state emerges immediately, and for $\delta = 0.05$, the solution forms a two bump metastable state. Note that none of these numerical examples illustrates convergence to equilibrium, since in each case, the equilibrium must be radially decreasing, as discussed in section \ref{metaatt} \cite{CHVY}.

\begin{figure}[h]
\vspace{-.5cm}
\begin{center}
 {\bf Metastability for Varying Interaction Range, $W(x) = e^{-x^2/4 \delta^2}/\sqrt{4 \pi \delta^2}$}  \\ 
\vspace{.1cm}
 {\bf $\boldsymbol\delta$ = 0.04} \\
  \includegraphics[height=4cm,trim={.6cm .7cm .6cm .6cm},clip, valign=t]{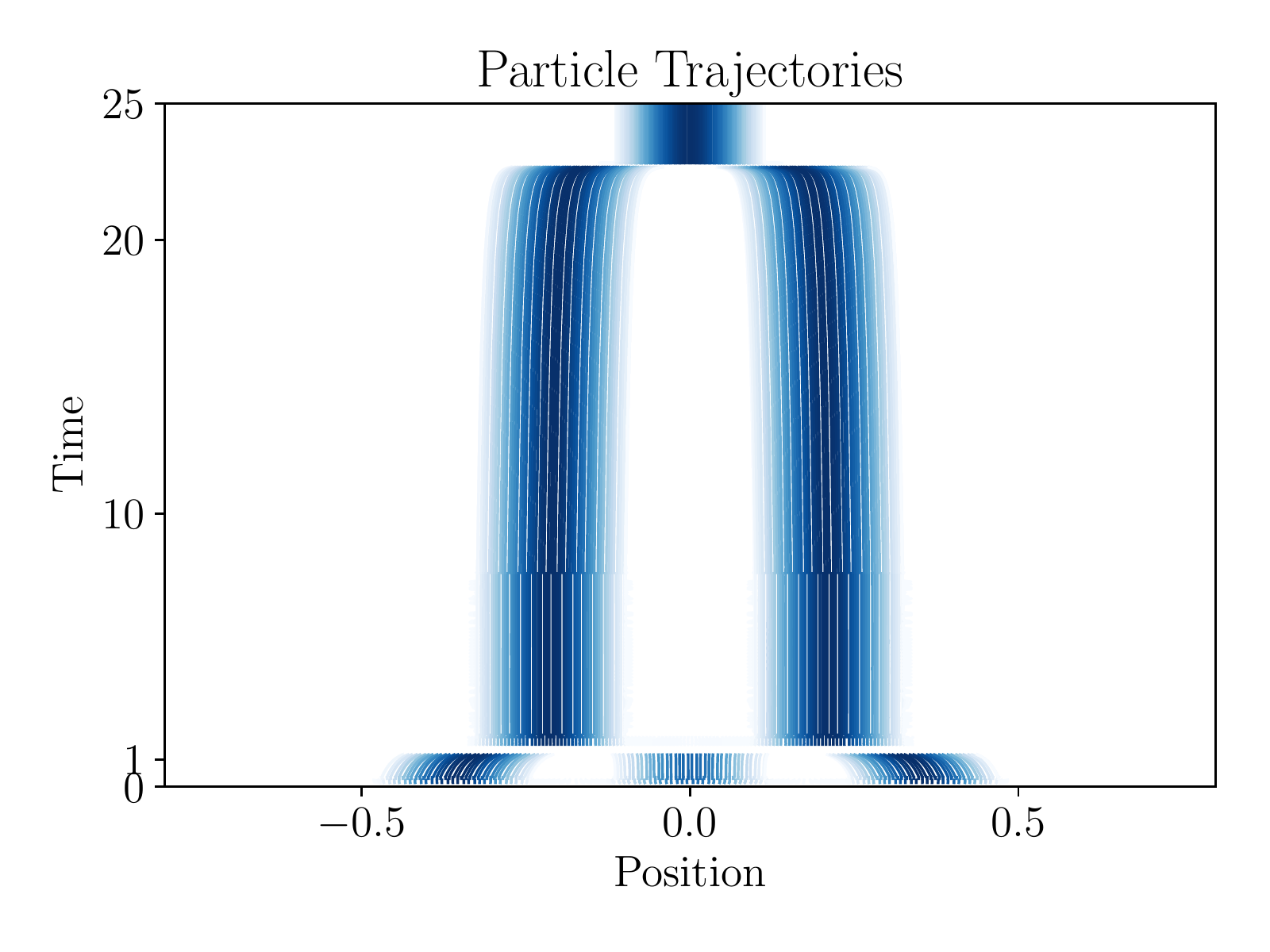}  
 \includegraphics[height=4cm,trim={.6cm .7cm .6cm .6cm},clip, valign=t]{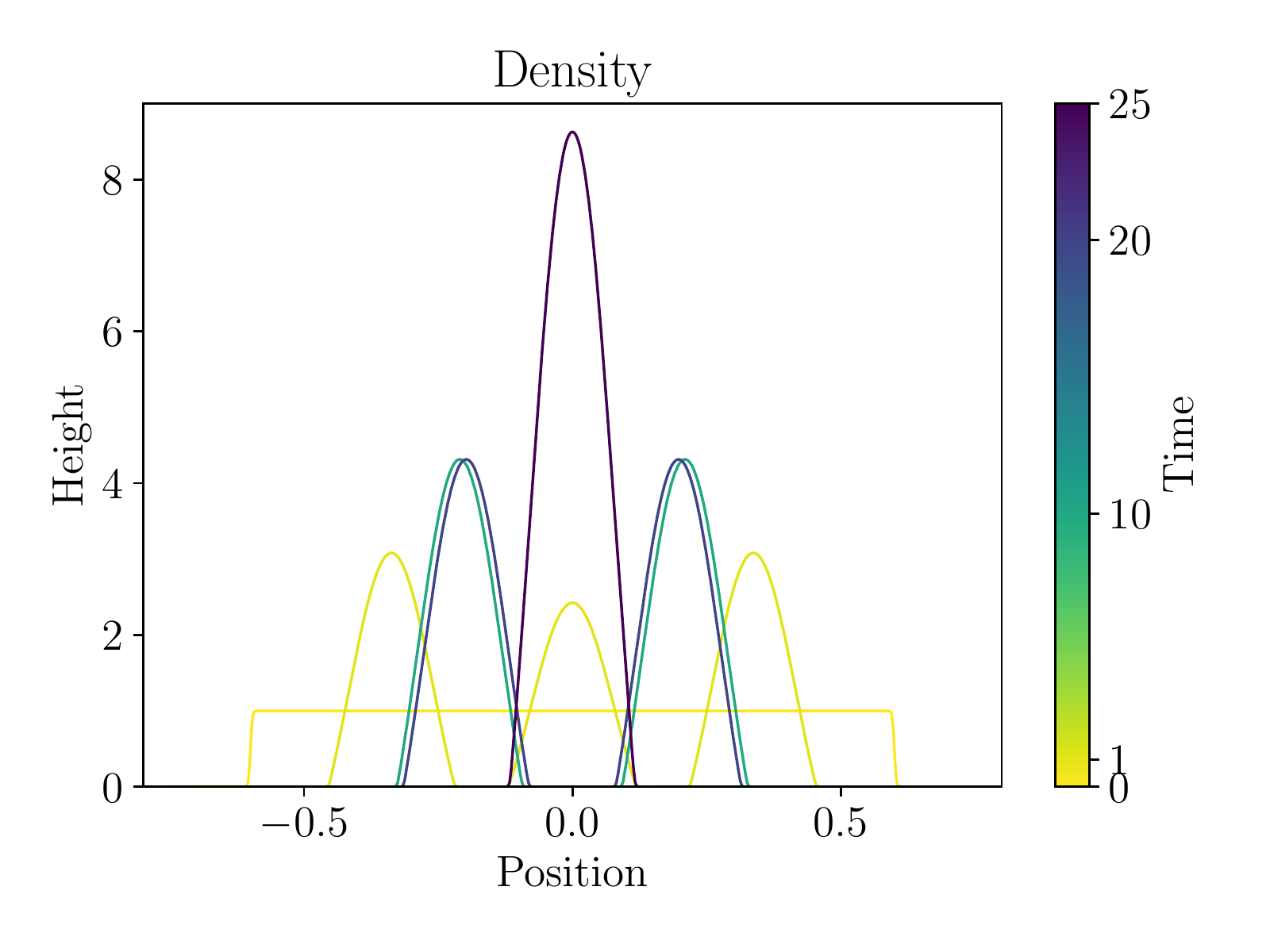} 

\vspace{-.2cm}
 {\bf $\boldsymbol\delta$ = 0.05} \\
  \includegraphics[height=4cm,trim={.6cm .7cm .6cm .6cm},clip, valign=t]{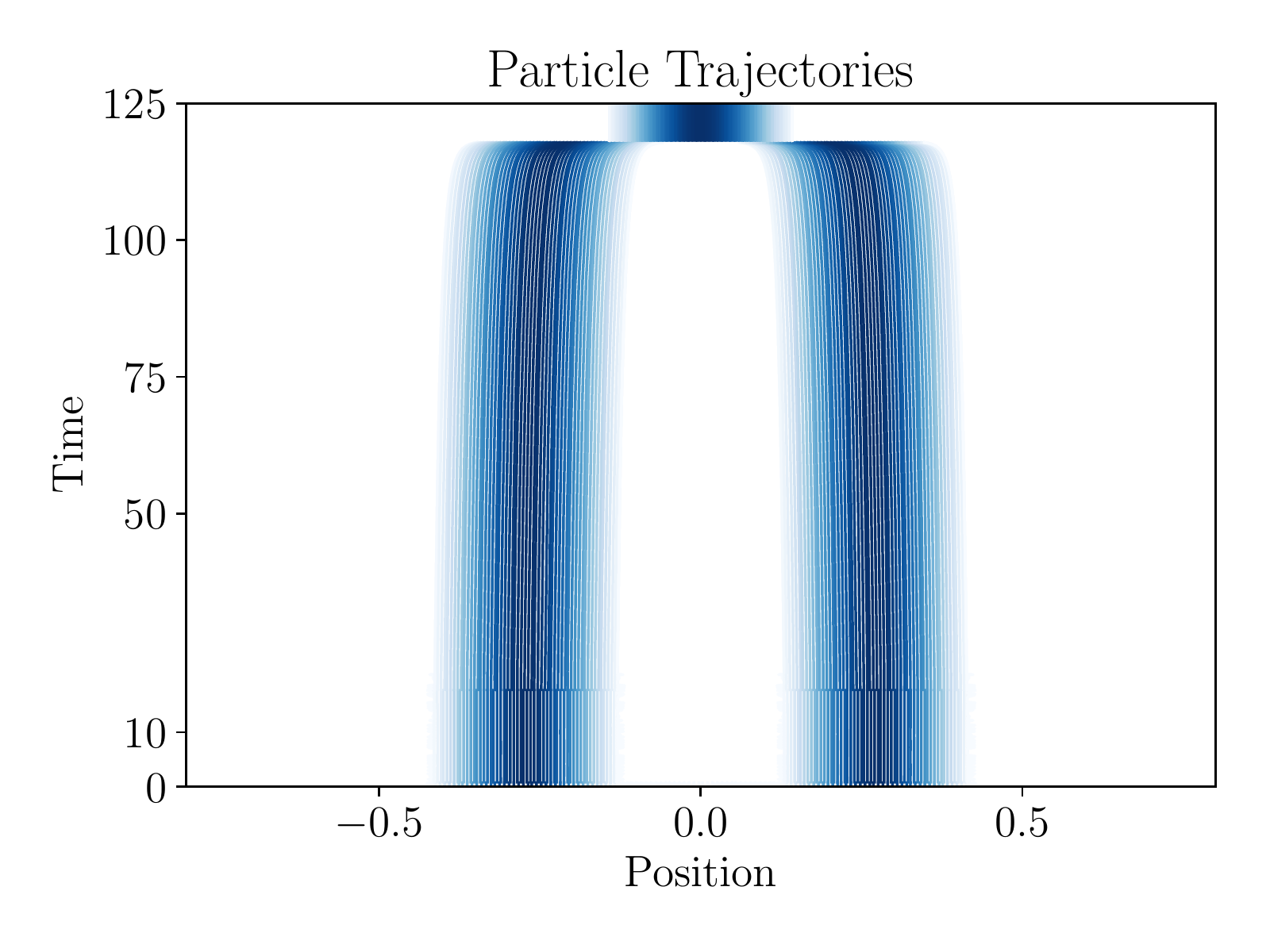}  
 \includegraphics[height=4cm,trim={.6cm .7cm .6cm .6cm},clip, valign=t]{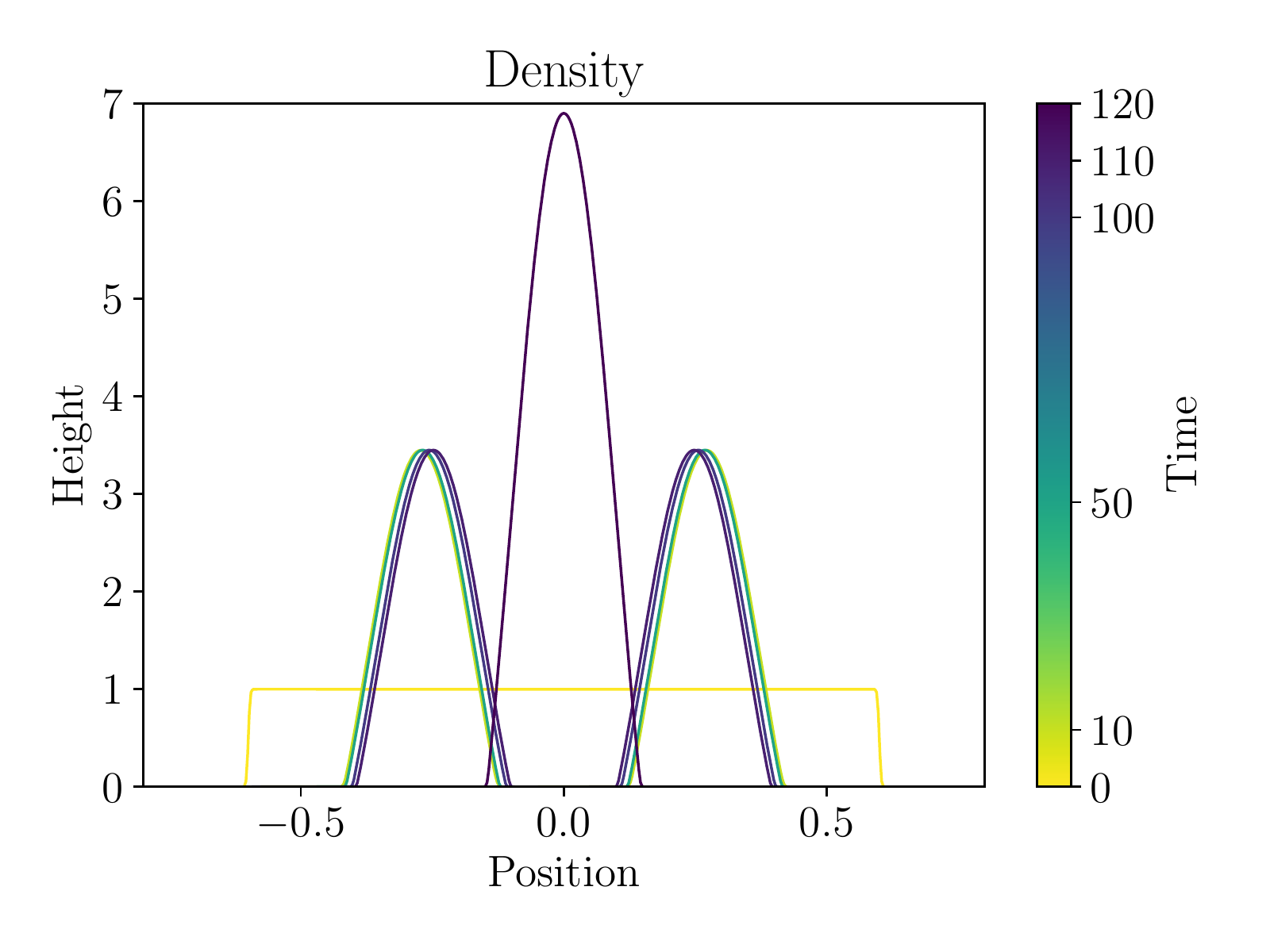} 
\vspace{-.25cm}
\caption{When a purely attractive interaction kernel is localized,  with porous medium diffusion ($m=2$), solutions of aggregation diffusion equations can form long-lived metastable states. In the case $\delta = 0.04$, the two bump metastable state lasts two orders of magnitude after the dynamics numerically appear to stabilize, and in the $\delta = 0.05$ case, the two bump metastable state lasts three orders of magnitude after apparent stabilization. Both solutions ultimately converge to a radially decreasing stationary state.}
\label{varyinginteractionrangelongtime}
\end{center}
\vspace{-.75cm}
\end{figure}

In Figure \ref{varyinginteractionrangelongtime}, we compute the convergence to equilibrium for two examples from the previous figure: $\delta = 0.04$ and $\delta = 0.05$. In the former case, the intermediate two bump metastable state is preserved on a timescale of order $\Delta t \sim 20$, before a near instantaneous transition to a single bump steady state. In the latter case, the two bump steady state, which Figure \ref{varyinginteractionrange} showed to form by time $t = 0.15$, is preserved past time $t = 120$ before rapidly transitioning to a single bump equilibrium. The long lived metastable state in the $\delta = 0.05$ case can be easily confused numerically with a steady state. Indeed, it was only the theoretical result ensuring that the equilibrium configuration must be radially decreasing that prompted us to run the simulations three orders of magnitude after the dynamics seem to stabilize ($t \in [0.15, 125.00]$).

\begin{figure}[H]
\begin{center}
\vspace{-.5cm}
 {\bf Metastability for Varying Diffusion Coefficient, $W(x) = e^{-x^2/4 (.05)^2}/\sqrt{4 \pi (.05)^2}$}  \\ 
\vspace{.1cm}

 {\bf $\boldsymbol\nu$ = 0.15} \\

 \includegraphics[height=4cm,trim={.6cm .7cm .6cm .6cm},clip, valign=t]{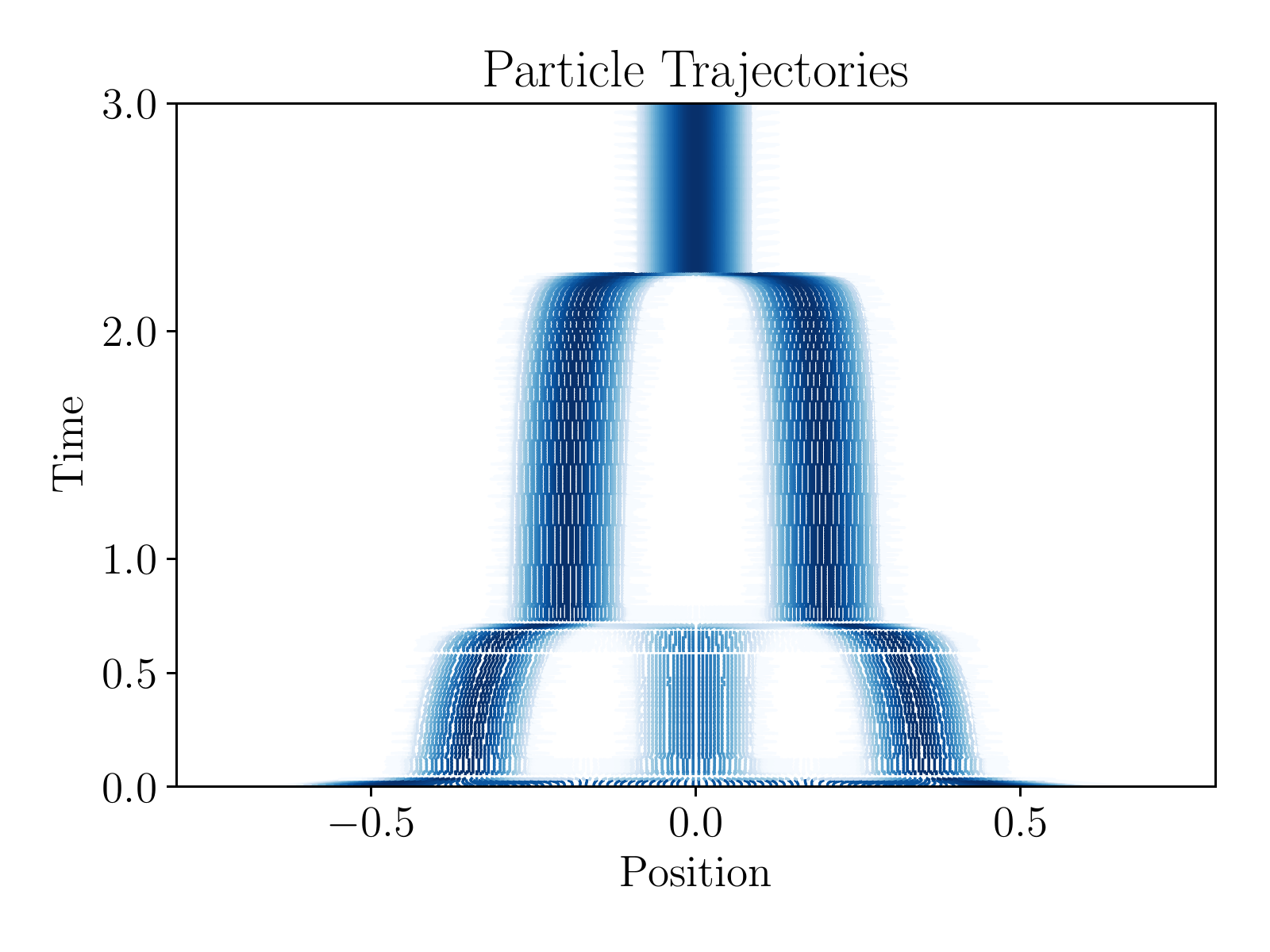} 
  \includegraphics[height=4cm,trim={.6cm .7cm .6cm .6cm},clip, valign=t]{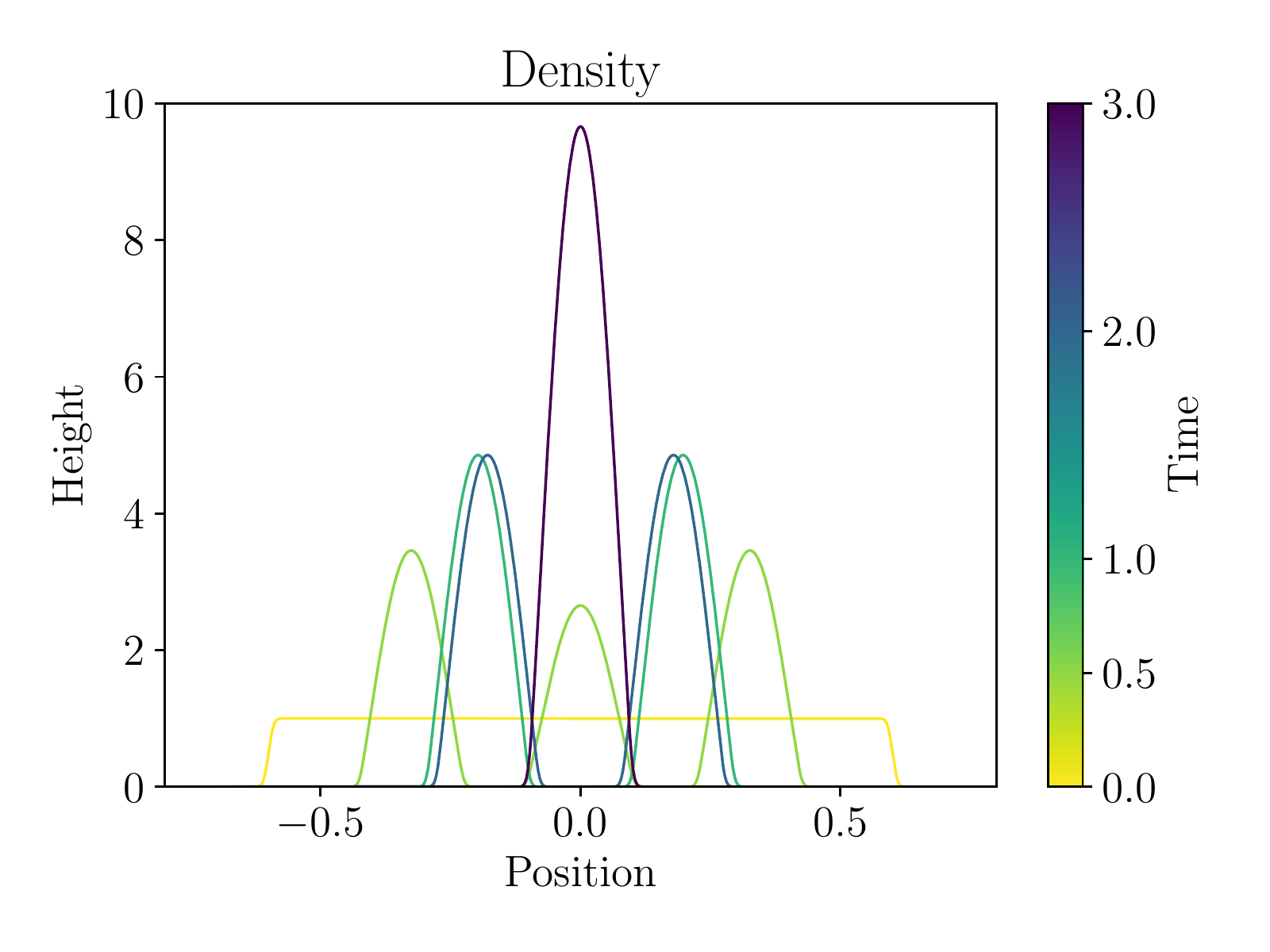}  

\vspace{-.2cm}
 {\bf $\boldsymbol\nu$ = 0.25} \\

  \includegraphics[height=4cm,trim={.6cm .7cm .6cm .6cm},clip, valign=t]{DiffMeta_gaussian_p05_particles_longtime.pdf}  
 \includegraphics[height=4cm,trim={.6cm .7cm .6cm .6cm},clip, valign=t]{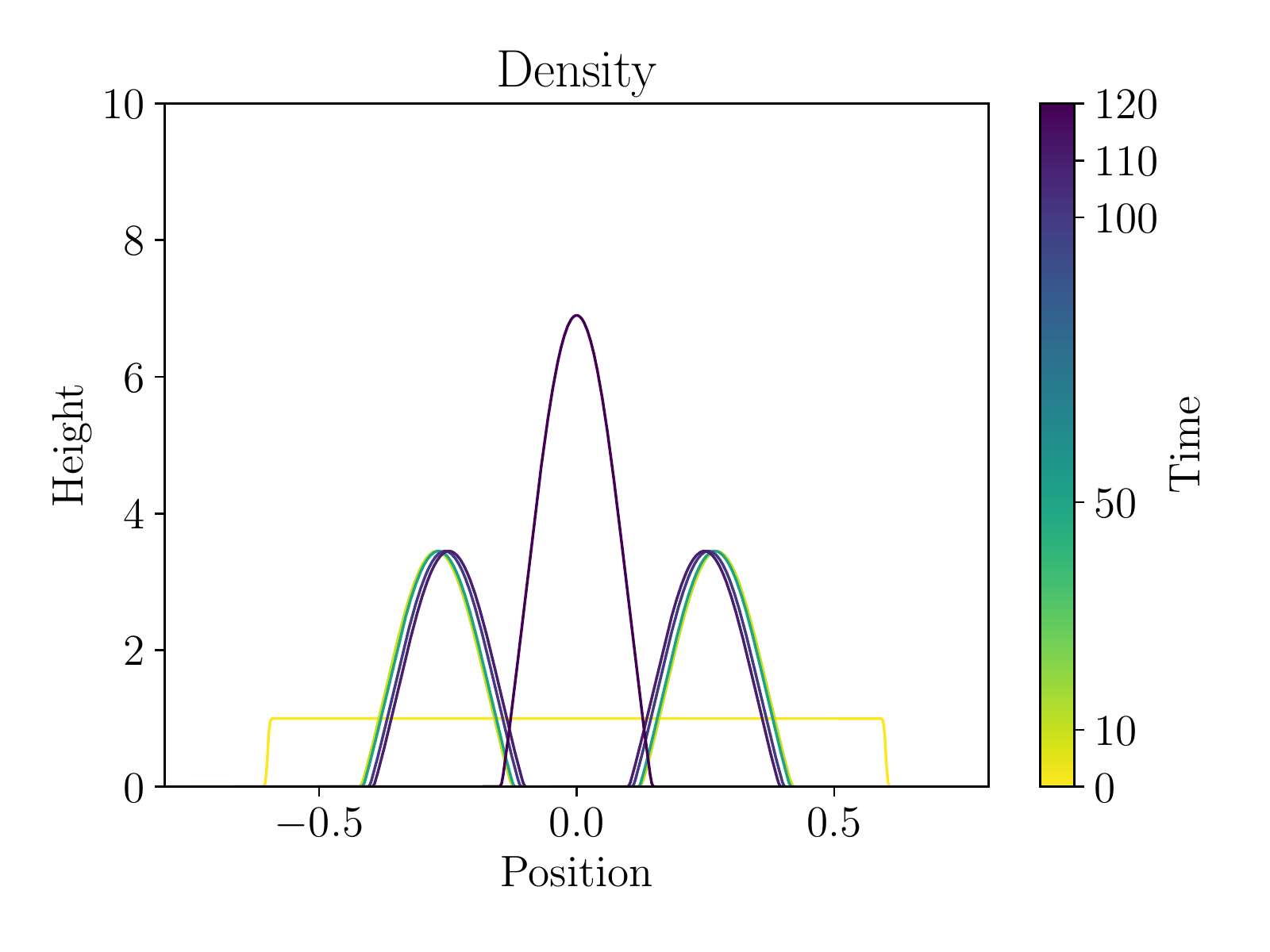} 

\vspace{-.2cm}
 {\bf $\boldsymbol\nu$ = 0.35} \\  
  
 \includegraphics[height=4cm,trim={.6cm .7cm .6cm .6cm},clip, valign=t]{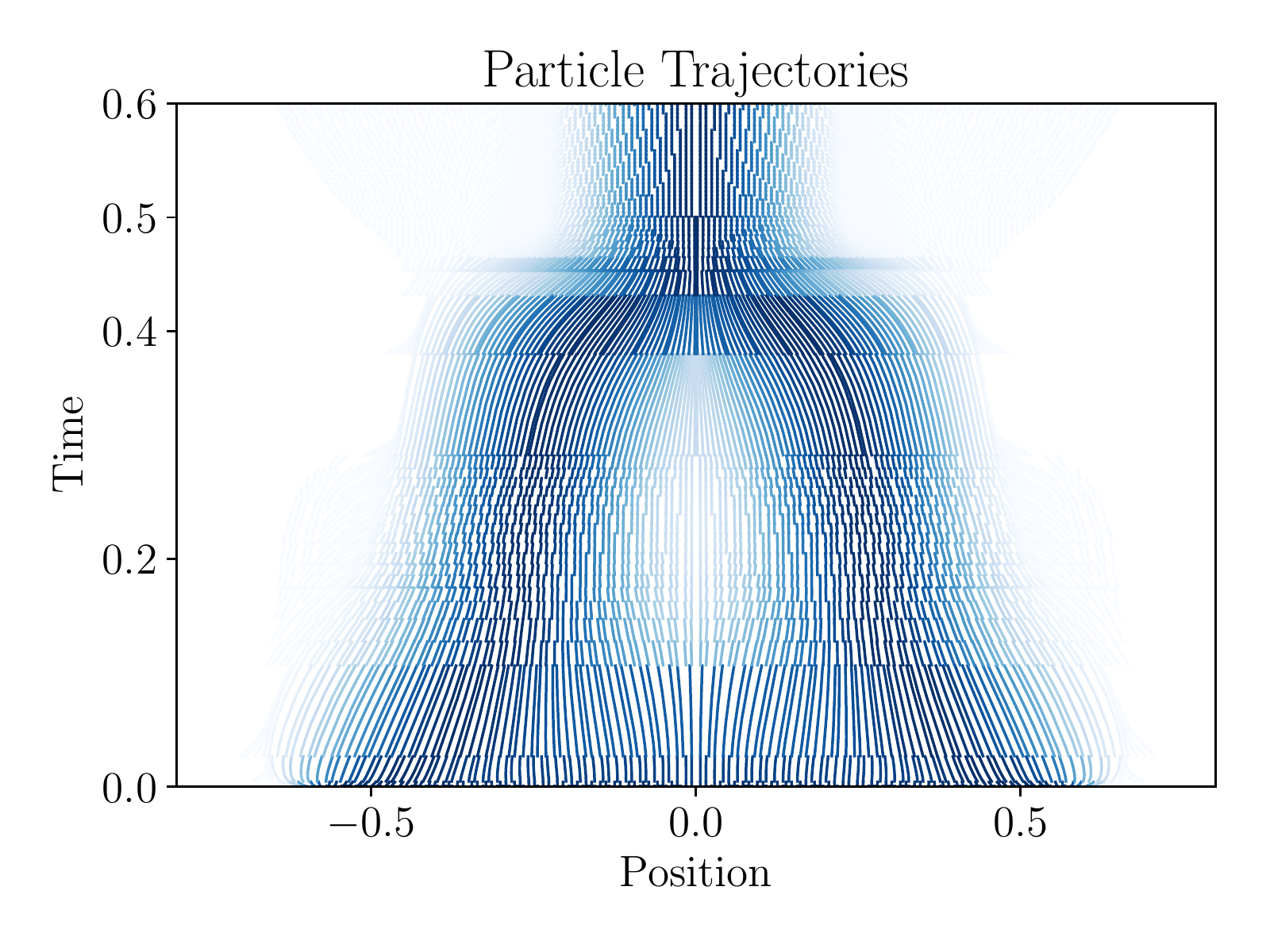} 
  \includegraphics[height=4cm,trim={.6cm .7cm .6cm .6cm},clip, valign=t]{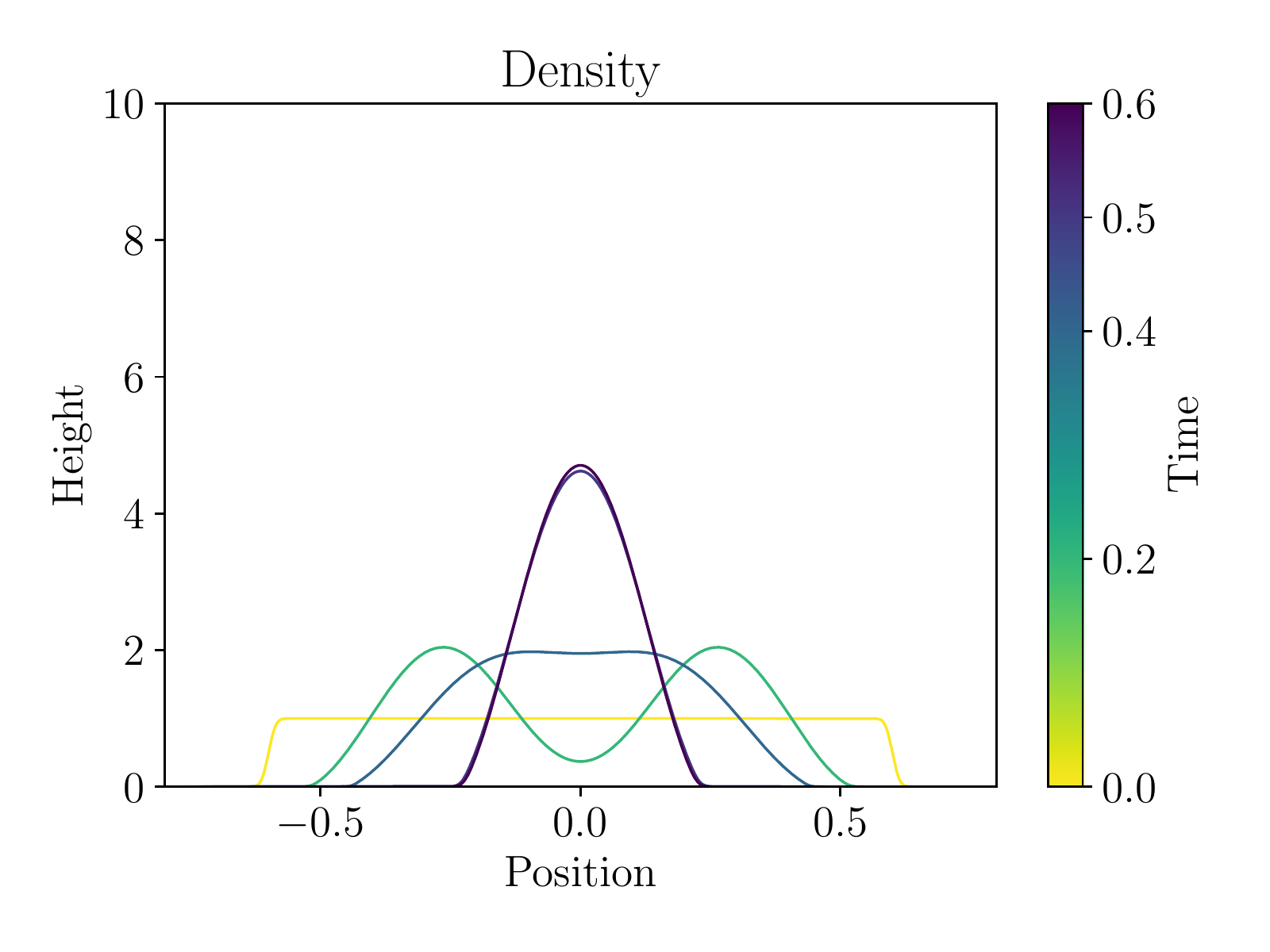}  
  
  \vspace{-.2cm}
   {\bf $\boldsymbol\nu$ = 0.45} \\  
  
 \includegraphics[height=4cm,trim={.6cm .7cm .6cm .6cm},clip, valign=t]{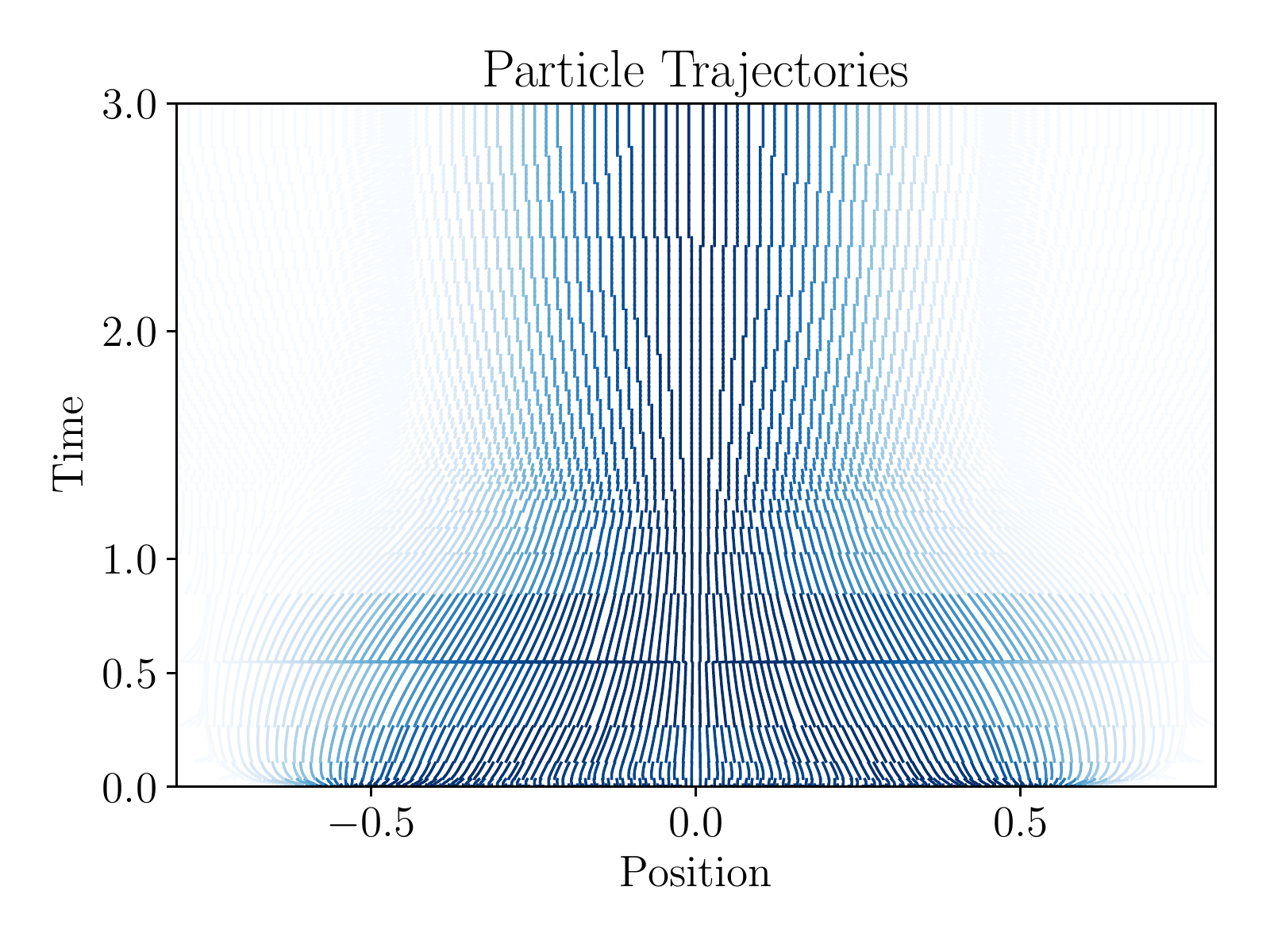} 
  \includegraphics[height=4cm,trim={.6cm .7cm .6cm .6cm},clip, valign=t]{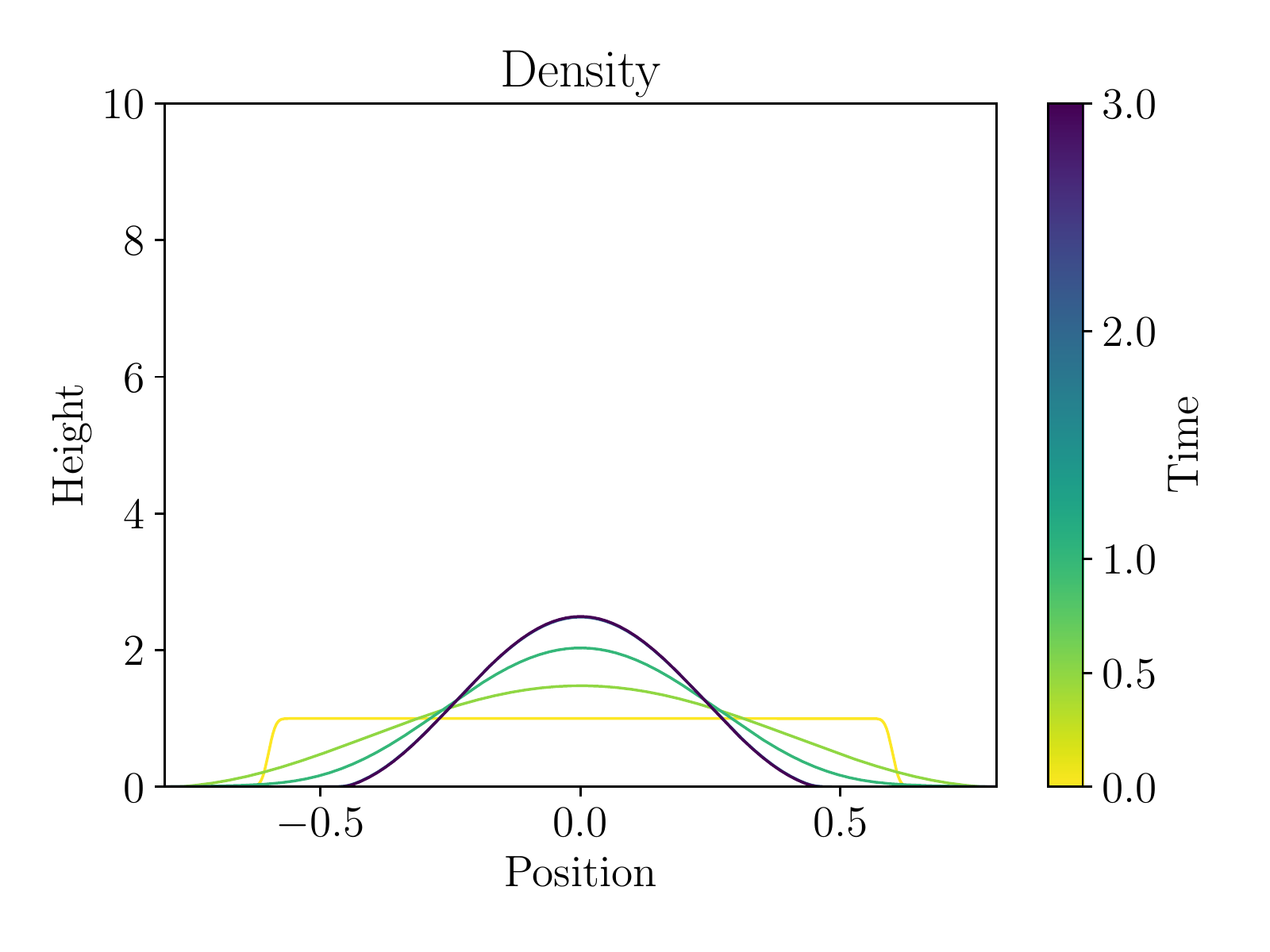}  

\vspace{-.3cm}
\caption{For aggregation diffusion equations with porous medium diffusion ($m=2$),   varying the diffusion coefficient $\nu$ leads to a variety of metastable   profiles, ranging from one to three bumps, with duration varying  more than two orders of magnitude, from $t = 0.6$ to $t = 120$.}
\label{varyingdiffcoeff}
\vspace{-2cm}
\end{center}

\end{figure}

We again consider metastable behavior in Figure \ref{varyingdiffcoeff}, this time considering an attractive Gaussian interaction kernel with fixed localization   $\delta = 0.05$ and varying the diffusion coefficient $\nu$ for the porous medium type diffusion ($m=2$). As in the previous examples,  the structure and duration of the metastable steady states depend strongly on the competition between the aggregation and diffusion terms. In each case, we take the initial data to be $1_{[-0.6,0.6]}$, and we choose $N = 500$ particles, with regularization $\epsilon = 0.9$, and maximal time step $k = 10^{-4}$. 
 For $\nu = 0.15$ in Figure \ref{varyingdiffcoeff},  we observe a three bump metastable state, a two bump metastable state, and ultimately a one bump steady state, which is reached by time $t = 2.5$. For $\nu = 0.25$, we recover the fourth example from Figure \ref{varyinginteractionrange}, which is unique in that it is not until approximately $t \sim 120$ until the metastable steady state subsides and the solution reaches a radially decreasing profile.
For $\nu = 0.35$, the solution has a two bump metastable state and a one bump steady state, which is reached by time $t =0.6$, and for $\nu = 0.45$, the solution remains radially decreasing for all time and approaches a one bump steady state by time $t = 2.5$. For $\nu > 0.5$, DiFrancesco and Jaafra showed that steady states do not exist, as all solutions decay to zero locally in $L^2(\Rd)$ \cite{DJ}.

\begin{figure}[H]
\begin{center}
\vspace{-.25cm}
 {\bf Metastability for Linear Difusion,  $W(x) = e^{-x^2/4 (.1)^2}/\sqrt{4 \pi (.1)^2}$}  \\ 
\vspace{.1cm}

 {\bf $T_{max} = 0.2$} \\

 \includegraphics[height=4cm,trim={.6cm .7cm .6cm .6cm},clip, valign=t]{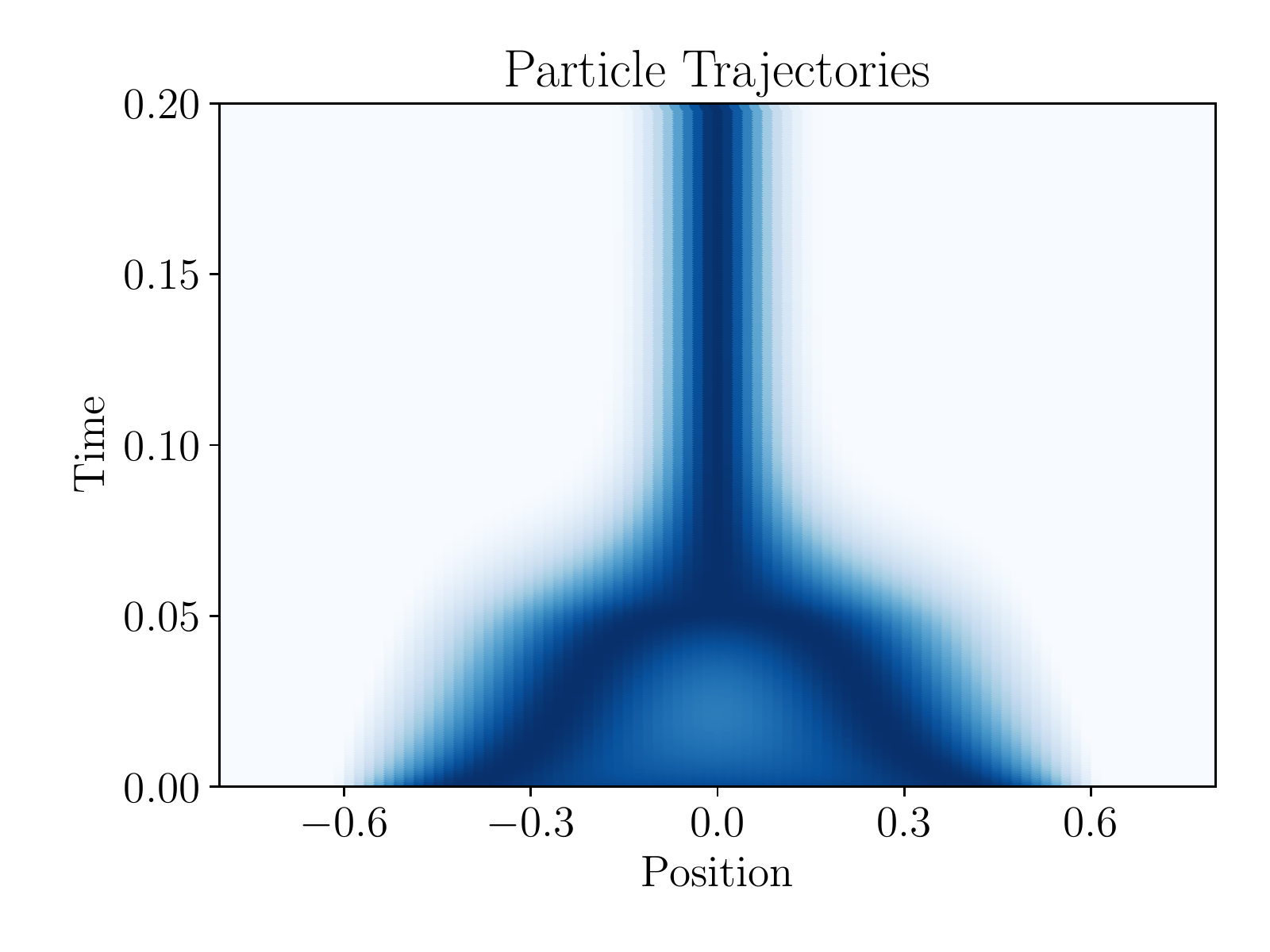} 
  \includegraphics[height=4cm,trim={.6cm .7cm .6cm .6cm},clip, valign=t]{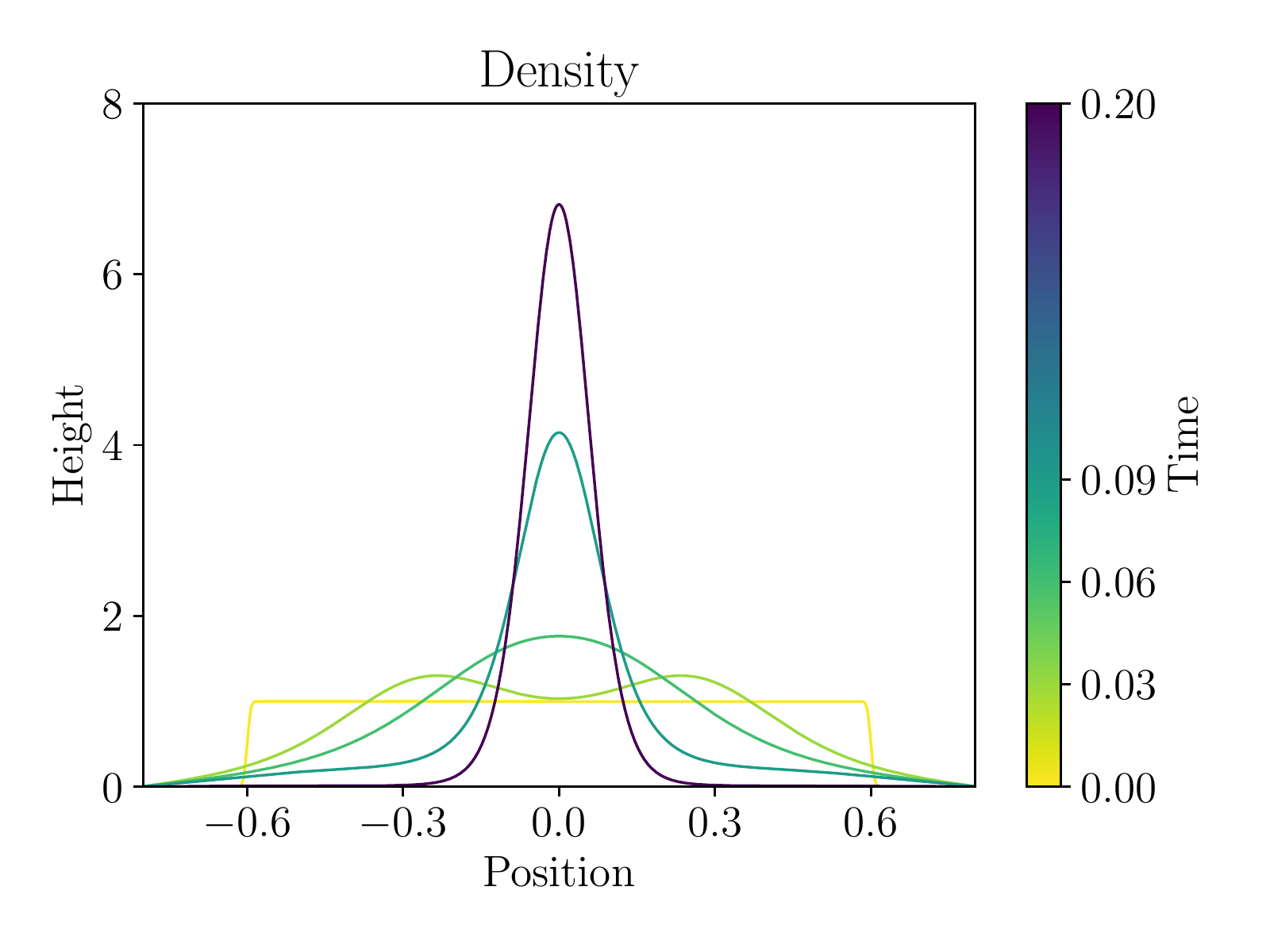}  

\vspace{.1cm}

 {\bf $T_{max} = 3.0$} \\
 
  \includegraphics[height=4cm,trim={.6cm .7cm .6cm .6cm},clip, valign=t]{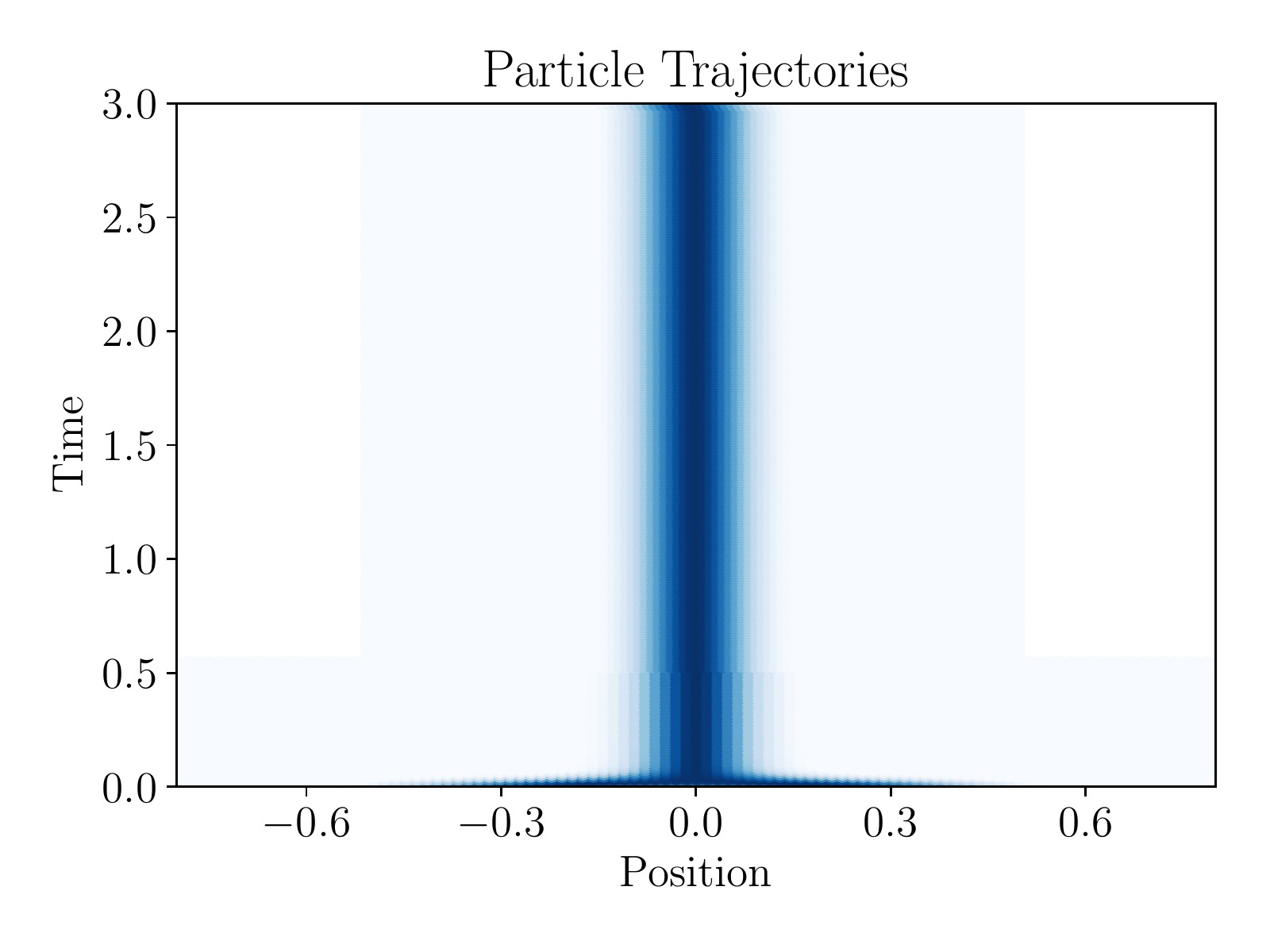} 
  \includegraphics[height=4cm,trim={.6cm .7cm .6cm .6cm},clip, valign=t]{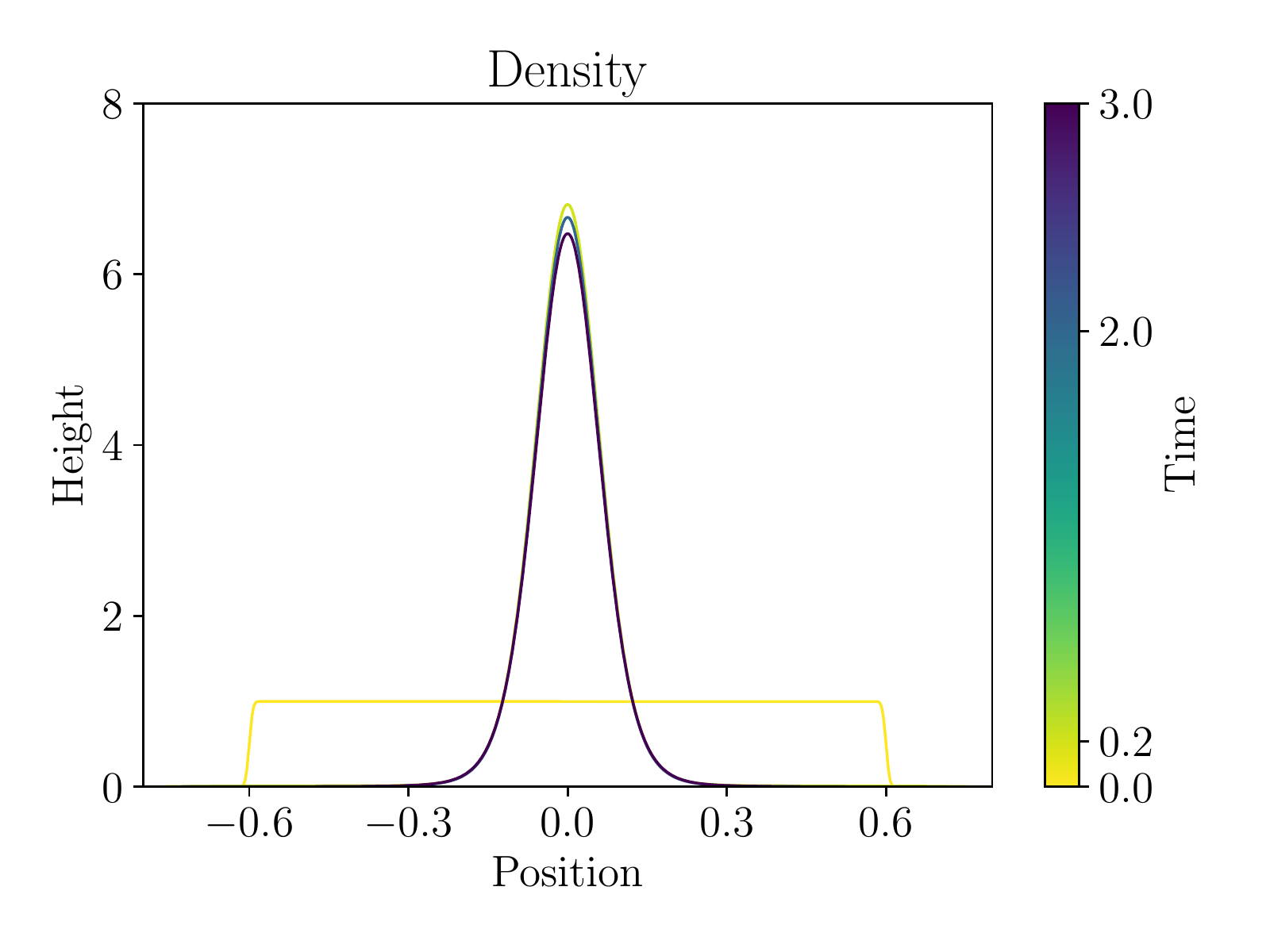}  
\vspace{-.3cm}
\caption{For aggregation diffusion equations with a small amount of linear diffusion ($m=1, \nu = 0.2$),  we observe rapid convergence to a single bump metastable state (top), followed by very slow spreading as the bump diffuses away (bottom).}
 \label{lineardiffusionmetastablefig}
\vspace{-.8cm}
\end{center}

\end{figure}

In Figure \ref{lineardiffusionmetastablefig}, we consider how linear diffusion affects metastability behavior in aggregation diffusion equations. As in the previous Figure \ref{varyingdiffcoeff}, we consider an  attractive Gaussian interaction kernel with fixed localization   $\delta = 0.1$ and   diffusion coefficient $\nu =0.2$. Similarly to the porous medium case, we observe rapid convergence to a metastable state, in this case  forming a single bump. However, unlike in the porous medium case, a steady state does not exist \cite{CDP}. Instead, we observe that the metastable state  spreads slowly as the bump diffuses away. In our simulation, we take the initial data to be $1_{[-0.6,0.6]}$, and we choose $N = 500$ particles, with regularization $\epsilon = 0.9$, and maximal time step $k = 10^{-3}$.

\begin{figure}[H]
\vspace{-.5cm}
\begin{center}
 {\bf Metastability Equilibrium for Varying Interaction Range, $W(x) = \log(|x/\delta|)/\delta$}  \\ 

 {\bf $\boldsymbol\delta$ = 0.1} \\

 \includegraphics[height=4cm,trim={.6cm .7cm .6cm .6cm},clip, valign=t]{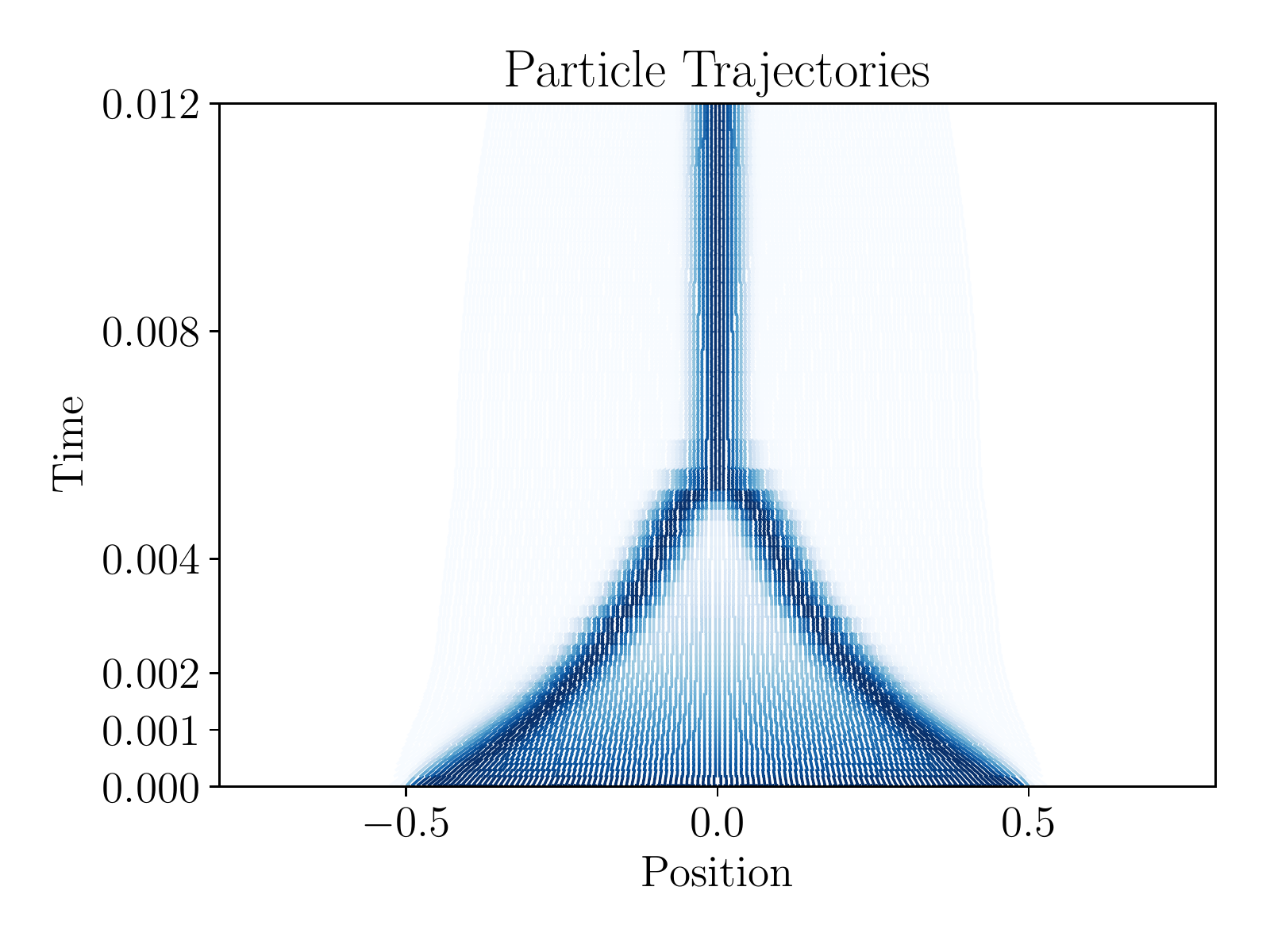} 
  \includegraphics[height=4cm,trim={.3cm .7cm .6cm .6cm},clip, valign=t]{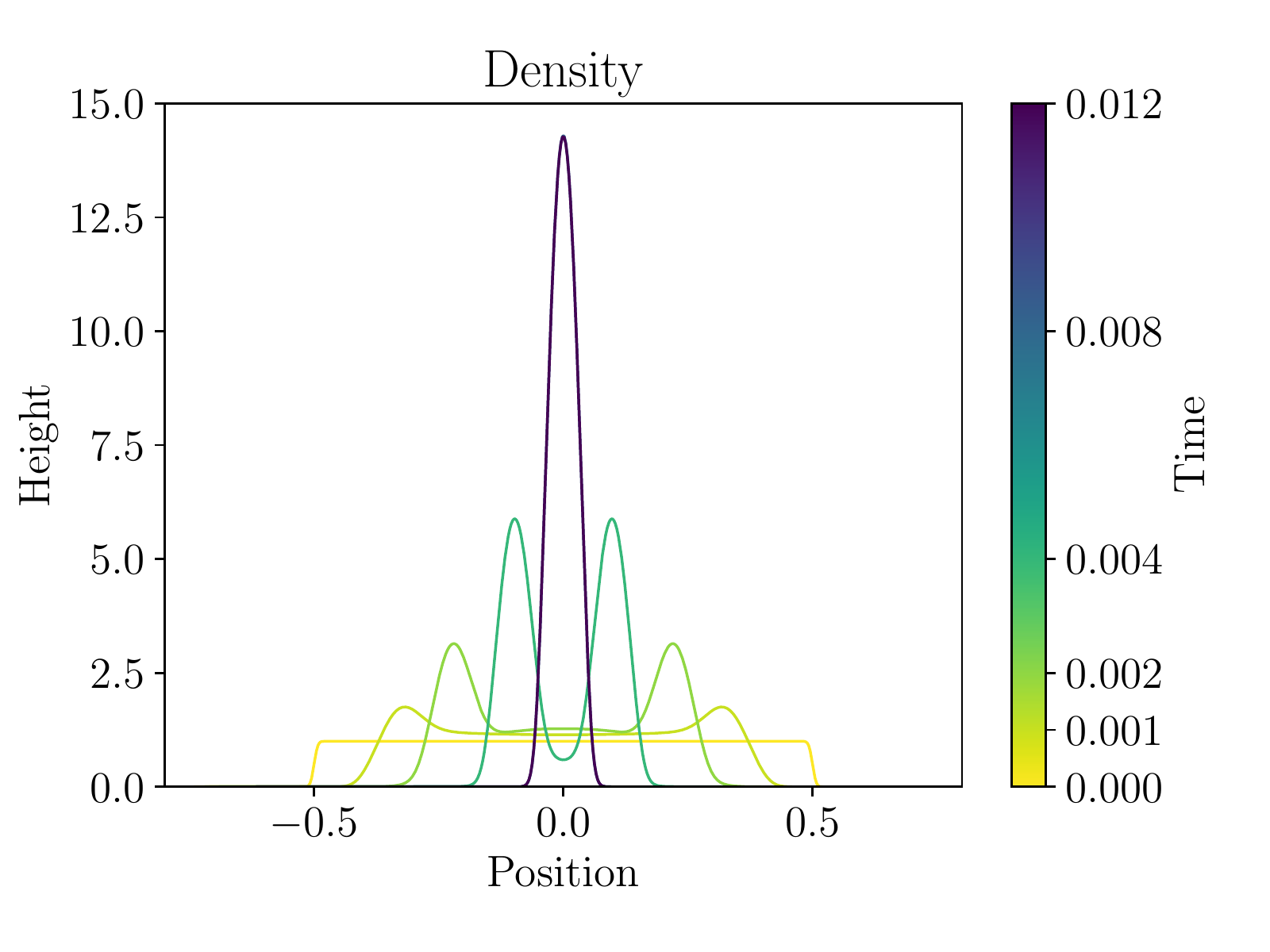}

\vspace{-.2cm}
 {\bf $\boldsymbol\delta$ = 0.5} \\

 \includegraphics[height=4cm,trim={.6cm .7cm .6cm .6cm},clip, valign=t]{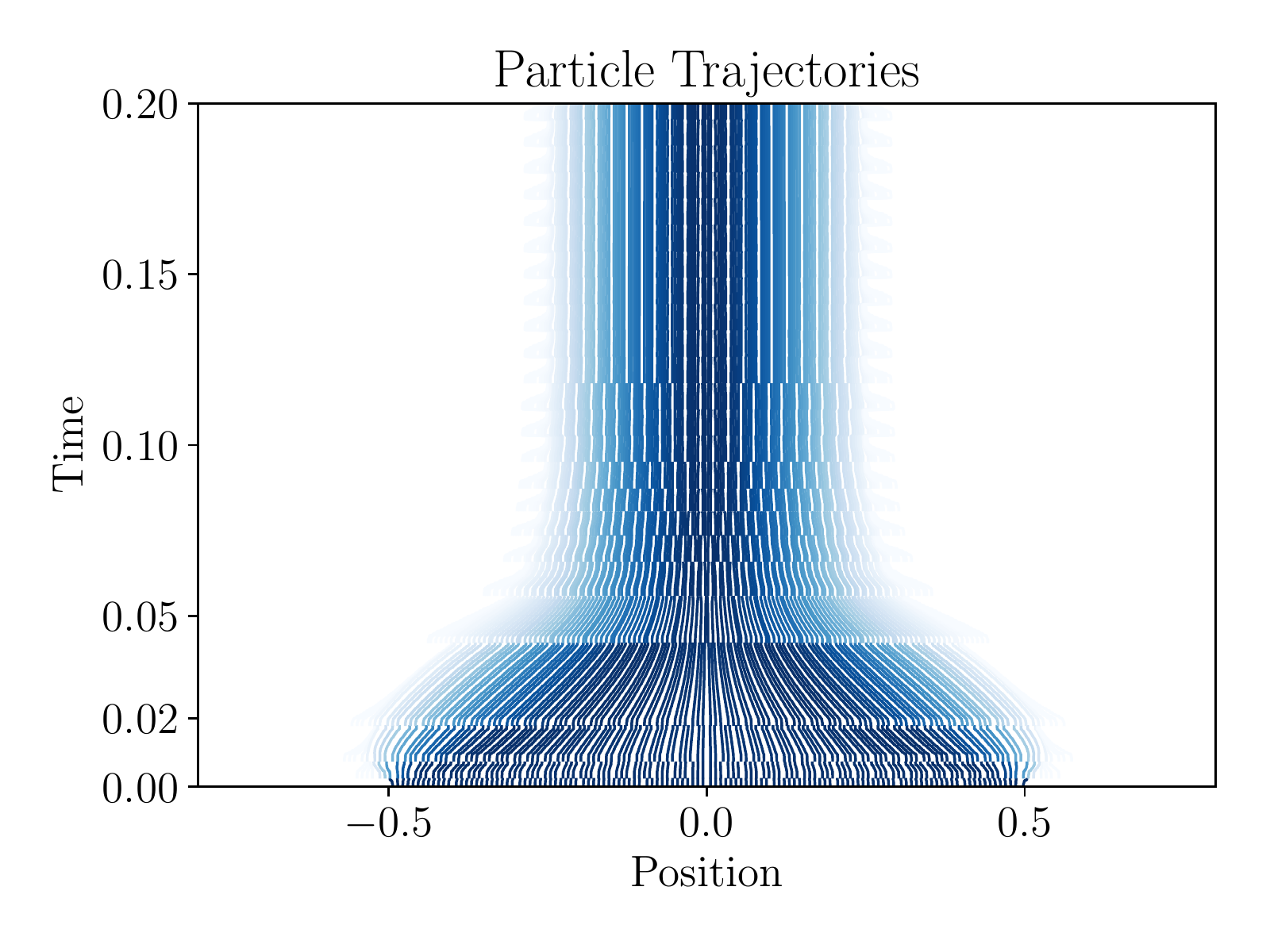} 
  \includegraphics[height=4cm,trim={.3cm .7cm .6cm .6cm},clip, valign=t]{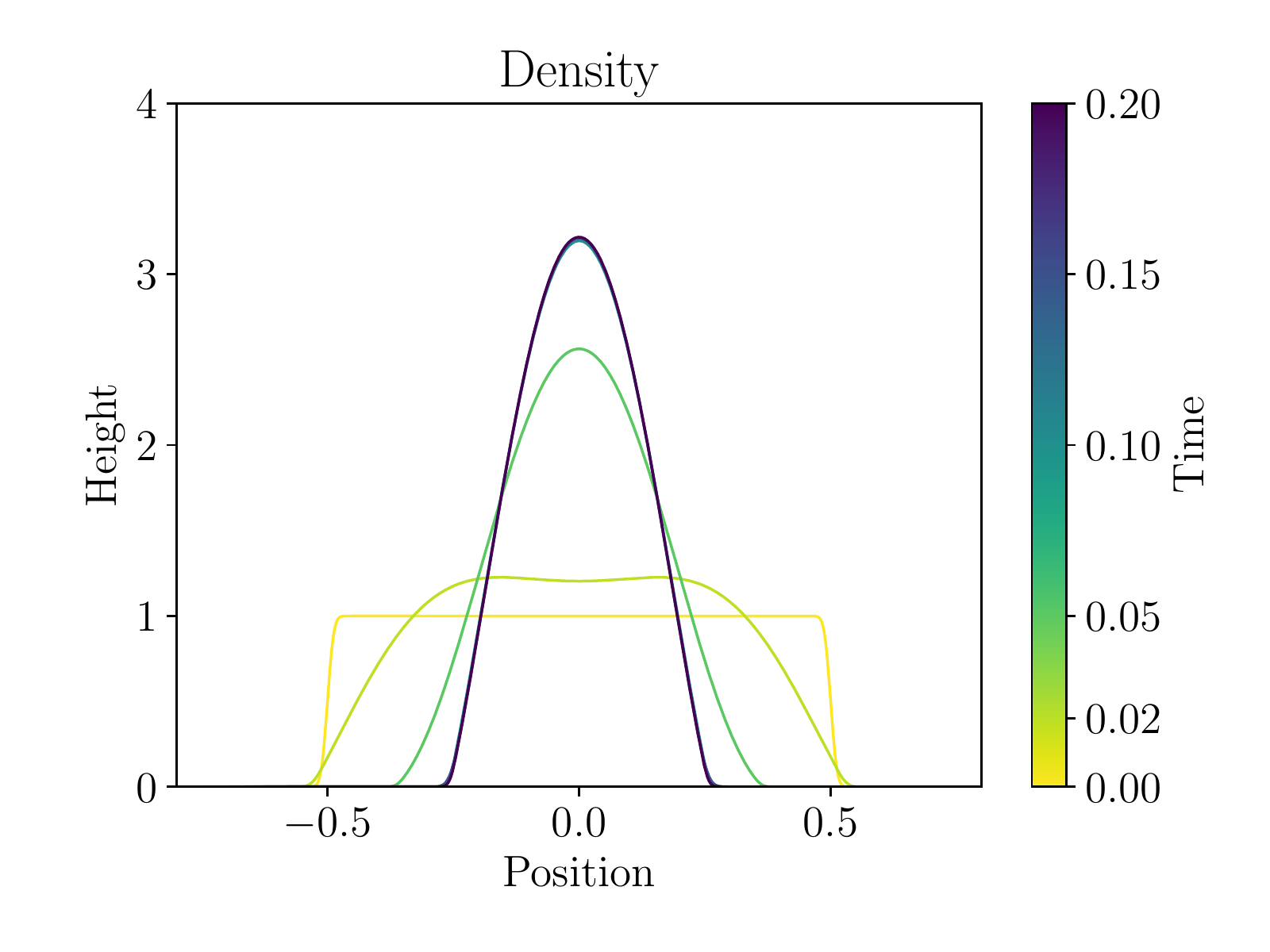}  

\vspace{-.2cm}
 {\bf $\boldsymbol\delta$ = 1.0} \\

 \includegraphics[height=4cm,trim={.6cm .7cm .6cm .6cm},clip, valign=t]{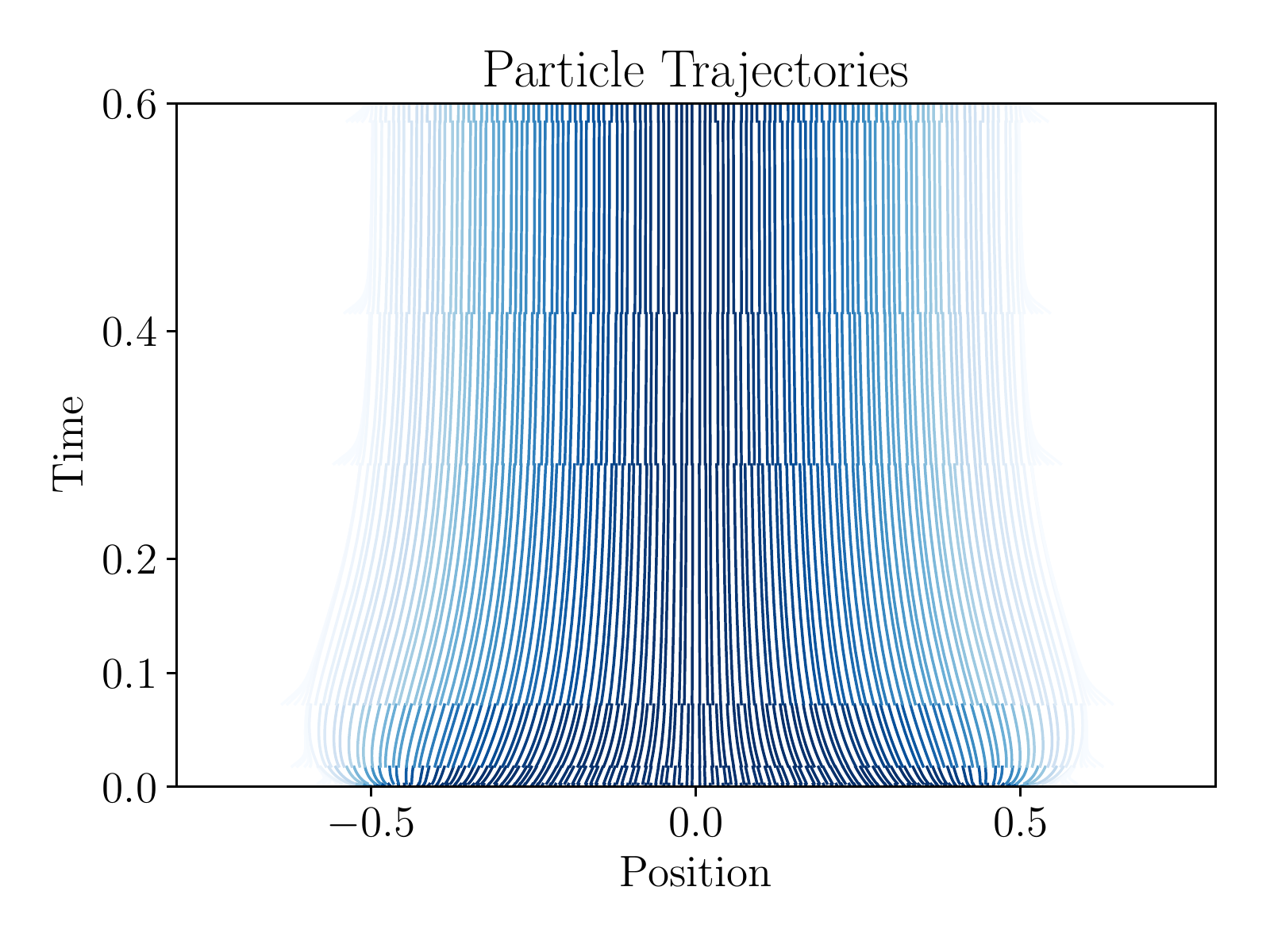} 
  \includegraphics[height=4cm,trim={.3cm .7cm .6cm .6cm},clip, valign=t]{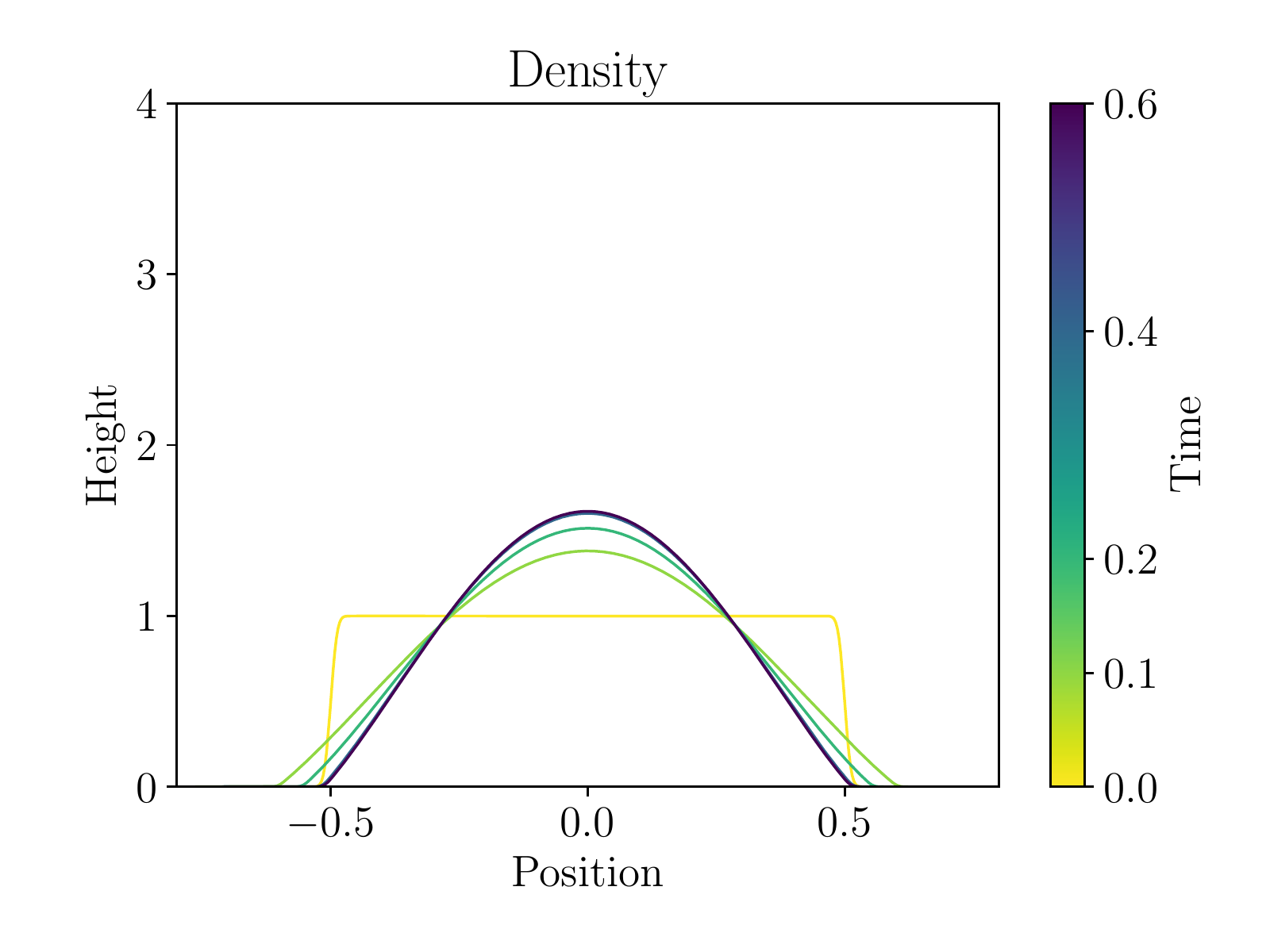}  

\caption{For aggregation diffusion equations with a singular, attractive interaction kernel and porous medium diffusion ($m=2$), the strong short range aggregation prevents the solution from forming distinct, metastable clusters, as observed in  Figure \ref{varyinginteractionrange}). }
\vspace{-.7cm}
\label{metalogfig}
\end{center}

\end{figure}
 
In Figure \ref{metalogfig}, we consider the behavior of localized aggregation with a singular interaction potential. Unlike in Figure \ref{varyinginteractionrange}, in which we considered smooth Gaussian interactions, the strong short range aggregation prevents the solution from breaking apart into distinct clusters and causes the solution to quickly approaches a radially decreasing equilibrium. In each simulation, we consider porous medium diffusion $m=2$ with diffusion coefficient $\nu = 0.4$, and we take $N = 500$, $\epsilon = 0.9$, and the maximum time step $k = 10^{-5}$. We consider the initial data $1_{[-0.5,0.5]}$ and restrict to the computational domain $[-0.55,0.55]$. 

\begin{figure}[H]
\begin{center}
 {\bf Evolution to Equilibrium for Varying Interaction Range, $W(x) = |x|/\delta^2$}  \\ 
\vspace{.1cm}

 {\bf $\boldsymbol\delta$ = 0.05} \\
  \includegraphics[height=4cm,trim={.6cm .7cm .6cm .6cm},clip, valign=t]{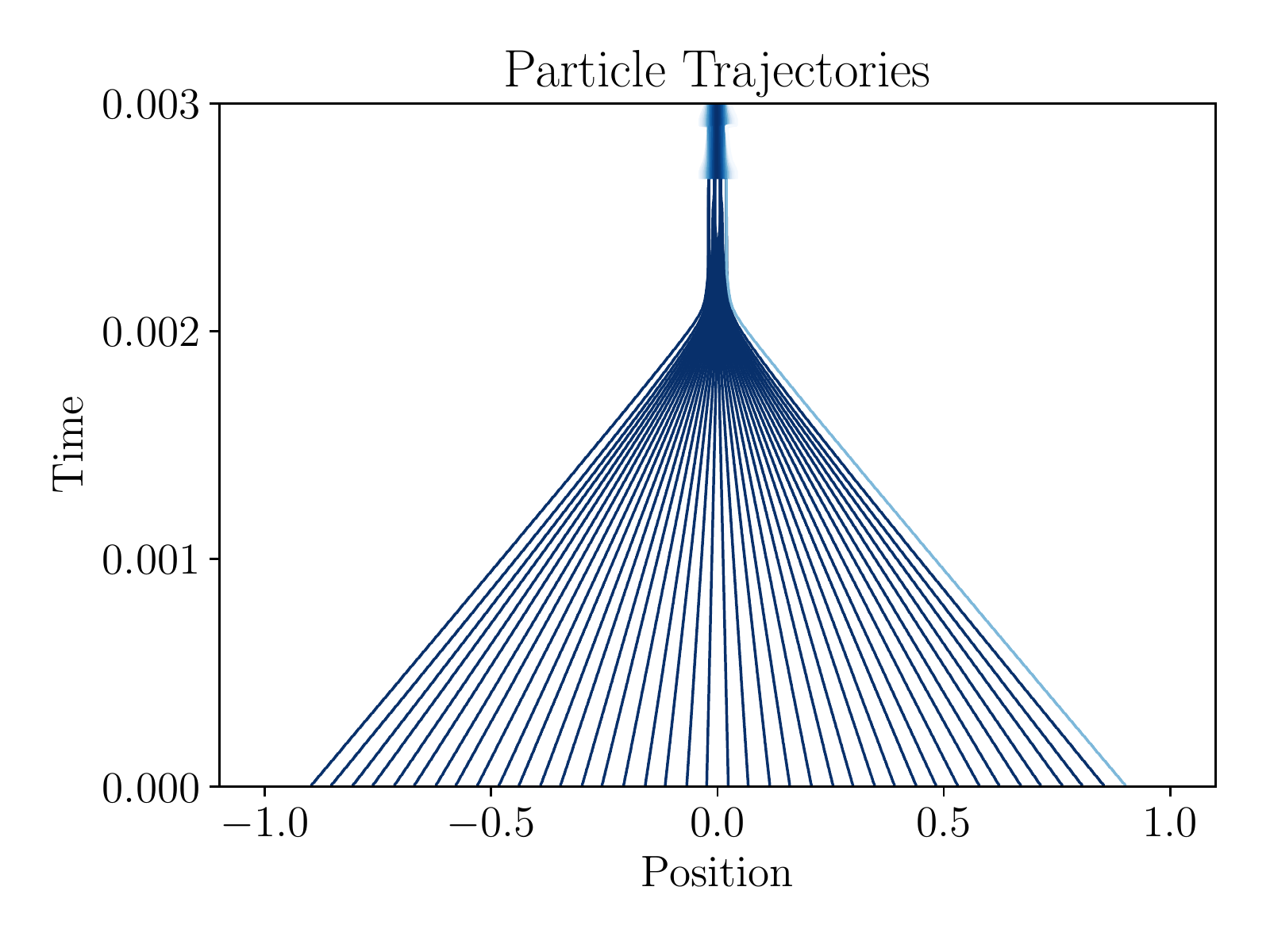}  
 \includegraphics[height=4cm,trim={.6cm .7cm .6cm .6cm},clip, valign=t]{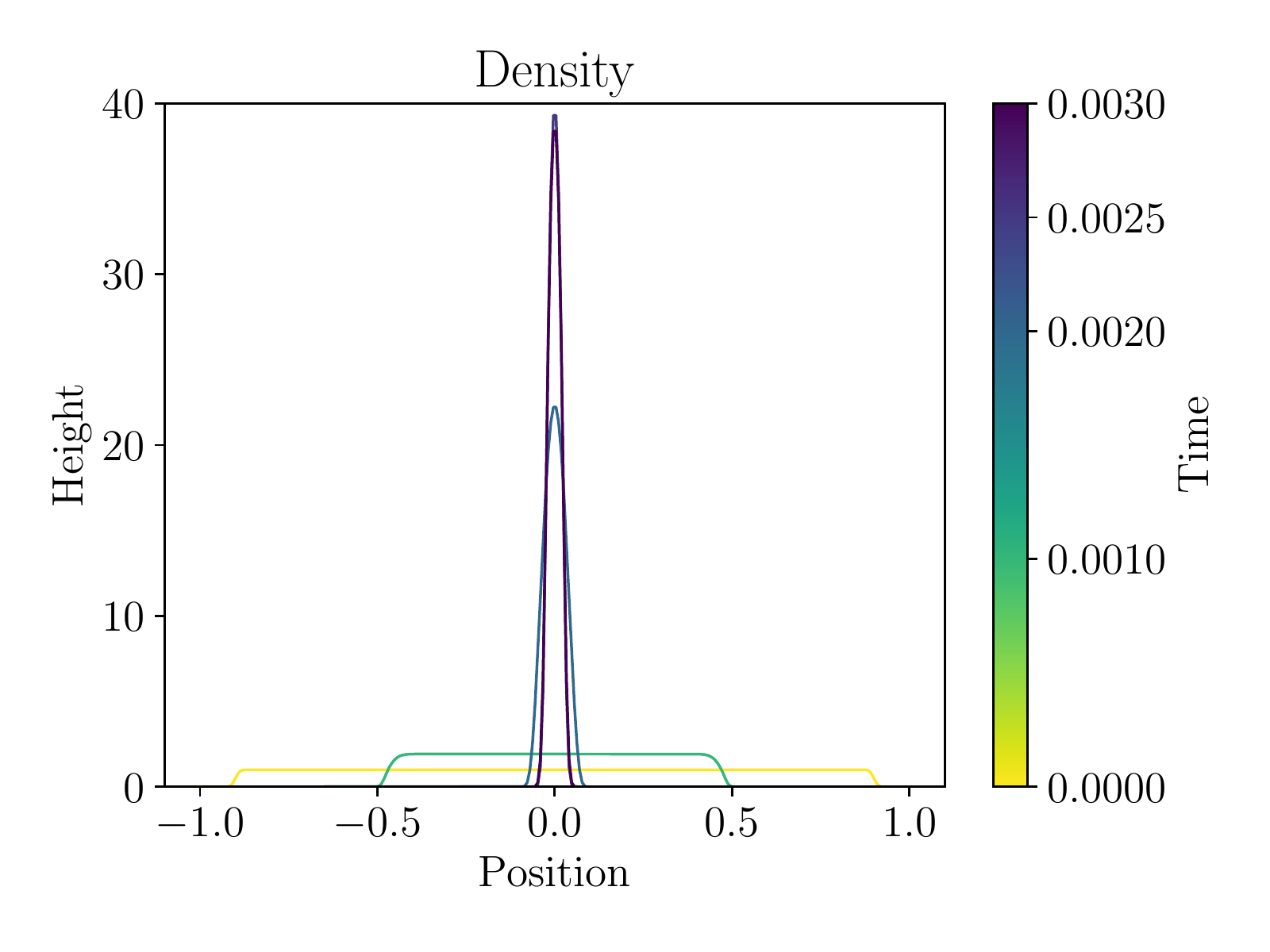} 
 
 \vspace{-.2cm}
  {\bf $\boldsymbol\delta$ = 0.1} \\
  \includegraphics[height=4cm,trim={.6cm .7cm .6cm .6cm},clip, valign=t]{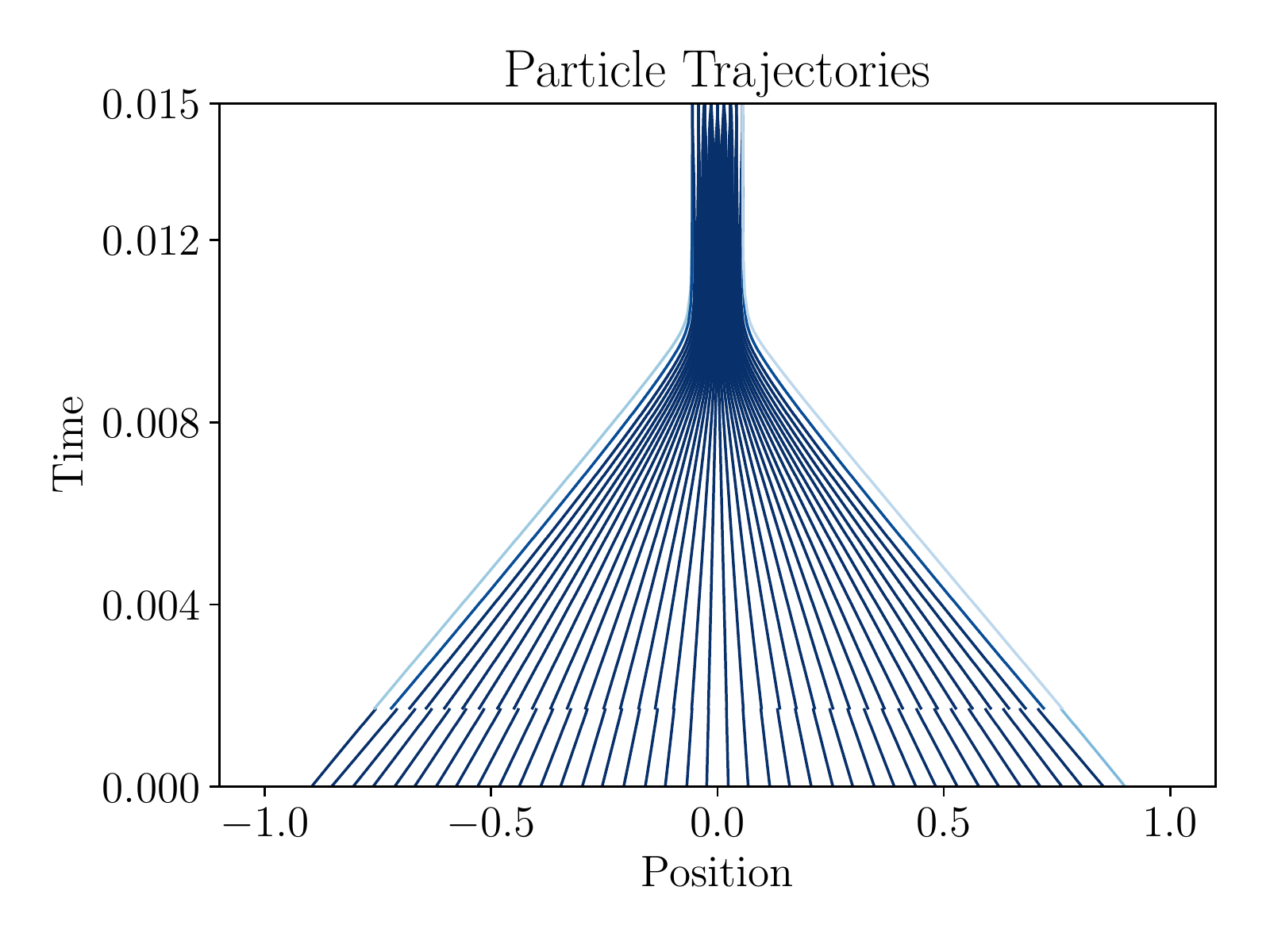}  
 \includegraphics[height=4cm,trim={.6cm .7cm .6cm .6cm},clip, valign=t]{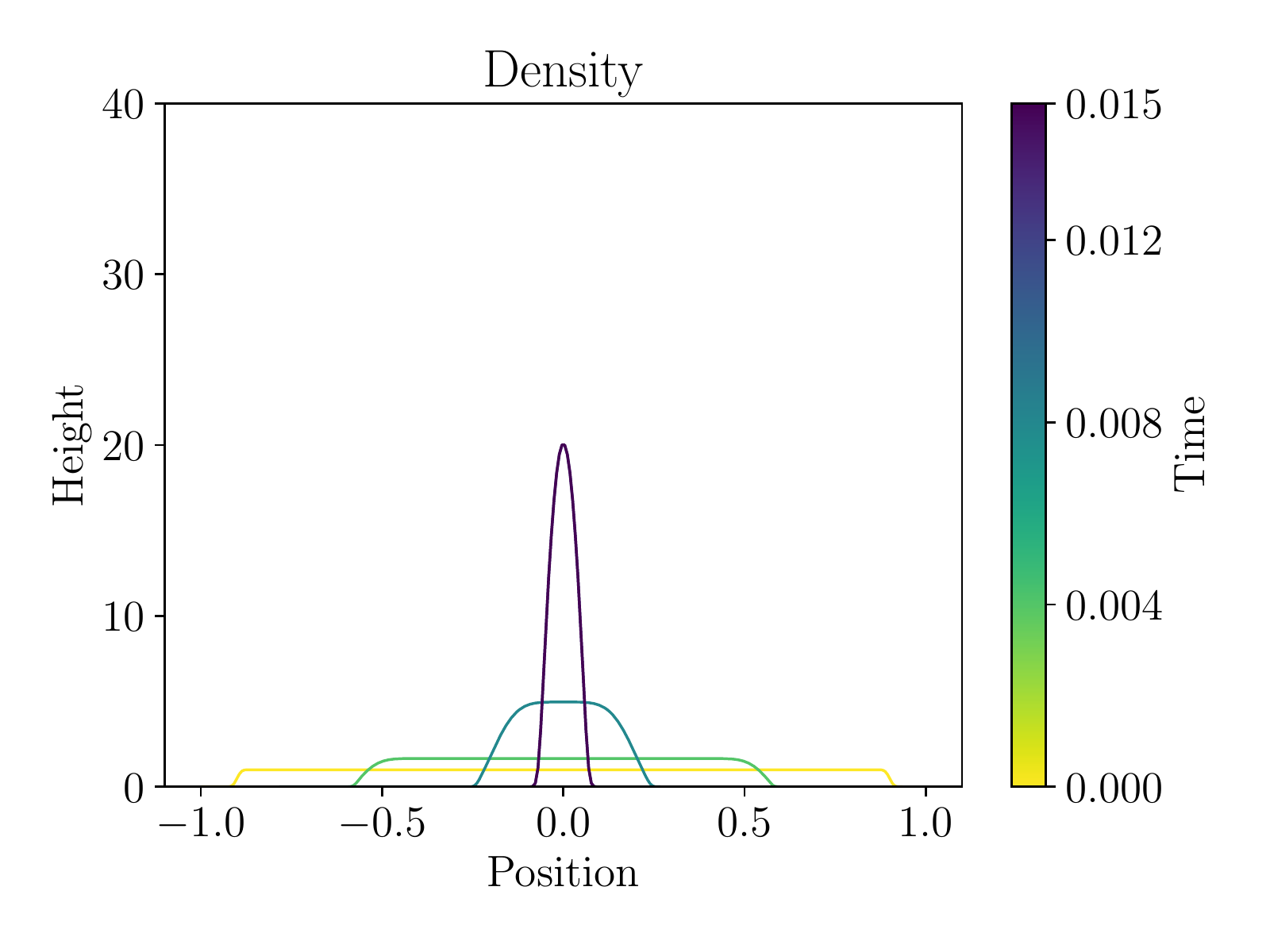} 

 \vspace{-.2cm}  
  {\bf $\boldsymbol\delta$ = 0.5} \\
  \includegraphics[height=4cm,trim={.6cm .7cm .6cm .6cm},clip, valign=t]{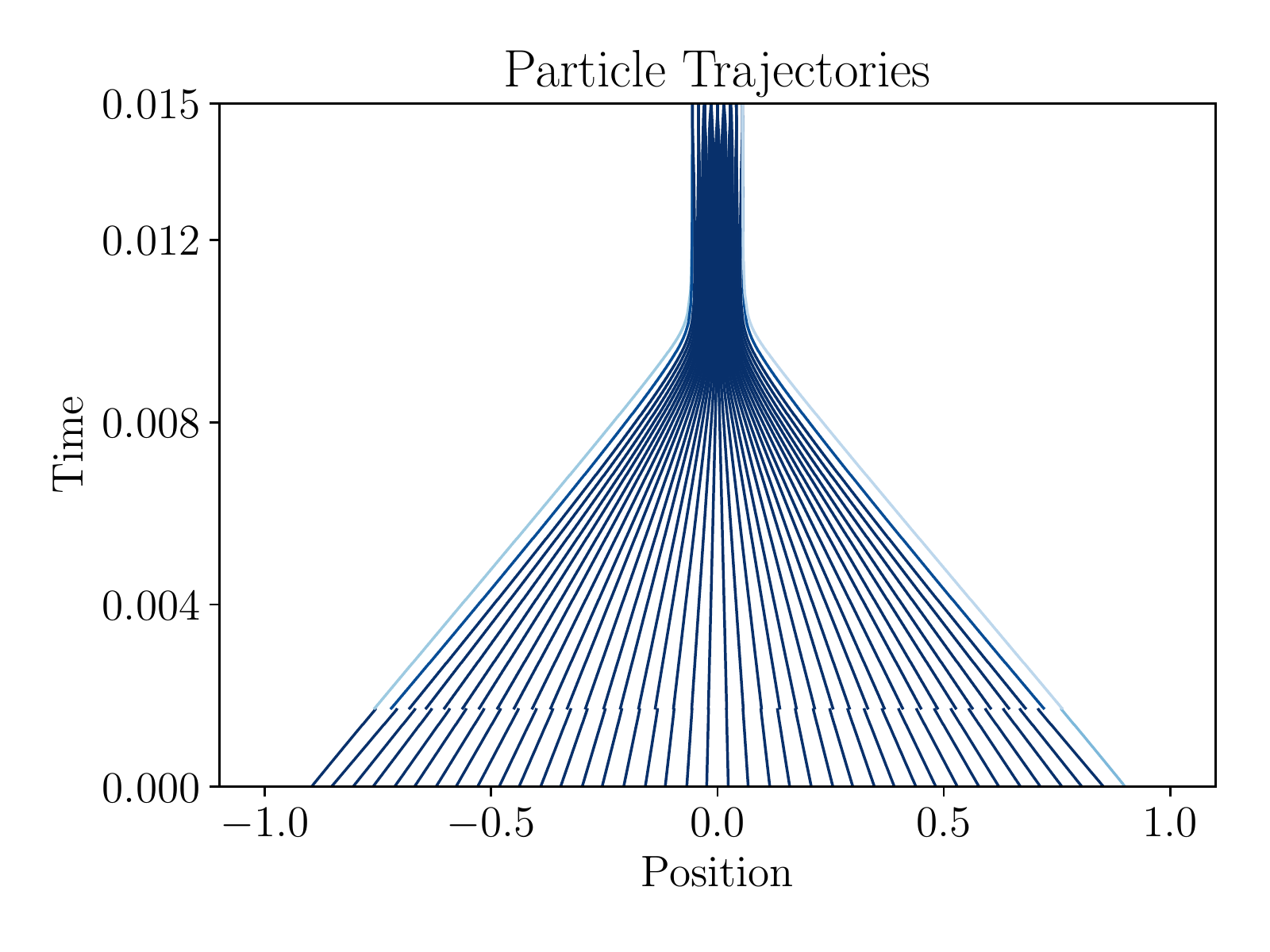}  
 \includegraphics[height=4cm,trim={.6cm .7cm .6cm .6cm},clip, valign=t]{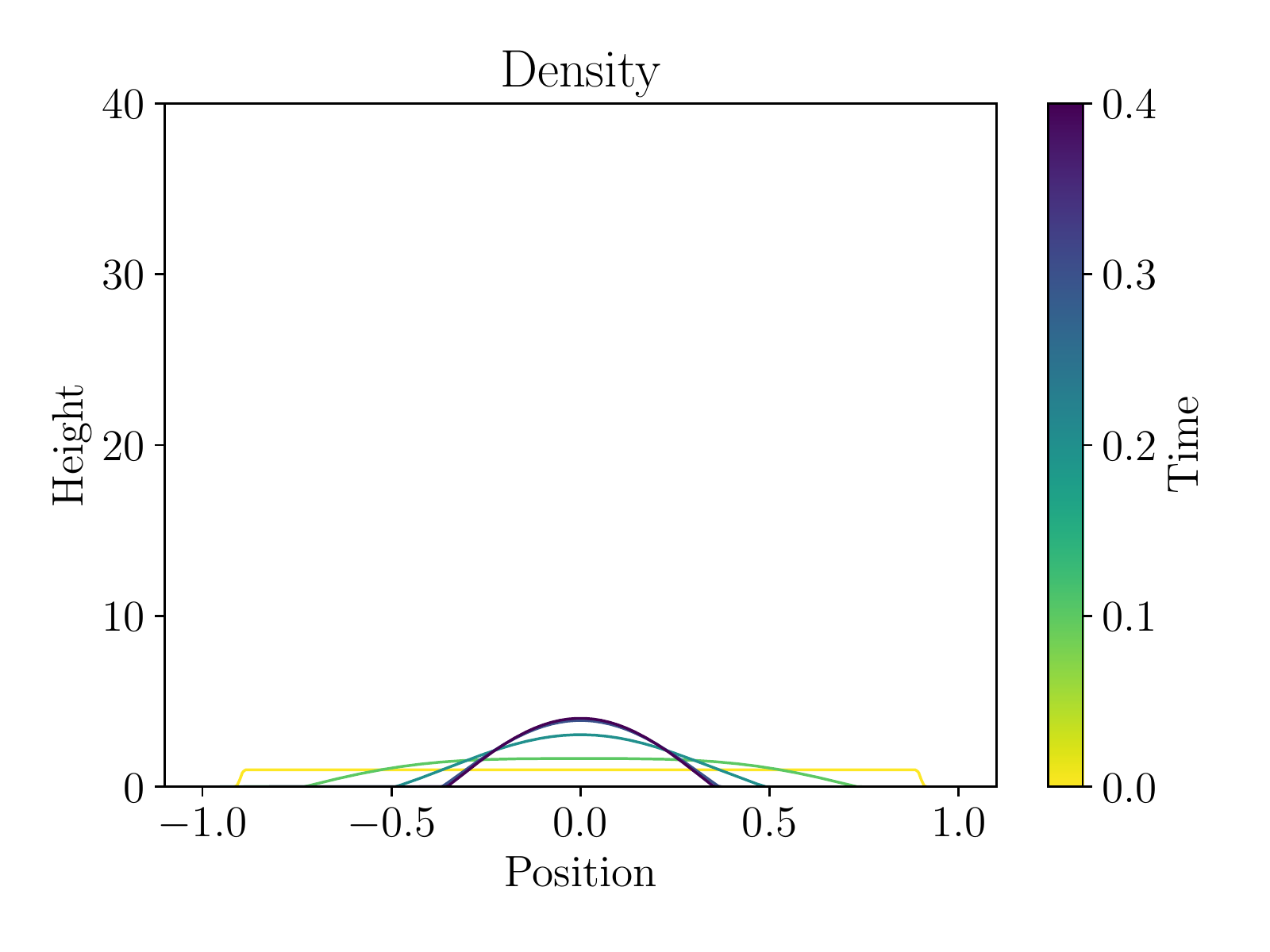} 
 
\caption{Unlike in the cases of aggregation diffusion equations with either a Gaussian or logarithmic localized interaction kernel (see Figures \ref{varyinginteractionrange} and \ref{metalogfig}), for a Newtonian interaction kernel, the qualitative dynamics remain the same as the scale of the interaction varies: solutions with radially decreasing initial data remain radially decreasing for all time.}
 \label{metanewtfig}
 \end{center}
  \vspace{-.7cm}
 \end{figure}

In Figure \ref{metanewtfig}, we investigate the competition between aggregation and porous medium diffusion $m=2$ as the interaction localizes for a third choice of interaction kernel: the attractive Newtonian potential. In terms of its singularity at the origin, this kernel lies between the very singular logarithmic kernel from figure \ref{metalogfig} and the smooth Gaussian kernel from Figure \ref{varyinginteractionrange}. However, unlike in both of these cases, we observe that the solution remains radially decreasing for all times, quickly approaches a single bump steady state. We do not observe any metastability behavior.   For each simulation, we consider diffusion coefficient $\nu = 0.1$, $N = 600$, $\epsilon = 0.9$, and the maximum time step $k = 10^{-5}$. We consider the initial data $1_{[-0.9,0.9]}$.  

\begin{figure}[H]
\begin{center}
 {\bf Metastability for Varying Diffusion, $W(x) = |x|^4/4 - |x|^2/2$}  \\ 
\vspace{.1cm}

 {\bf m=1, $ \bf \text{\boldmath$\nu$} = 0.075^2/2$} \\
 \includegraphics[height=4cm,trim={.6cm .7cm .6cm .6cm},clip, valign=t]{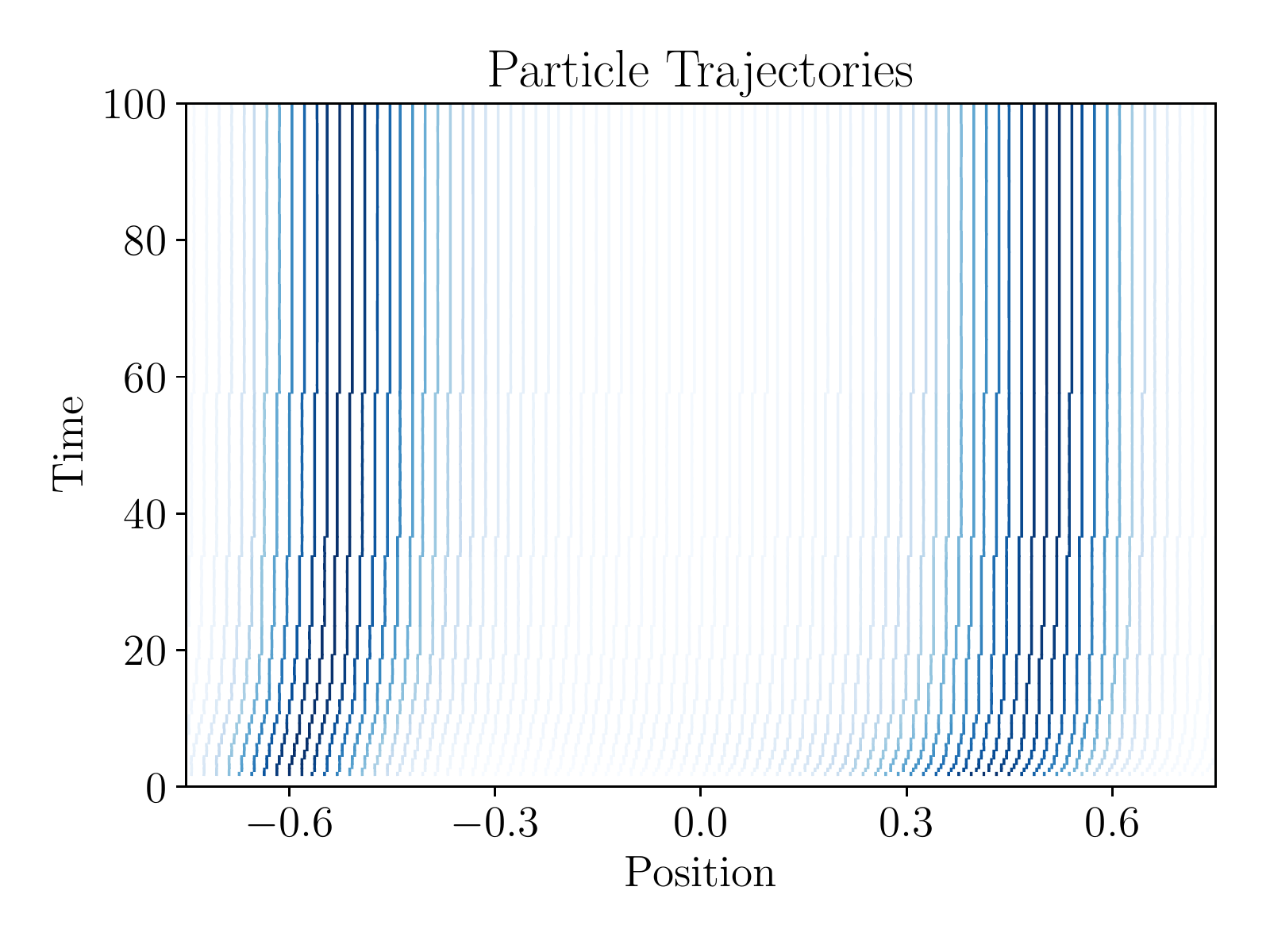} 
 \includegraphics[height=4cm,trim={.6cm .7cm .6cm .6cm},clip, valign=t]{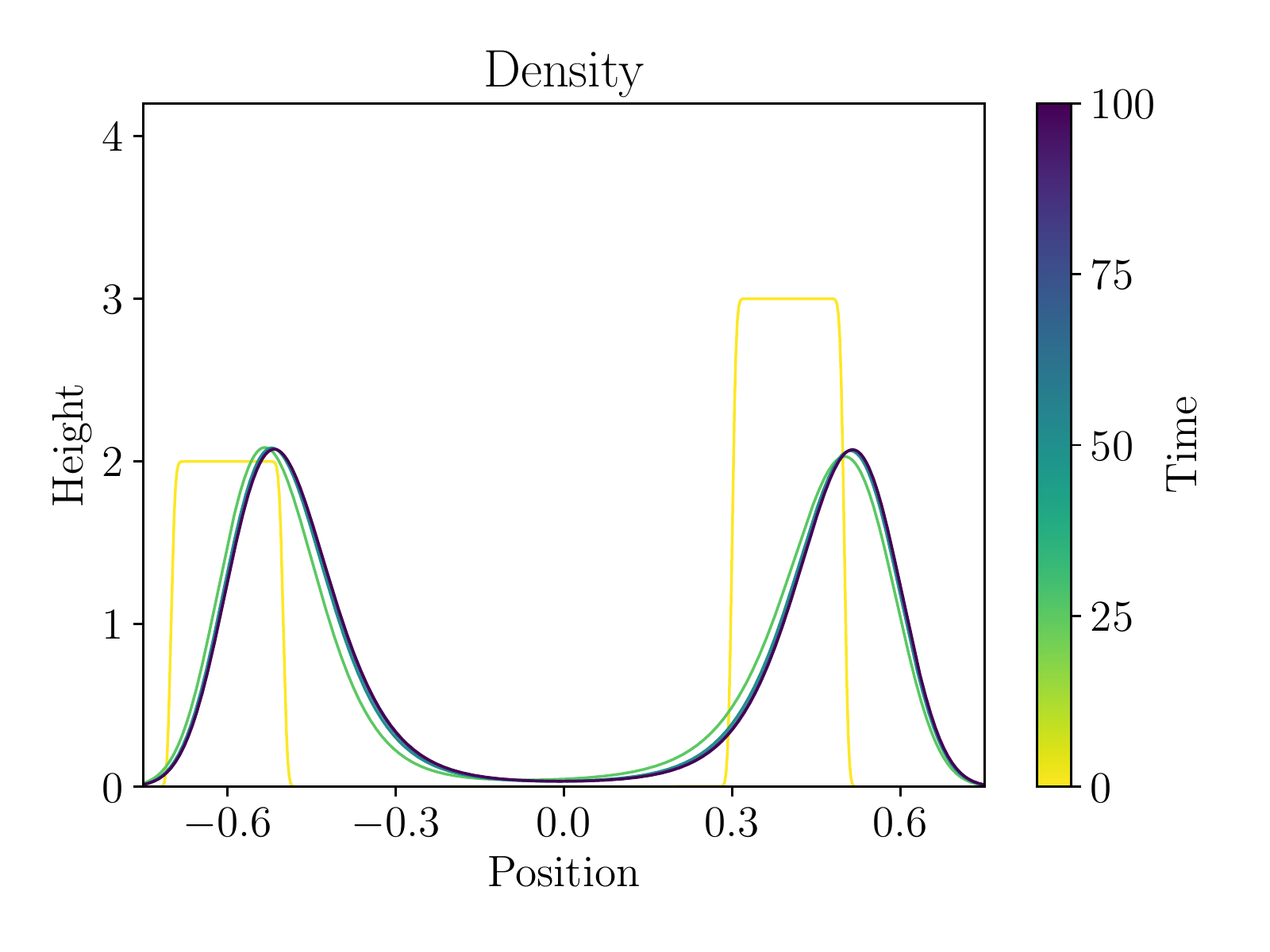} 
 
  {\bf m=2, $ \bf \text{\boldmath$\nu$} = 0.075^2/2$} \\
 \includegraphics[height=4cm,trim={.6cm .7cm .6cm .6cm},clip, valign=t]{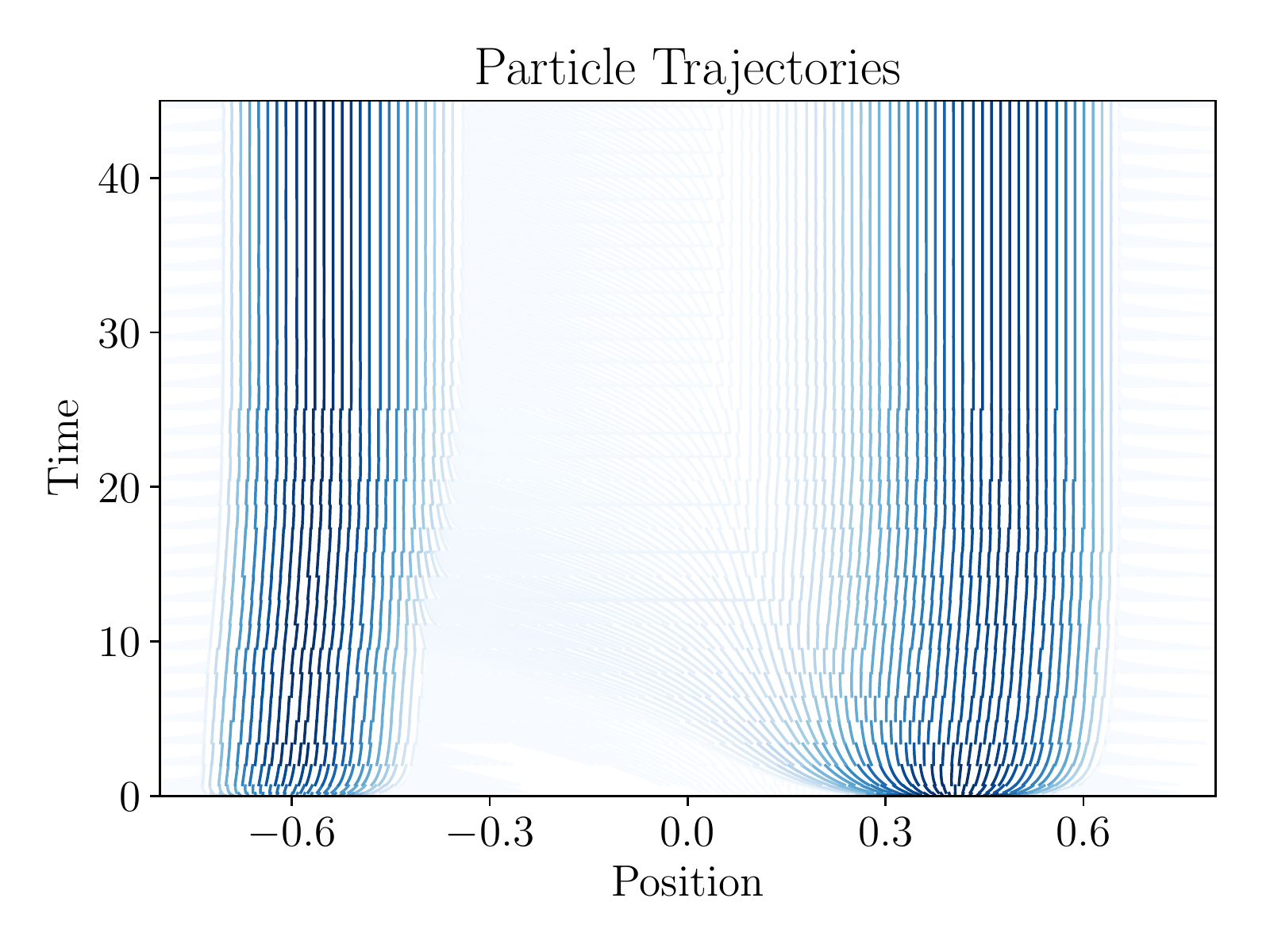}
 \includegraphics[height=4cm,trim={.6cm .7cm .6cm .6cm},clip, valign=t]{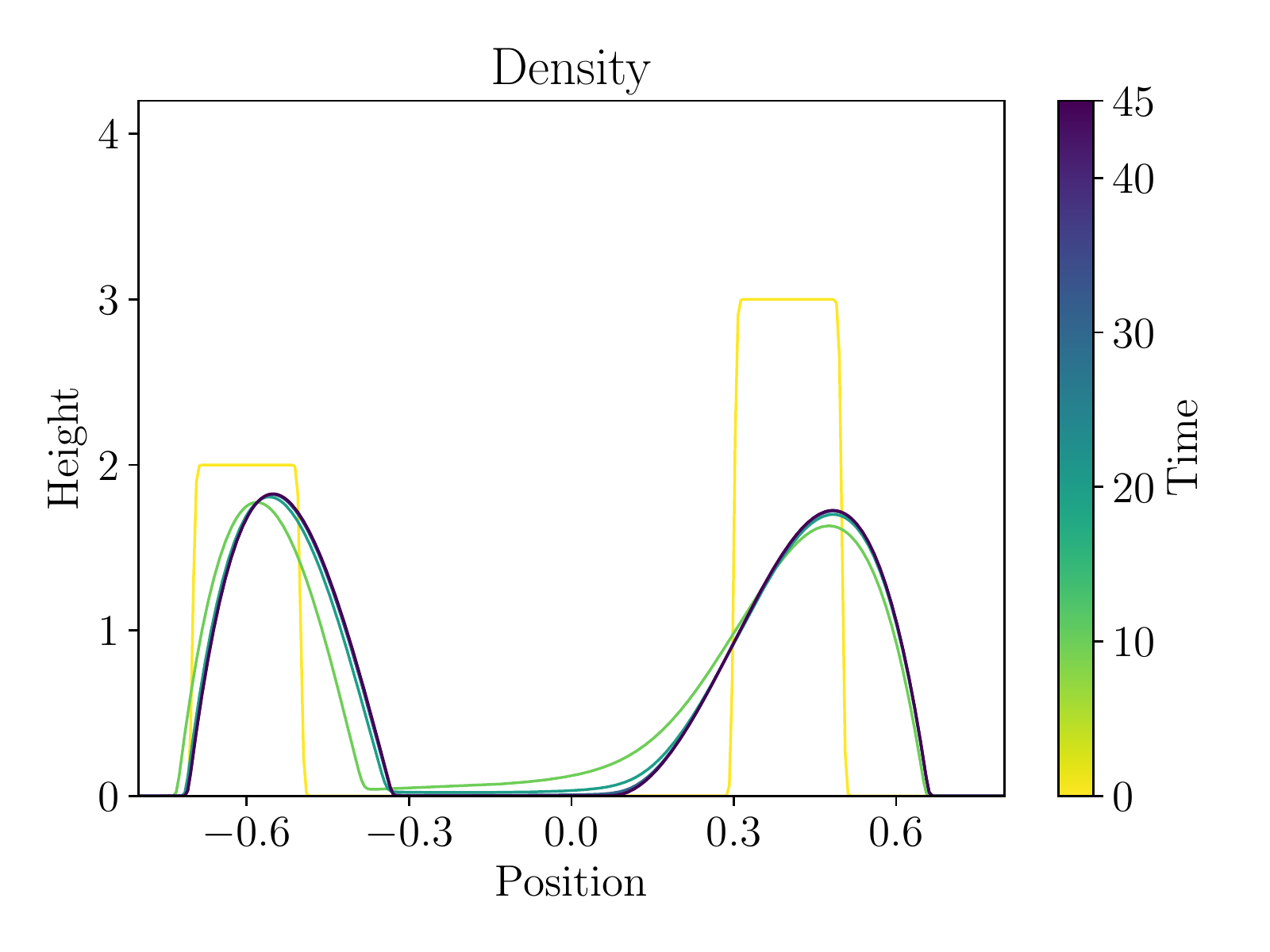} 
 
   {\bf m=2, $ \bf \text{\boldmath$\nu$} = 0.01$} \\
 \includegraphics[height=4cm,trim={.6cm .7cm .6cm .6cm},clip, valign=t]{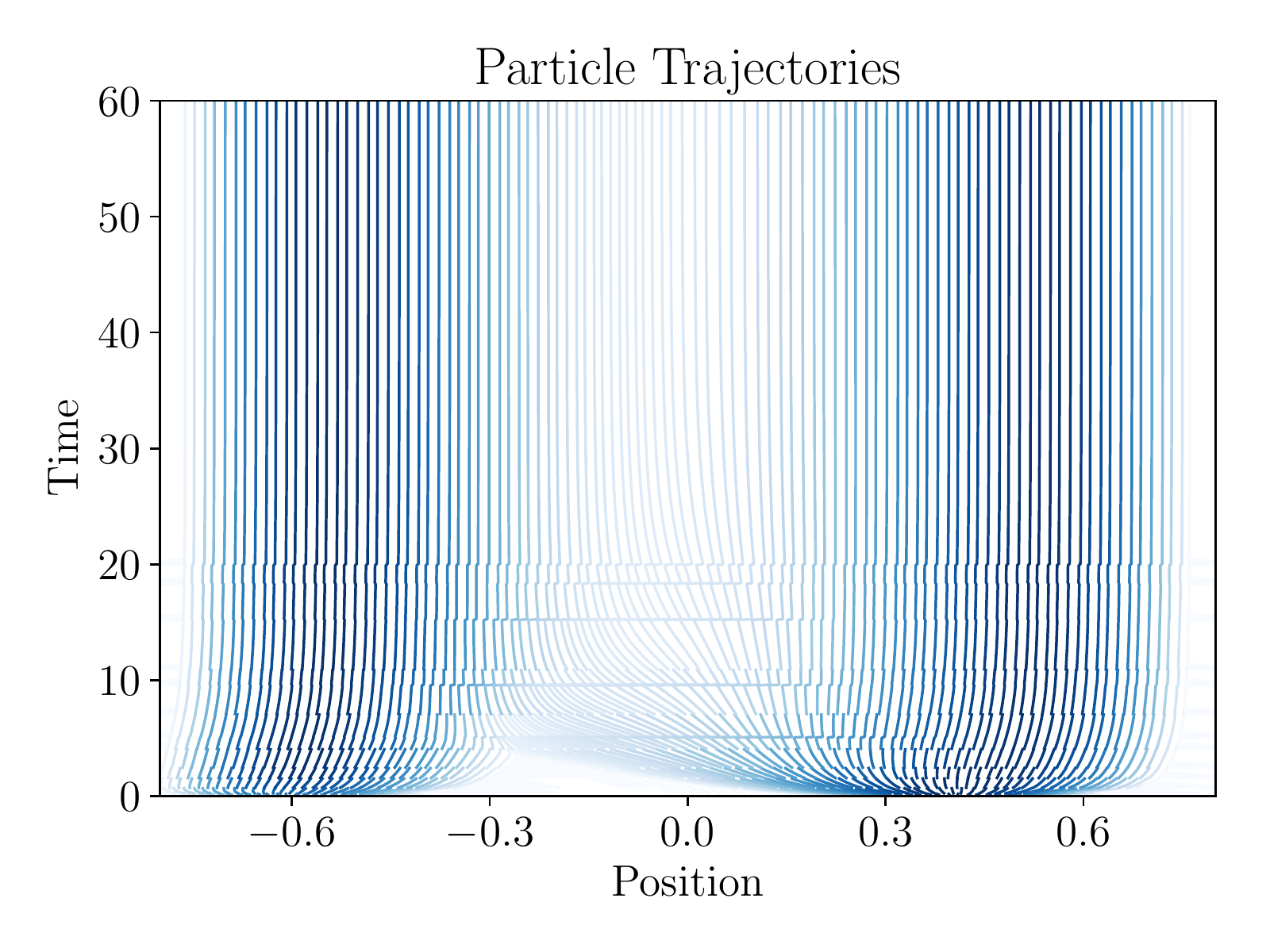}  \includegraphics[height=4cm,trim={.6cm .7cm .6cm .6cm},clip, valign=t]{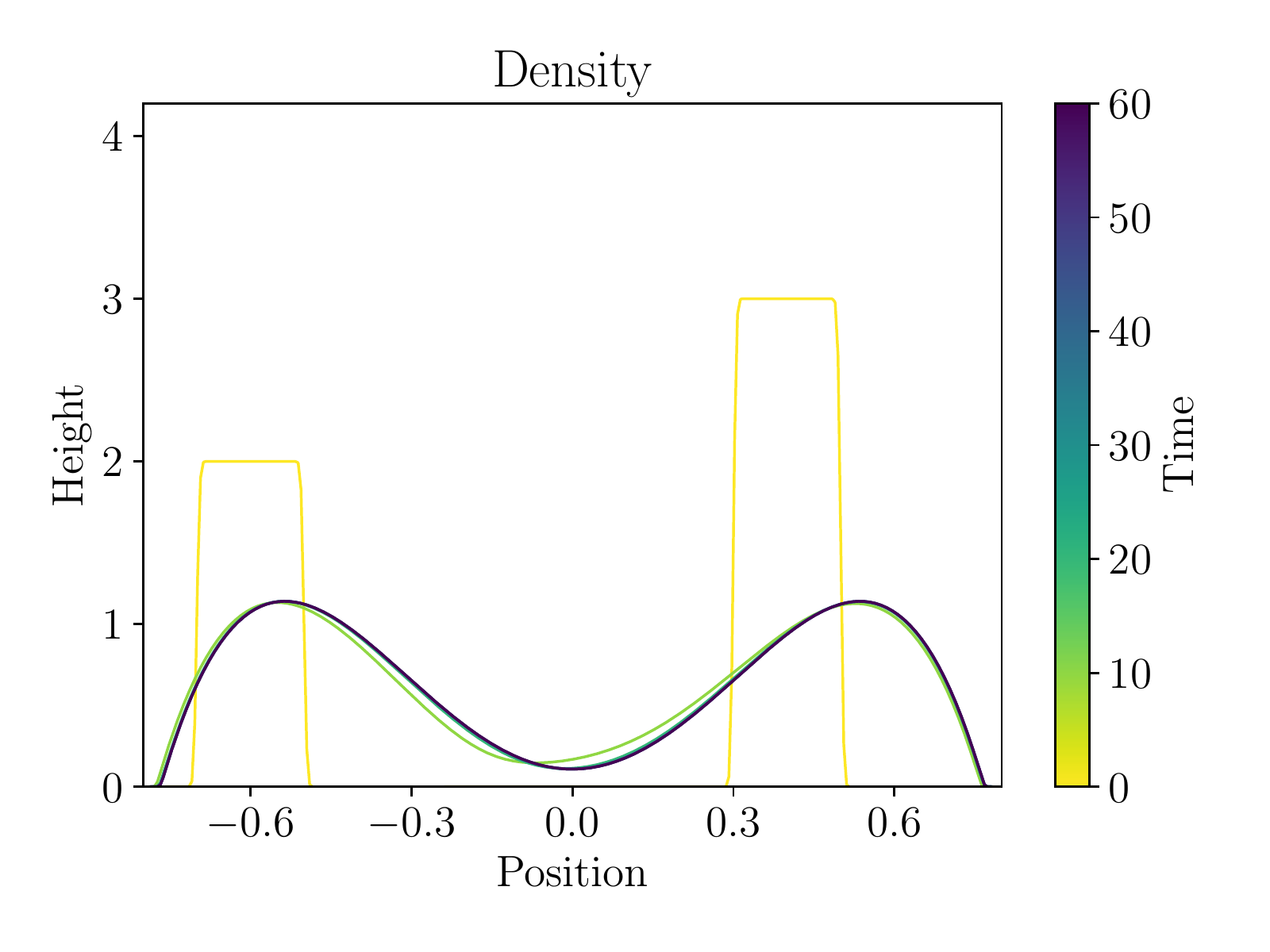} 
 
   {\bf m=800, $ \bf \text{\boldmath$\nu$} = 1$} \\
 \includegraphics[height=4cm,trim={.6cm .7cm .6cm .6cm},clip, valign=t]{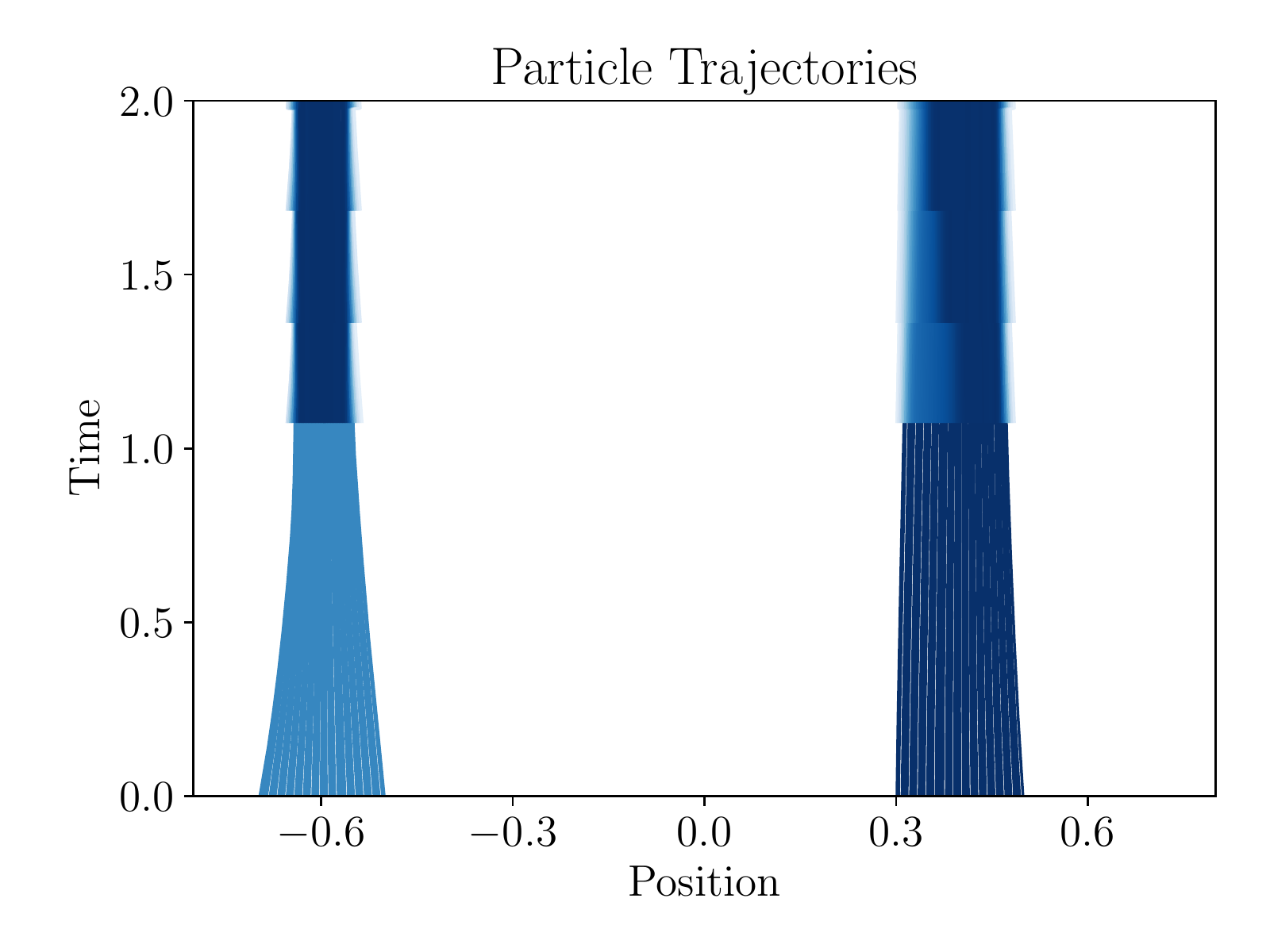} \includegraphics[height=4cm,trim={.6cm .7cm .6cm .6cm},clip, valign=t]{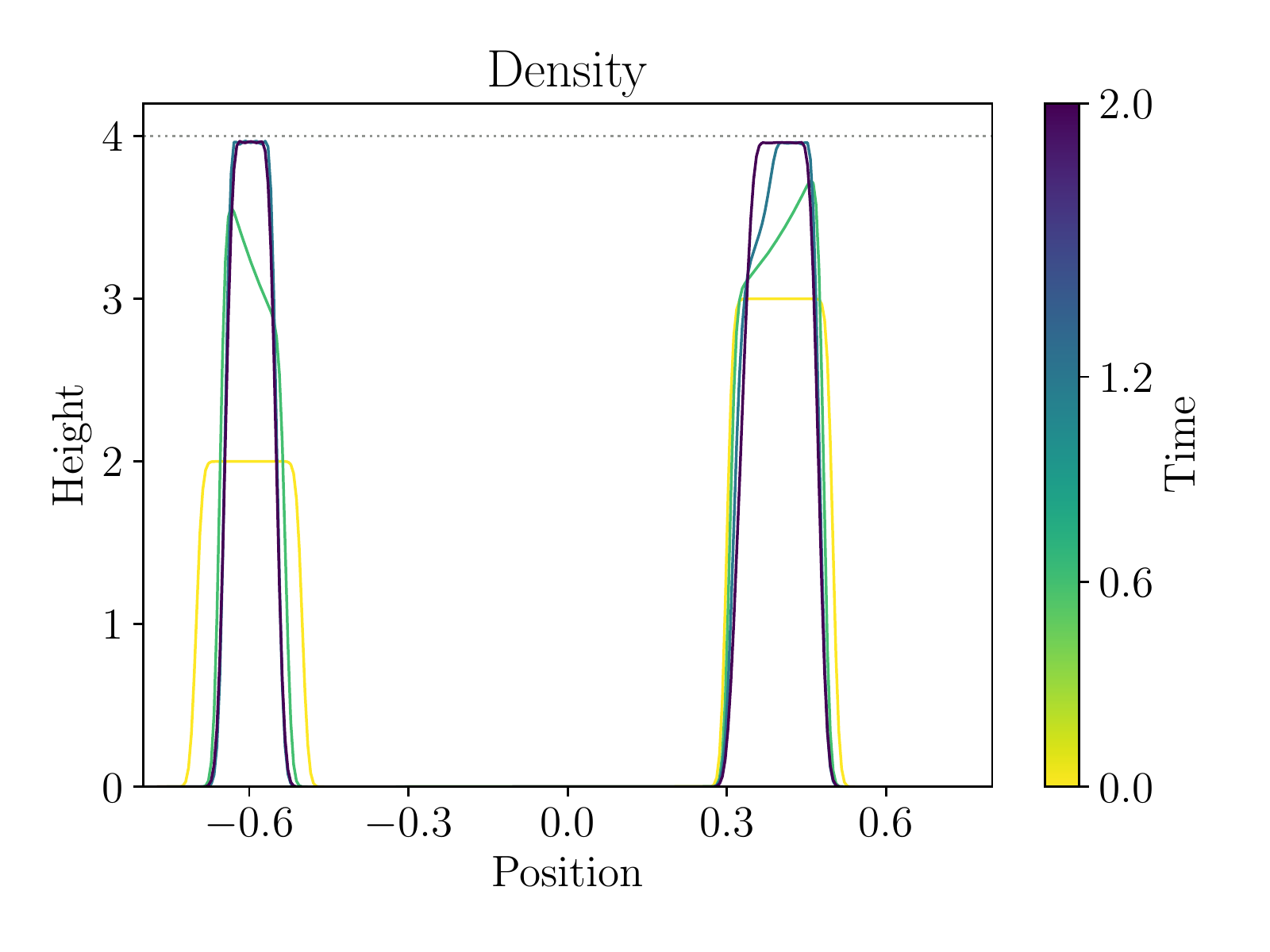}

\caption{For a repulsive-attractive interaction kernel and linear diffusion ($m=1$), solutions quickly evolve to a metastable state of two bumps with unequal mass and then the mass between the two bumps slowly equilibrates over time; see \cite{EK}. For quadratic diffusion ($m=2$), the solution   quickly forms two unequal bumps, which do not symmetrize for small diffusion coefficient ($\nu = 0.075^2/2$), but do symmetrize for larger diffusion coefficient ($\nu = 0.01$). For  the height constrained problem, the solution does not equilibrate, as there is  no diffusion to mediate this process. }
\label{metax4mx2}
\end{center}

\end{figure}

In our final example,  Figure \ref{metax4mx2}, we illustrate metastability behavior for aggregation diffusion exponents with the repulsive attractive interaction kernel $W(x) = |x|^4-|x|^2$ and various types of diffusion. For diffusion exponent $m=1$, weighted with diffusion coefficient $\nu = 0.075^2/2$, we recover the metastability behavior observed by Evers and Kolokolnikov \cite{EK}: a solution that is initially given by characteristic functions with weight $0.2$ centered at $-0.6$ and weight $0.3$ centered at $0.4$ quickly smooths to two bumps of unequal weights and then equilibrates slowly over time. We take $N= 600$ particles, with regularization $\epsilon = 0.99$, and maximal time step $k = 10^{-3}$. Interestingly, the equilibration behavior appears to be driven entirely by particles with extremely small mass on the tails of the bumps and is therefore nearly indistinguishable at the particle level.


Next, in Figure \ref{metax4mx2}, we consider  diffusion exponent $m=2$, weighted with  diffusion coefficient $\nu = 0.075^2/2$ and observe that the solution starting from characteristic functions with unequal weights does not symmetrize: the larger mass on the right hand side is preserved asymptotically. However, when we increase the  diffusion exponent to $\nu = 0.01$ in the third simulation, the solution symmetrizes  quickly. In both cases, we take $N = 800$ particles, with regularization $\epsilon = 0.99$, and maximal time step $10^{-3}$.

In the fourth simulation of Figure \ref{metax4mx2}, we contrast  the previous examples of degenerate diffusion $m=1,2$   with the height constrained aggregation equation. Naturally, without the mechanism of diffusion to spread mass, the unequal distribution of mass in the initial data is preserved asymptotically. We discretize the interval $[-0.9,0.7]$ with $N = 600$ particles, with regularization $\epsilon = 0.85$, and maximal times step $k=10^{-4}$.

\section*{Acknowledgements}
\vspace{-.3cm}
We thank Matias Delgadino, Franca Hoffmann, Jingwei Hu, Francesco Patacchini, Ihsan Topaloglu, and Li Wang  for useful discussions. JAC was partially supported by the EPSRC grant number  EP/P031587/1. JAC is grateful to the Mittag-Leffler Institute for providing a fruitful working environment during the special semester \emph{Mathematical Biology.} KC was supported by NSF DMS-1811012. YY was supported by NSF DMS-1715418. The authors acknowledge the American Institute of Mathematics (AIM)  for supporting a visit during the early stages of this work. This work used the Extreme Science and Engineering Discovery Environment (XSEDE) Comet at the San Diego Supercomputer Center through allocation DMS180023.

\bibliographystyle{abbrv}
\bibliography{AggDiffBookChapterBiblioFile}
\end{document}